\documentclass[11pt,letterpaper]{amsart}

\usepackage{epic,eepic,latexsym, amssymb, amscd, amsfonts, xypic, floatflt,ifthen, graphicx}

\usepackage{psfrag,xspace}
\usepackage{fancybox}
\usepackage{psboxit}
\newcommand{\whichos}{"Mac"}
\ifthenelse{\equal{\whichos}{"PC"}}
{\usepackage{graphics}}
{\usepackage{graphicx}}

\usepackage{tabularx}

\newcommand{\puteps}[2][0.5]
{\ifthenelse{\equal{\whichos}{"PC"}}
{\scalebox{#1}{%
\includegraphics{ #2.eps}}}%
{\includegraphics[scale=#1]{#2.eps}}}

\input xy
\xyoption{all}

 \newlength{\baseunit}               
 \newcount{\numlines}                
\newcommand{\bpf}{\noindent {\em Proof.  }}
\newcommand{\epf}{\qed \vspace{+10pt}}

\newtheorem{tm}{Theorem}[section]
\newtheorem{pr}[tm]{Proposition}
\newtheorem{lm}[tm]{Lemma}
\newtheorem{co}[tm]{Corollary}
\newtheorem{con}[tm]{Conjecture}
\newtheorem{df}[tm]{Definition}
\newtheorem{re}[tm]{Remark}

\newcommand{\Faber}{\mathsf{Fab}}

\newcommand{\E}{\mathbb{E}}
\newcommand{\F}{\mathbb{F}}
\renewcommand{\H}{\mathbb{H}}
\newcommand{\Hyp}{\mathbb{H}yp}
\renewcommand{\P}{\mathbb{P}}

\newcommand{\Q}{\mathbb{Q}}

\newcommand{\C}{\mathbb{C}}
\newcommand{\G}{\mathbb{G}}
\renewcommand{\zeta}{\mathbb{T}}

\newcommand{\ka}{\kappa}
\newcommand{\proj}{\mathbb P}

\newcommand{\oh}{{\mathcal{O}}}

\newcommand{\cm}{{\mathcal{M}}}
\newcommand{\cmbar}{\overline{\cm}}

\newcommand{\al}{\alpha}

\newcommand{\be}{\beta}
\newcommand{\ga}{\gamma}

\newcommand{\de}{\delta}
\newcommand{\De}{\Delta}
\newcommand{\ep}{\epsilon}
\newcommand{\si}{\sigma}
\newcommand{\la}{\lambda}

\newcommand{\Sym}{\operatorname{Sym}}

\newcommand{\virt}{{\rm{vir}}}
\newcommand{\wt}{{\rm{wt}}}
\newcommand{\ts}{\mathrm}
\newcommand{\f}{\frac}
\newcommand{\La}{{\Lambda}}
\newcommand{\Om}{{\Omega}}
\newcommand{\vep}{{\varepsilon}}

\newcommand{\Ups}{{\mathsf{T}}}
\newcommand{\ld}{{\ldots}}
\newcommand{\cd}{{\cdots}}

\newcommand{\sV}{{\mathcal{V}}}
\newcommand{\sT}{{\mathcal{T}}}
\newcommand{\sE}{{\mathcal{E}}}
\newcommand{\sA}{{\mathcal{A}}}
\newcommand{\sB}{{\mathcal{B}}}
\newcommand{\cS}{{\mathcal{S}}}
\newcommand{\oE}{{\mathsf{E}}}
\newcommand{\sC}{{\mathsf{C}}}
\newcommand{\sS}{{\mathsf{S}}}

\newcommand{\SYM}{{\mathsf{sym}}}
\newcommand{\sP}{{\mathcal{P}}}

\newcommand{\bfp}{{\bf p}}
\newcommand{\bfq}{{\bf q}}
\newcommand{\bfH}{{\bf H}}

\newcommand{\pa}{{\partial}}

\newcommand{\widH}{{\widehat{H}}}
\newcommand{\widbfH}{\widehat{\bf H}}

\newcommand{\whY}{{\widehat{Y}}}
\newcommand{\raFab}{{\rangle_g^{\F \mathrm{ab}}}}
\newcommand{\laFab}{{\langle}}

\newcommand{\cV}{{\mathcal{V}}}
\newcommand{\bij}{\stackrel{\sim}{\longrightarrow}}
\newcommand{\rt}{\mathfrak{t}}
\newcommand{\rw}{\mathfrak{w}}
\newcommand{\rv}{\mathfrak{v}}

\newcommand{\cT}{{\mathcal{T}}}

\newcommand{\cU}{{\mathcal{U}}}
\newcommand{\ombr}{\Omega^{\mathsf{br}}}

\newcommand{\bfM}{{\mathbf{M}}}
\newcommand{\bfy}{{\mathbf{y}}}

\newcommand{\cited}{}

\newcommand{\lremind}[1]{{\bf[label:  #1]}}
\newcommand{\notation}[1]{}
\renewcommand{\lremind}[1]{{}}

\newcommand{\secretnote}[1]{}
\newcommand{\cut}[1]{}

\newcommand{\dmjdel}[1]{}  

\begin{document}
\pagestyle{plain} \title{{\large {
The moduli space of curves, double Hurwitz numbers,
and Faber's intersection number conjecture
}}
}
\author{I. P. Goulden, D. M. Jackson and R. Vakil}
\address{Department of Combinatorics and
Optimization, University of Waterloo
}
\address{Department of Combinatorics and
Optimization, University of Waterloo
}
\address{Department of Mathematics, Stanford University
}

\thanks{The first two authors are partially supported by
NSERC grants.
The third author is partially supported by NSF PECASE/CAREER grant DMS--0238532.
}
\date{Tuesday, November 21, 2006.}

\begin{abstract}
  We define the dimension $2g-1$ Faber-Hurwitz Chow/homology classes
  on the moduli space of curves, parametrizing curves expressible as
  branched covers of $\proj^1$ with given ramification over $\infty$
  and sufficiently many fixed ramification points elsewhere.
  Degeneration of the target and judicious localization expresses such
  classes in terms of localization trees weighted by ``top
  intersections'' of tautological classes and genus $0$ double Hurwitz
  numbers.  This identity of generating series can be inverted,
  yielding a ``combinatorialization'' of top intersections of
  $\psi$-classes.  As genus $0$ double Hurwitz numbers with at most
  $3$ parts over $\infty$ are well understood, we obtain Faber's
 Intersection Number Conjecture for up to $3$ parts, and an approach
  to the Conjecture in general (bypassing the Virasoro Conjecture).
  We also recover other geometric results in a unified manner,
  including Looijenga's theorem, the socle theorem for curves with
  rational tails, and the hyperelliptic locus in terms of $\ka_{g-2}.$
\end{abstract}

\maketitle


\part{INTRODUCTION AND SUMMARY OF RESULTS}\label{partA}

Since we shall be using arguments from geometry and combinatorics,
we have separated the material into three parts to
assist the reader.  Part~\ref{partA} gives the background to the topic
and a summary of our results.  Part~\ref{partB} contains the geometry
that uses degeneration to obtain a recursion for the Faber-Hurwitz classes,
and localization to express these as tree sums involving the Faber symbol.
Part~\ref{partC} contains an approach through
algebraic combinatorics  to  transform and then  solve the formal partial
differential equations and functional equations that originate from
degeneration and localization in Part~\ref{partB} and thence to obtain the top
intersection numbers.  We have sought to
make the transition from the geometry of Part~\ref{partB} to the
combinatorics of Part~\ref{partC} pellucid.

\section{Summary of results}
\label{sor}\lremind{sor}The purpose of this paper is to give a geometrico-combinatorial
approach that is  direct and, we hope, enlightening,
 to the three known results that are listed below. 
 We give a summary of results for those quite familiar with moduli
spaces of curves, and Faber's foundational  conjectures on their
cohomology or Chow rings.  
A more detailed introduction  to the paper is given in
Section \ref{introduction}, and most readers should turn immediately to this.

The three results are:

\begin{itemize}
\item[\textbf{(I)}]
$R_{2g-1}(\cm_{g,n}^{rt})$ is generated by a single element, which we 
denote $\G_{g,1}$
(Theorem~\ref{onedim}).
This argument was promised in \cite[Sec.~5.7]{thmstar}.
(Since $R_{2g-1}(\cm_g)$ was shown earlier by Faber
to be non-zero \cite[Thm.~2]{fconj},   this single element
is also non-zero by Remark~\ref{summary}(iii) below.)

\item [\textbf{(II)}]
A combinatorial description of $\psi_1^{a_1} \cdots
\psi_n^{a_n} \in R_{2g-1}(\cm_{g,n}^{rt})$ as a multiple of this
generator, in terms of genus $0$ double Hurwitz numbers
(Theorem~\ref{bigyuck}).  (These
intersections of $\psi$-classes determine all top intersections in the
tautological ring, and are the subject of Faber's Intersection Number
Conjecture.)

\item [\textbf{(III)}]
Hence a proof of Faber's
Intersection Number Conjecture for up to three points,
and arbitrary genus (Theorem~\ref{whatwedo}).
\end{itemize}

Past proofs of some of these results are described in Section \ref{background}.
 In the above statements we have used the following notation.
Let $\cm^{rt}_{g,n}$ be the moduli
space of $n$-pointed genus $g$ stable curves with ``rational tails'',
and $R^*(\cm^{rt}_{g,n})$  its tautological ring.
Genus $0$ double Hurwitz numbers enumerate
branched covers of the sphere
by another sphere, with branching over $0$ and $\infty$ specified by
partitions $\al$ and $\be$ respectively, and the simplest non-trivial branching 
over an appropriate number of other given points.  
They are well understood
in the case where one of the partitions has at most
three parts \cite{gjv2}.

The ``one-part analysis'' yields as corollaries new proofs of a number
of important facts, such as the class of the hyperelliptic locus in
$\cm_g$ as a multiple of $\ka_{g-2}$.  These results are collected as
Corollaries~\ref{cgformula} through \ref{Fnice}.  A direct proof of
these consequences by our methods would be much shorter;
much of our effort will be to develop techniques to deal with more
parts.

We note that Faber's Intersection Number Conjecture for {\em bounded
  genus} involves a finite amount of information, thanks to the string
and dilaton equation (Prop.~ \ref{sd}); for any (reasonably small)
genus, this finite amount of information can be directly computed
\cite{fprog}.   (We point out that although Faber's verification of his
conjecture up to genus $21$, \cite{femail}, using work of
Pandharipande, involves a finite amount of information, it is a 
large amount, and is very difficult.)
However, Faber's Intersection Number Conjecture for a
{\em bounded number of points}, the case considered here, involves an infinite amount of
information.  

The methods here can readily be adapted to deal with a larger number
of points.  For example, the case $n=4$ follows from the formulae for
genus $0$ double Hurwitz numbers from \cite{gjv2}, and a
\textsf{Maple} computation.  The case $n=5$ should also be
computationally tractable.  However, we content ourselves with
calculations that could be done by hand, as the route to a natural
proof of Faber's Intersection Number Conjecture is clearly not through
any hoped-for closed form description of the genus $0$ double Hurwitz
generating series in general --- such series might be expected to get
quite complicated.  Instead, we hope for a proof using the structure
of the double Hurwitz generating series as a whole, and there are
indications that this may be tractable.  We point out in particular
the very recent preprint \cite{ssv}, giving a good description of
genus $0$ double Hurwitz numbers in general, and a particularly
elegant description in many cases.

\subsection{Motivation}
The ELSV formula provides a remarkable link between Hurwitz numbers
(counting branched covers of $\proj^1$ or, combinatorially, transitive
factorizations of elements of the symmetric group into transpositions)
and  the intersection theory on the moduli space of curves.  
See
\cite{elsv1, elsv2}, and \cite{gvelsv} for a proof in the context of
Gromov-Witten theory.    The ELSV formula describes
Hurwitz numbers in terms of top intersections of the moduli space
of curves.  This relation can be ``inverted'' to prove
results on top intersections on the moduli space of curves.
This has been used by a large number of authors, including
Ekedahl, Kazarian, Lando, Okounkov, Pandharipande, Shadrin,
Shapiro, Vainshtein and Zvonkine, to great effect.
Notable examples are the proofs of 
Witten's Conjecture \cite{op, kl}.

This paper relies on the observation that similar methods can be
applied to the {\em (non-compact)} moduli space of {\em smooth}
curves.  In this case, localization turns question about top
intersections into ``genus $0$ combinatorics'' (involving trees rather
than general graphs), which should in principle be simpler than the
genus $g$ combinatorics that arises in the case of compactified moduli space.
Indeed, we immediately get some results (for example, results {\bf
  (I)} and {\bf (II)}) from the statement of the analogue of the ELSV
formula.  Explicitly inverting this ELSV-analogue requires more work,
but then our knowledge of genus $0$ double Hurwitz numbers with few
parts translates directly into Faber's Intersection Number Conjecture
with few parts (result {\bf (III)}).

We point out that these techniques of algebraic combinatorics are
useful in geometry in a wider context.  A further example of their use
appears in \cite{gjvlambdag} where we give a direct proof of Getzler
and Pandharipande's
$\lambda_g$-Conjecture (without Gromov-Witten theory), and related
techniques have been used by a number of the authors mentioned above.

\subsection{Background}
\label{background}\lremind{background}The one-dimensionality of {\bf (I)}
was established in \cite{fabp}, and may also follow
from Looijenga \cite{looijenga}. 
The argument given  here  can be seen as an extension of the tautological 
vanishing
theorem of \cite{thmstar}.  The
entire argument is outlined in a few pages in \cite[Sec.~7]{cime}.

Getzler and Pandharipande showed that Faber's Intersection Number
Conjecture is a formal consequence of the Virasoro Conjecture for the
projective plane, \cite{gp}, in fact just the degree $0$ part (the
``large volume'' limit).  Givental proved the Virasoro Conjecture for
projective space (and more generally Fano toric manifolds) \cite{gi1,
givental}.  Y.-P. Lee and Pandharipande are writing a book \cite{lp}
giving details.  Givental's result is one of the most important
results in Gromov-Witten theory, and is a marvelous feat.  However, it
seems circuitous to prove the Intersection Number Conjecture by means
of the Virasoro Conjecture. The latter is a very heavy instrument
which conceals the combinatorial structure that lies behind the
intersection numbers.  As noted by K.~Liu and Xu \cite{kefeng}, it is
very desirable to have a shorter and direct explanation.  (Liu and Xu
show how the conjecture cleanly follows from another attractive
conjectural identity.)  For this reason, we give such an argument,
paralleling our understanding of top intersection numbers on
$\cmbar_{g,n}$ \textit{via} Hurwitz numbers.

\subsection{Outline of paper}
The strategy of the paper is as follows.  We define dimension $2g-1$
Faber-Hurwitz Chow/homology classes on the moduli space of curves, by
considering (``virtually'') branched covers of $\proj^1$ with given
ramification over $\infty$ and sufficiently many fixed ramification
points elsewhere.  We consider only curves with rational tails, which
simplifies the combinatorics dramatically.  We use the two most
effective techniques of Gromov-Witten theory, degeneration and
localization.  

$\bullet$ Degeneration of the target
yields a recursion for such classes, which we can solve explicitly.
The story is analogous to that of genus $0$ Hurwitz numbers (which are
combinatorially straightforward), {\em not} genus $g$ Hurwitz numbers.

$\bullet$ Localization expresses such classes in
terms of the desired top intersections and genus $0$ double Hurwitz
numbers.  The rational tails constraint forces the resulting
localization graphs to be trees. 

More precisely, localization expresses Faber-Hurwitz classes as a sum
over certain decorated trees (Sec.~\ref{Dloctr}) of linear
combinations of top intersections (including those that are the
subject of Faber's conjecture) and genus $0$ double Hurwitz numbers.
The relation between the top intersections and Faber-Hurwitz classes
can be easily seen to be invertible, \textit{i.e.}\ the top intersections are
linear combinations of Faber-Hurwitz classes (which we already
understood \textit{via} degeneration).  From this, Looijenga's
Theorem  drops out quickly
(from Theorem~\ref{onedim}), for example. 
The central idea
of this paper
is that this inversion may be done explicitly: the
localization sum gives an expression which readily can be unwound,
and from which the top intersections can be extracted.
This is done by a change of variables.

Using this approach, we quickly recover various geometric results in a
unified manner (Sec.~\ref{sigc}).  Also, as genus $0$ double Hurwitz
numbers with at most $3$ parts over $\infty$ are well understood, we
obtain Faber's Intersection Number Conjecture for up to $3$ parts, and
an approach to the conjecture in general.

In Section~\ref{introduction} we state the Faber Intersection
Number Conjecture and give the geometric and algebraic combinatorics
background to our approach. The Faber-Hurwitz classes are defined in
Section~\ref{DandL}, and we obtain
a Join-cut Recursion for them by
\emph{degeneration} of the target. This is the first form of the Degeneration
Theorem, and the equivalent Join-cut (partial differential) Equation 
for the Faber-Hurwitz series is given as the second form of this Theorem.
In addition, using \emph{localization}, we obtain an
expression for the Faber-Hurwitz classes as a weighted sum over
(localization) trees, giving the first form of the Localization
Tree  Theorem.
Section~\ref{s:genseries} addresses the weighted tree sum,
using the combinatorics of rooted, labelled trees and
exponential generating series in a countable set of indeterminates.
This results in the second form of the Localization Tree Theorem,
which is purely algebraic, giving an expression for the
Faber-Hurwitz series in terms of Faber's intersection numbers 
and the unique solution to a functional equation.
Section~\ref{sssSOD} introduces the \emph{fundamental transformation},
as the composition of three operators. The first two operators are a
symmetrization and an implicit change of variables, which gives
\emph{polynomials} in the new indeterminates. The third operator restricts
these polynomials to ``top'' terms --- those of maximum total degree.
Our strategy of proof for Faber's Intersection Number Conjecture is to
apply the fundamental transformation to the Localization Tree Theorem
and to the Degeneration Theorem, and eliminate the transformed
Faber-Hurwitz series. The key to the proof is that
only top intersection numbers remain when we do so,
and they  appear in enough linearly independent equations to
uniquely determine them.
Section~\ref{Sinop} applies the  strategy
for the first time, to the one-part case (of Faber top intersection
numbers) which has some 
non-trivial geometric consequences.
The methodology to that point requires us to consider an equation
for each genus $g$ separately, so in 
Section~\ref{Safhst} we refine the methodology by creating generating series in genus.
This means
that a single generating series equation suffices to prove
Faber's Intersection Number Conjecture for each number of parts.
In Section~\ref{sFTINC23Pm}, we apply this refined methodology to establish Faber's (top)
Intersection Number Conjecture for $2$ and $3$ parts, and we include
remarks about the general case.

Appendix~\ref{glossary} is a glossary of notation for the reader's convenience.

\noindent {\bf Acknowledgments.}
The third author has benefited from discussions with Renzo Cavalieri,
Carel Faber, Y.-P. Lee, Rahul Pandharipande, Hsian-Hua Tseng, and
especially Tom Graber, who developed many of the algebro-geometric
foundations on which this paper is based.  We also thank
Ezra Getzler, Kefeng Liu, and Melissa Liu.

\section{Introduction}\label{introduction}

\subsection{Geometric background}\label{ssTR1}
We work over the complex numbers.  Throughout, the genus $g$ is at
least $1$.  We assume some knowledge of the moduli space of curves.
An overview is given in \cite{cime}, which outlines the necessary
background and ends with a sketch of many of the results of this
paper; see also \cite{notices}.  We also assume familiarity with
Gromov-Witten theory, in particular the theory of relative stable
maps, and with localization on their moduli space (``relative virtual
localization'') \cite{li1, li2, thmstar, llz}.  An introduction to
many of the Gromov-Witten ideas we shall use may be found in
\cite{cime} and \cite{mirsym}.  We prefer to work in the Chow ring
$A^*$ rather than the cohomology ring $H^{2*}$, because our arguments
apply in this more refined setting, but there is no loss should the
reader wish to work in cohomology.  Chow/homology classes will often
be written in blackboard bold font (\emph{e.g.}\ $\F$) in order to
distinguish them from numbers.

Faber's conjectures on the topology of the moduli space of smooth
curves $\cm_g$ (given in \cite{fconj}) are a striking description of
the ``tautological'' part of the cohomology ring, the part of the
cohomology ring arising ``naturally from geometry''.  The ``top
intersections'' in this ring have a particularly remarkable
combinatorial structure.  We begin by describing these conjectures.

On the moduli space of stable $n$-pointed genus $g$ stable curves
$\cmbar_{g,n}$ (or any open subset thereof), let $\psi_i$ ($1 \leq i
\leq n$) be the first Chern class of the line bundle corresponding to
the cotangent space of the universal curve at the $i$th marked point.

We shall denote all forgetful morphisms $\cmbar_{g,n} \rightarrow
\cmbar_{g,n'}$ by $\pi$, where the source and target will be clear
from the context.  For example if $\pi: \cm_{g,1} \rightarrow \cm_g$,
then the $i$th Mumford-Morita-Miller ``$\ka$-class'' is defined by
$\ka_i := \pi_* \psi_1^{i+1}$.

Given an $(n_1+1)$-pointed curve of genus $g_1$, and
an $(n_2+1)$-pointed curve of genus $g_2$, gluing
the first curve to the second along the last point of each
yields an $(n_1+n_2)$-pointed curve of genus $g_1+g_2$.
This gives a map
$$
\cmbar_{g_1, n_1+1} \times \cmbar_{g_2, n_2+1} \rightarrow \cmbar_{g_1+g_2, n_1+n_2}.$$
Similarly, we can take a single $(n+2)$-pointed curve of genus $g$,
and glue its last two points together to get an $n$-pointed curve
of genus $g+1$.  This gives a map
$$
\cmbar_{g,n+2} \rightarrow \cmbar_{g+1, n}.$$
We call these two types of maps {\em gluing morphisms}.
We call the forgetful and gluing morphisms the {\em natural morphisms}
between moduli spaces of curves.

\subsection{The tautological ring}\label{ssTR}
There are many equivalent definitions of the tautological ring in the literature.
The following will be convenient for our purposes.

\begin{df}\cite[Def.\ 4.2]{thmstar} \label{tautdef}\lremind{tautdef}
The system of tautological rings $( R^*(\cmbar_{g,n}) \subset
A^*(\cmbar_{g,n}))_{g,n}$ is the smallest system of $\Q$-vector spaces
closed under pushforwards by the natural morphisms, such that all
monomials in $\psi_1$, \dots, $\psi_n$ lie in $R^*(\cmbar_{g,n})$.
\end{df}

If $\cm$ is an open subset
of $\cmbar_{g,n}$, let 
$R^*(\cm) := R^*(\cm_{g,n})|_{\cm}$,
$\;R_j(\cm)  := R^{\dim \cm - j}(\cm)$.  (Here $R^k(\cm)$ is
the codimension $k$ part of $R^*(\cm)$.)
We shall often use the notation $R_*$ instead of $R^*$, because we wish
to think of classes as homology classes.  

We take this opportunity to introduce a third sort of tautological
class: let $\E_{g,n}$ be the {\em Hodge bundle} on $\cmbar_{g,n}$.  It
has rank $g$, and \lremind{lambdapullback}
\begin{equation*}\label{lambdapullback}
\pi^* \E_{g,0} = \E_{g,n}
\end{equation*}  
where $\pi$ is the forgetful morphism $\pi: \cmbar_{g,n} \rightarrow \cmbar_g$.
Over a point $[(C, p_1, \dots, p_n)] \in
\cm_{g,n}$, the fiber of $\E_{g,n}$ is the vector space of differentials on $C$.  The
$\la$-classes are defined by $\la_i = c_i(\E_{g,n})$.  
By  the above relation, they behave well with respect to
pullback by forgetful morphisms (``$\pi^* \la_i = \la_i$''). 
Note that $\la_0=1$.

It is not hard to show that the tautological ring of the moduli
space of {\em smooth curves} $\cm_g$ is
generated by the $\ka$-classes, and indeed this is essentially the original
definition of $R^*(\cm_g)$ in \cite{fconj}.

We now describe three predictions of Faber, namely
the Vanishing Conjecture, the Perfect Pairing Conjecture, and the Intersection Number
Conjecture, which is a central subject of this paper.

\noindent\underline{\em Vanishing Conjecture.}
$R^i(\cm_g) = 0$ for $i>g-2$, and $R^{g-2}(\cm_g) \cong \Q$.  This was
proved by Looijenga and Faber.  Looijenga's
Theorem \cite{looijenga} is that $R^{g-2}(\cm_g)$ 
is
generated as a vector space by a single element
(this will also follow from  our analysis, see
Thm.~\ref{onedim}).  Faber proved the following ``non-vanishing theorem''.

\begin{tm}[Faber \cite{fconj}, Thm.~2]
$R^{g-2}(\cm_g) \neq 0$. \label{nv}\lremind{nv}
\end{tm}

(For other proofs, see \cite[Thm.~6.5]{bp}, and \cite[Thm.~0.2]{bct},
which is based on \cite{renzo}.  More precisely, Faber described a
linear functional $\int \cdot \la_{g-1} \la_{g}: R^{g-2}(\cm_g)
\rightarrow \Q$, and showed that the image was non-zero by computing a
certain intersection number.  Bryan-Pandharipande~\cite{bp} and later
Bertram-Cavalieri-Todorov~\cite{bct} show the image is non-zero by different
enlightening computations.)

\noindent\underline{\em Perfect Pairing Conjecture.}
The analogue
of Poincar\'e duality holds: for $0 \leq i \leq g-2$, the
cup product $R^i(\cm_g) \times R^{g-2-i}(\cm_g) \rightarrow
R^{g-2}(\cm_g) \cong \Q$ is a perfect pairing.  This conjecture is 
known only in special cases, and is essentially completely open.

Informally, these two conjectures state that ``$R^*(\cm_g)$ behaves like
the ($(p,p)$-part of the) cohomology ring of a $(g-2)$-dimensional
complex projective manifold.''
They imply that the entire structure of the ring is
determined by the top intersections of the $\ka$-classes, \textit{i.e.}\ by
$\prod \ka_{d_i}$ (where $\sum_i d_i = g-2$) in terms of some fixed
generator of $R_{2g-1}(\cm_g)$.   

\noindent\underline{\em Faber's Intersection Number Conjecture.}
This gives a combinatorial
description of these top intersections. 
We note that this conjecture is useful even without knowing
the perfect pairing conjecture; Faber's algorithm \cite{falg}
reduces all ``top intersections'' in the tautological ring
to intersections of the form described in his
Intersection Number Conjecture.

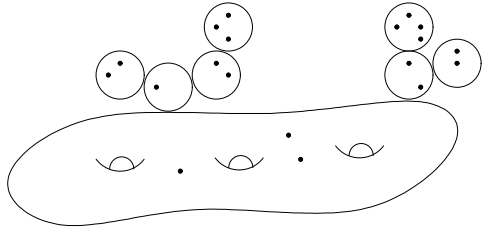
\begin{figure}[ht]
\begin{center}
\setlength{\unitlength}{0.00083333in}
\begingroup\makeatletter\ifx\SetFigFont\undefined%
\gdef\SetFigFont#1#2#3#4#5{%
  \reset@font\fontsize{#1}{#2pt}%
  \fontfamily{#3}\fontseries{#4}\fontshape{#5}%
  \selectfont}%
\fi\endgroup%
{\renewcommand{\dashlinestretch}{30}
\begin{picture}(2972,1422)(0,-10)
\path(577,687)(604,693)(631,697)
	(660,701)(689,705)(718,708)
	(749,710)(781,713)(813,715)
	(846,716)(880,717)(915,718)
	(950,718)(986,719)(1022,719)
	(1059,718)(1096,718)(1133,718)
	(1170,717)(1207,716)(1244,715)
	(1281,714)(1317,714)(1353,713)
	(1389,712)(1423,711)(1457,711)
	(1491,711)(1524,711)(1557,711)
	(1589,711)(1620,712)(1652,712)
	(1685,714)(1718,716)(1751,718)
	(1784,720)(1817,722)(1851,725)
	(1885,729)(1919,732)(1953,736)
	(1987,739)(2021,743)(2055,747)
	(2089,751)(2123,755)(2156,759)
	(2188,763)(2220,767)(2252,770)
	(2282,774)(2311,777)(2340,779)
	(2368,782)(2394,784)(2419,786)
	(2444,787)(2467,788)(2490,789)
	(2511,789)(2532,788)(2552,787)
	(2572,786)(2592,784)(2611,781)
	(2630,778)(2648,774)(2665,769)
	(2682,763)(2698,757)(2713,750)
	(2727,743)(2741,735)(2754,726)
	(2765,716)(2776,706)(2785,695)
	(2794,684)(2801,672)(2807,660)
	(2812,648)(2815,635)(2817,622)
	(2819,609)(2818,595)(2817,582)
	(2815,568)(2812,554)(2807,540)
	(2802,525)(2794,509)(2786,492)
	(2776,475)(2764,458)(2751,440)
	(2737,422)(2721,403)(2703,384)
	(2684,364)(2663,345)(2641,325)
	(2618,306)(2594,287)(2568,269)
	(2541,250)(2514,233)(2486,216)
	(2457,201)(2427,186)(2398,172)
	(2368,159)(2338,148)(2307,137)
	(2276,128)(2245,120)(2214,112)
	(2186,107)(2158,103)(2129,99)
	(2099,96)(2069,93)(2037,91)
	(2005,90)(1972,89)(1938,88)
	(1904,89)(1869,89)(1833,90)
	(1796,91)(1760,93)(1723,94)
	(1686,96)(1649,98)(1613,100)
	(1576,103)(1540,105)(1505,107)
	(1470,108)(1435,110)(1402,111)
	(1369,113)(1337,113)(1305,114)
	(1274,114)(1244,113)(1214,112)
	(1182,111)(1151,109)(1119,106)
	(1087,103)(1056,100)(1024,95)
	(992,91)(959,86)(927,81)
	(894,76)(862,70)(829,65)
	(797,59)(765,53)(734,48)
	(703,43)(673,37)(643,33)
	(614,28)(586,24)(558,21)
	(532,18)(506,15)(481,13)
	(457,12)(434,12)(411,12)
	(389,12)(363,14)(338,17)
	(313,22)(289,27)(264,33)
	(240,40)(217,49)(194,58)
	(172,69)(150,81)(130,93)
	(111,107)(94,121)(77,135)
	(63,151)(50,167)(39,183)
	(30,199)(23,216)(17,232)
	(14,249)(12,266)(12,283)
	(14,300)(17,315)(22,330)
	(28,345)(35,361)(44,377)
	(54,393)(66,409)(80,426)
	(95,442)(111,459)(129,476)
	(149,493)(169,509)(191,525)
	(215,541)(239,557)(264,572)
	(290,586)(317,600)(344,613)
	(372,625)(400,636)(428,647)
	(457,657)(487,666)(516,674)
	(546,681)(577,687)
\put(724.000,365.000){\arc{150.000}{3.1416}{6.2832}}
\put(2207.000,620.500){\arc{375.000}{0.6435}{2.4981}}
\put(2217.000,448.000){\arc{150.000}{3.1416}{6.2832}}
\put(1457.000,545.500){\arc{375.000}{0.6435}{2.4981}}
\put(1467.000,373.000){\arc{150.000}{3.1416}{6.2832}}
\put(2514,950){\ellipse{300}{300}}
\put(2514,1250){\ellipse{300}{300}}
\put(2814,1025){\ellipse{300}{300}}
\put(1014,875){\ellipse{300}{300}}
\put(1314,950){\ellipse{300}{300}}
\put(1389,1250){\ellipse{300}{300}}
\put(714,950){\ellipse{300}{300}}
\put(1089,350){\blacken\ellipse{24}{24}}
\put(1089,350){\ellipse{24}{24}}
\put(1764,575){\blacken\ellipse{24}{24}}
\put(1764,575){\ellipse{24}{24}}
\put(1839,425){\blacken\ellipse{24}{24}}
\put(1839,425){\ellipse{24}{24}}
\put(639,950){\blacken\ellipse{24}{24}}
\put(639,950){\ellipse{24}{24}}
\put(714,1025){\blacken\ellipse{24}{24}}
\put(714,1025){\ellipse{24}{24}}
\put(939,875){\blacken\ellipse{24}{24}}
\put(939,875){\ellipse{24}{24}}
\put(1314,1025){\blacken\ellipse{24}{24}}
\put(1314,1025){\ellipse{24}{24}}
\put(1389,950){\blacken\ellipse{24}{24}}
\put(1389,950){\ellipse{24}{24}}
\put(1389,1175){\blacken\ellipse{24}{24}}
\put(1389,1175){\ellipse{24}{24}}
\put(1389,1325){\blacken\ellipse{24}{24}}
\put(1389,1325){\ellipse{24}{24}}
\put(1314,1250){\blacken\ellipse{24}{24}}
\put(1314,1250){\ellipse{24}{24}}
\put(2439,1250){\blacken\ellipse{24}{24}}
\put(2439,1250){\ellipse{24}{24}}
\put(2514,1325){\blacken\ellipse{24}{24}}
\put(2514,1325){\ellipse{24}{24}}
\put(2589,1250){\blacken\ellipse{24}{24}}
\put(2589,1250){\ellipse{24}{24}}
\put(2589,1175){\blacken\ellipse{24}{24}}
\put(2589,1175){\ellipse{24}{24}}
\put(2514,1025){\blacken\ellipse{24}{24}}
\put(2514,1025){\ellipse{24}{24}}
\put(2589,875){\blacken\ellipse{24}{24}}
\put(2589,875){\ellipse{24}{24}}
\put(2814,1100){\blacken\ellipse{24}{24}}
\put(2814,1100){\ellipse{24}{24}}
\put(2814,1025){\blacken\ellipse{24}{24}}
\put(2814,1025){\ellipse{24}{24}}
\put(714.000,537.500){\arc{375.000}{0.6435}{2.4981}}
\end{picture}
}
\end{center}
\caption{A pointed curve with ``rational tails''.  Notice that the
dual graph is a tree.
\lremind{rt}}
\label{rt}
\end{figure}

Faber reformulated his
conjecture in the striking form given in Conjecture~\ref{ConjFinal}.  In order to give this
reformulation, we review the extension of Faber's Conjecture to curves
with ``rational tails'' (by Faber and Pandharipande, see for example
\cite{icm}).\label{rtdefhere}\lremind{rtdefhere} Recall that a nodal
genus $g$ nodal curve is said to be a {\em genus $g$ curve with
  rational tails} if one component is a smooth curve of genus $g$, and
hence the remaining components are genus $0$ (spheres), see
Figure~\ref{rt}.  Then the dual graph is a tree, a fact which will prove crucial
for us.  This is the reason
that the graphs arising from localization involves just tree
combinatorics.

Denote the corresponding open subset of $\cmbar_{g,n}$ 
(corresponding to stable curves
with rational tails) by $\cm_{g,n}^{rt}$.  If $\pi:
\cmbar_{g,n} \rightarrow \cmbar_g$, then $\cm_{g,n}^{rt} = \pi^{-1}
(\cm_g)$.  In particular, the forgetful morphism $\cm_{g,n}^{rt}
\rightarrow \cm_{g,n'}^{rt}$ for $n \geq n'$ is projective, and
$\cm_{g,1}^{rt} = \cm_{g,1}$.

We  recall some background.
\begin{re}\label{summary}
\begin{enumerate}
\item[i)] $R_j(\cm^{rt}_{g,n}) = 0$ for $j<2g-1$ and all $n$;
\item[ii)] for $n\geq 1$, \lremind{isoeq}
\begin{equation} \label{isoeq}
\pi_*: R_{2g-1}(\cm^{rt}_{g,n}) \rightarrow
R_{2g-1}(\cm^{rt}_{g,1}) = R_{2g-1}(\cm_{g,1})
\end{equation}
(where $\pi: \cm^{rt}_{g,n} \rightarrow \cm_{g,1}$ is the forgetful
morphism) is an isomorphism; and
\item[iii)] $\pi_*: R_{2g-1}(\cm_{g,1}) \rightarrow R_{2g-1}(\cm_g)$ 
(where $\pi: \cm_{g,1} \rightarrow \cm_g$ is the forgetful morphism)
is
surjective.
\end{enumerate}
\end{re}
Statements (i) and (ii) are
parts (a) and (b) of \cite[Prop.~5.8]{thmstar}.
Statement (iii) is immediate from an appropriate 
definition of the tautological ring, for example Definition~\ref{tautdef}.

\begin{con}[Faber's Intersection Number Conjecture, ``$\psi$-form'']\label{ConjFinal}
For any $n$-tuple of {\bf positive} integers $(d_1, \dots, d_n)$,
\begin{equation}\label{Fabconj}
\pi_* \left(\psi_1^{d_1} \cdots \psi_n^{d_n}\right) = 
\frac { (2g-3+n)! (2g-1)!!} { (2g-1)! \prod_{j=1}^n (2d_j-1)!!}
\psi_1^{g-1}\quad \text{for $\sum_i d_i = g-2+n$}
\end{equation}
where $\pi$ is the forgetful
morphism $\cm_{g,n}^{rt} \rightarrow \cm_{g,1}$.\label{finalform}\lremind{finalform}
(Recall that $(2a-1)!! = 1 \cdot 3  \cdots (2a-1) = (2a)! /
(2^a a!)$.)
\end{con}

Note that Faber's Intersection Number Conjecture is immediate for the
case $n=1$.  We shall prove the following.
\begin{tm}\label{Tfint}
Faber's Intersection Number
  Conjecture~\ref{finalform} is true for up to $3$ points (i.e.\ for
  $n \leq 3$).  \label{whatwedo}\lremind{whatwedo}
\end{tm}

We shall need to consider more general intersection numbers
involving one $\la$-class, arising in the virtual localization
formula.  Motivated by the isomorphism \eqref{isoeq}, and in analogy
with the ``Witten symbol'' used in Gromov-Witten theory (\textit{e.g.}
\cite[Sec.~26.2]{mirsym}), define the {\em Faber symbol} as
\begin{equation}\label{eWS}
\laFab \tau_{a_1} \cdots \tau_{a_n} \la_k \raFab :=
\pi_* \left( \psi_1^{a_1} \cdots \psi_n^{a_n} \la_k \right)
\in R_{2g-1}(\cm_{g,1})
\end{equation} 
where the pushforward $\pi_*$ is the isomorphism of \eqref{isoeq}
($\pi: \cm^{rt}_{g,n} \rightarrow \cm_{g,1}$ is the forgetful morphism).
\emph{Caveat}: this symbol is a Chow class, {\em not} a number; we repeat
that we will often indicate classes using blackboard bold font.  The symbol
is declared to be zero if some $a_i<0$, or if the product lies in
wrong degree, \textit{i.e.}\
$$k+\sum_i a_i  \neq \dim \cm_{g,n}^{rt} -(2g-1) = g-2+n.$$ 
Note that the case $k=0$ is equivalent to there being no $\la$-factor,
as $\la_0=1$. 
This symbol satisfies a version of the usual string and dilaton equations:

\begin{pr}
The following relations among the Faber symbols hold.\lremind{sd, string, dilaton}
\label{sd}
\begin{eqnarray*}
\laFab \tau_0 \tau_{a_1} \cdots \tau_{a_n} \la_k \raFab
&=& \sum_{i=1}^n \laFab \tau_{a_1} \cdots \tau_{a_i-1} \cdots \tau_{a_n} \la_k \raFab\quad \quad \text{(string equation),}
\label{string}
\\
\laFab \tau_1 \tau_{a_1} \cdots \tau_{a_n} \la_k \raFab
&=&  (2g-2+n) \laFab \tau_{a_1} \cdots \tau_{a_n} \la_k \raFab \label{dilaton}
\quad \quad \text{(dilaton equation).}
\end{eqnarray*}
\end{pr}

\noindent The proof is a variation of the proof of the usual string and dilaton
equations, and is left as an exercise to the reader (see for example
\cite[Sec.~3.13]{cime}).



\part{GEOMETRY}\label{partB}
In Part~\ref{partB}, we use degeneration and localization to obtain a topological recursion for the Faber-Hurwitz classes, and an expression for the Faber-Hurwitz classes as a sum over a class of weighted trees that are to be defined. The topological recursion is then transformed into a partial differential equation for the Faber-Hurwitz series. 

\section{Degeneration and localization}\label{DandL}

For a
partition $\al$, we use $|\al |$ and $l(\al)$ for the sum of parts
and number of parts, respectively, and write $\al\vdash |\al |$.
If $\al$ has $i_j$ parts equal to $j$, $j\geq 1$, then we also
write $\al =(1^{i_1}2^{i_2}\cd )$, where convenient.
The set of all non-empty partitions is denoted by $\sP$.

\subsection{Relative stable maps}\label{secRsm} We shall use the theory of stable
relative maps to $\proj^1$, following J.~Li's algebro-geometric
description in \cite{li1}, and his description of their
deformation-obstruction theory in \cite{li2}.  (We point out earlier
definitions of relative stable maps in the differentiable category due
to A.-M. Li and Y.~Ruan \cite{lr}, and Ionel and Parker \cite{ip1,
  ip2}, and Gathmann's work \cite{gathmann} in the algebraic category in genus $0$.)
We need the algebraic category for several reasons, most importantly
because we shall use virtual localization, and an explicit description
of the moduli space's deformation-obstruction theory.

\subsubsection{Relative to one point}\label{sssROPT}
The moduli space of genus $g$ relative stable maps to
$\proj^1$ relative to one point $\infty$, with branching above
$\infty$ given by the partition $\al \vdash d$, is denoted by
$\cmbar_{g,\al}(\proj^1)$ where $d$ is the degree of the cover.

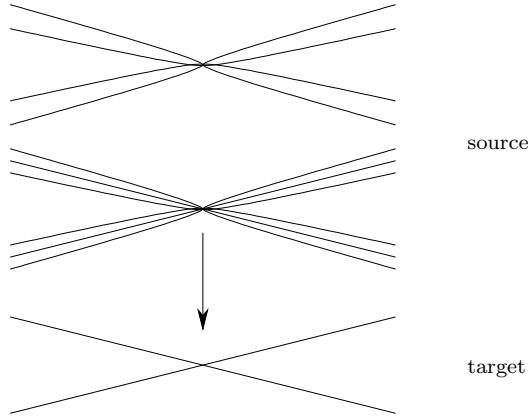
\begin{figure}[ht]
\begin{center}
\setlength{\unitlength}{0.00083333in}
\begingroup\makeatletter\ifx\SetFigFont\undefined%
\gdef\SetFigFont#1#2#3#4#5{%
  \reset@font\fontsize{#1}{#2pt}%
  \fontfamily{#3}\fontseries{#4}\fontshape{#5}%
  \selectfont}%
\fi\endgroup%
{\renewcommand{\dashlinestretch}{30}
\begin{picture}(3200,2589)(0,-10)
\put(2862,237){\makebox(0,0)[lb]{{\SetFigFont{8}{9.6}{\rmdefault}{\mddefault}{\updefault}target}}}
\path(2412,612)(12,12)
\path(12,1587)(2412,987)
\path(12,987)(2412,1587)
\path(1212,1137)(1212,537)
\blacken\path(1182.000,657.000)(1212.000,537.000)(1242.000,657.000)(1212.000,621.000)(1182.000,657.000)
\path(12,1662)(14,1661)(18,1660)
	(25,1658)(37,1655)(54,1650)
	(76,1644)(103,1636)(136,1627)
	(174,1616)(216,1604)(263,1591)
	(314,1576)(367,1561)(422,1545)
	(479,1529)(537,1512)(595,1495)
	(652,1478)(709,1462)(765,1445)
	(819,1429)(871,1413)(920,1398)
	(968,1383)(1012,1368)(1053,1355)
	(1091,1342)(1125,1330)(1154,1318)
	(1178,1308)(1196,1300)(1208,1292)
	(1212,1287)(1208,1284)(1196,1282)
	(1178,1283)(1154,1285)(1125,1288)
	(1091,1293)(1053,1299)(1012,1306)
	(968,1314)(920,1322)(871,1332)
	(819,1342)(765,1353)(709,1364)
	(652,1376)(595,1387)(537,1400)
	(479,1412)(422,1424)(367,1435)
	(314,1447)(263,1458)(216,1468)
	(174,1477)(136,1485)(103,1492)
	(76,1498)(54,1503)(37,1506)
	(25,1509)(18,1511)(14,1512)(12,1512)
\path(2412,1662)(2410,1661)(2406,1660)
	(2399,1658)(2387,1655)(2370,1650)
	(2348,1644)(2321,1636)(2288,1627)
	(2250,1616)(2208,1604)(2161,1591)
	(2110,1576)(2057,1561)(2002,1545)
	(1945,1529)(1887,1512)(1829,1495)
	(1772,1478)(1715,1462)(1659,1445)
	(1605,1429)(1553,1413)(1504,1398)
	(1456,1383)(1412,1368)(1371,1355)
	(1333,1342)(1299,1330)(1270,1318)
	(1246,1308)(1228,1300)(1216,1292)
	(1212,1287)(1216,1284)(1228,1282)
	(1246,1283)(1270,1285)(1299,1288)
	(1333,1293)(1371,1299)(1412,1306)
	(1456,1314)(1504,1322)(1553,1332)
	(1605,1342)(1659,1353)(1715,1364)
	(1772,1376)(1829,1387)(1887,1400)
	(1945,1412)(2002,1424)(2057,1435)
	(2110,1447)(2161,1458)(2208,1468)
	(2250,1477)(2288,1485)(2321,1492)
	(2348,1498)(2370,1503)(2387,1506)
	(2399,1509)(2406,1511)(2410,1512)(2412,1512)
\path(2412,1062)(2410,1062)(2406,1063)
	(2399,1065)(2387,1068)(2370,1071)
	(2348,1076)(2321,1082)(2288,1089)
	(2250,1097)(2208,1106)(2161,1116)
	(2110,1127)(2057,1139)(2002,1150)
	(1945,1162)(1887,1175)(1829,1187)
	(1772,1198)(1715,1210)(1659,1221)
	(1605,1232)(1553,1242)(1504,1252)
	(1456,1260)(1412,1268)(1371,1275)
	(1333,1281)(1299,1286)(1270,1289)
	(1246,1291)(1228,1292)(1216,1290)
	(1212,1287)(1216,1282)(1228,1274)
	(1246,1266)(1270,1256)(1299,1244)
	(1333,1232)(1371,1219)(1412,1206)
	(1456,1191)(1504,1176)(1553,1161)
	(1605,1145)(1659,1129)(1715,1112)
	(1772,1096)(1829,1079)(1887,1062)
	(1945,1045)(2002,1029)(2057,1013)
	(2110,998)(2161,983)(2208,970)
	(2250,958)(2288,947)(2321,938)
	(2348,930)(2370,924)(2387,919)
	(2399,916)(2406,914)(2410,913)(2412,912)
\path(12,1062)(14,1062)(18,1063)
	(25,1065)(37,1068)(54,1071)
	(76,1076)(103,1082)(136,1089)
	(174,1097)(216,1106)(263,1116)
	(314,1127)(367,1139)(422,1150)
	(479,1162)(537,1175)(595,1187)
	(652,1198)(709,1210)(765,1221)
	(819,1232)(871,1242)(920,1252)
	(968,1260)(1012,1268)(1053,1275)
	(1091,1281)(1125,1286)(1154,1289)
	(1178,1291)(1196,1292)(1208,1290)
	(1212,1287)(1208,1282)(1196,1274)
	(1178,1266)(1154,1256)(1125,1244)
	(1091,1232)(1053,1219)(1012,1206)
	(968,1191)(920,1176)(871,1161)
	(819,1145)(765,1129)(709,1112)
	(652,1096)(595,1079)(537,1062)
	(479,1045)(422,1029)(367,1013)
	(314,998)(263,983)(216,970)
	(174,958)(136,947)(103,938)
	(76,930)(54,924)(37,919)
	(25,916)(18,914)(14,913)(12,912)
\path(12,2562)(14,2561)(18,2560)
	(25,2558)(37,2555)(54,2550)
	(76,2544)(103,2536)(136,2527)
	(174,2516)(216,2504)(263,2491)
	(314,2476)(367,2461)(422,2445)
	(479,2429)(537,2412)(595,2395)
	(652,2378)(709,2362)(765,2345)
	(819,2329)(871,2313)(920,2298)
	(968,2283)(1012,2268)(1053,2255)
	(1091,2242)(1125,2230)(1154,2218)
	(1178,2208)(1196,2200)(1208,2192)
	(1212,2187)(1208,2184)(1196,2182)
	(1178,2183)(1154,2185)(1125,2188)
	(1091,2193)(1053,2199)(1012,2206)
	(968,2214)(920,2222)(871,2232)
	(819,2242)(765,2253)(709,2264)
	(652,2276)(595,2287)(537,2300)
	(479,2312)(422,2324)(367,2335)
	(314,2347)(263,2358)(216,2368)
	(174,2377)(136,2385)(103,2392)
	(76,2398)(54,2403)(37,2406)
	(25,2409)(18,2411)(14,2412)(12,2412)
\path(2412,1962)(2410,1962)(2406,1963)
	(2399,1965)(2387,1968)(2370,1971)
	(2348,1976)(2321,1982)(2288,1989)
	(2250,1997)(2208,2006)(2161,2016)
	(2110,2027)(2057,2039)(2002,2050)
	(1945,2062)(1887,2075)(1829,2087)
	(1772,2098)(1715,2110)(1659,2121)
	(1605,2132)(1553,2142)(1504,2152)
	(1456,2160)(1412,2168)(1371,2175)
	(1333,2181)(1299,2186)(1270,2189)
	(1246,2191)(1228,2192)(1216,2190)
	(1212,2187)(1216,2182)(1228,2174)
	(1246,2166)(1270,2156)(1299,2144)
	(1333,2132)(1371,2119)(1412,2106)
	(1456,2091)(1504,2076)(1553,2061)
	(1605,2045)(1659,2029)(1715,2012)
	(1772,1996)(1829,1979)(1887,1962)
	(1945,1945)(2002,1929)(2057,1913)
	(2110,1898)(2161,1883)(2208,1870)
	(2250,1858)(2288,1847)(2321,1838)
	(2348,1830)(2370,1824)(2387,1819)
	(2399,1816)(2406,1814)(2410,1813)(2412,1812)
\path(2412,2562)(2410,2561)(2406,2560)
	(2399,2558)(2387,2555)(2370,2550)
	(2348,2544)(2321,2536)(2288,2527)
	(2250,2516)(2208,2504)(2161,2491)
	(2110,2476)(2057,2461)(2002,2445)
	(1945,2429)(1887,2412)(1829,2395)
	(1772,2378)(1715,2362)(1659,2345)
	(1605,2329)(1553,2313)(1504,2298)
	(1456,2283)(1412,2268)(1371,2255)
	(1333,2242)(1299,2230)(1270,2218)
	(1246,2208)(1228,2200)(1216,2192)
	(1212,2187)(1216,2184)(1228,2182)
	(1246,2183)(1270,2185)(1299,2188)
	(1333,2193)(1371,2199)(1412,2206)
	(1456,2214)(1504,2222)(1553,2232)
	(1605,2242)(1659,2253)(1715,2264)
	(1772,2276)(1829,2287)(1887,2300)
	(1945,2312)(2002,2324)(2057,2335)
	(2110,2347)(2161,2358)(2208,2368)
	(2250,2377)(2288,2385)(2321,2392)
	(2348,2398)(2370,2403)(2387,2406)
	(2399,2409)(2406,2411)(2410,2412)(2412,2412)
\path(12,1962)(14,1962)(18,1963)
	(25,1965)(37,1968)(54,1971)
	(76,1976)(103,1982)(136,1989)
	(174,1997)(216,2006)(263,2016)
	(314,2027)(367,2039)(422,2050)
	(479,2062)(537,2075)(595,2087)
	(652,2098)(709,2110)(765,2121)
	(819,2132)(871,2142)(920,2152)
	(968,2160)(1012,2168)(1053,2175)
	(1091,2181)(1125,2186)(1154,2189)
	(1178,2191)(1196,2192)(1208,2190)
	(1212,2187)(1208,2182)(1196,2174)
	(1178,2166)(1154,2156)(1125,2144)
	(1091,2132)(1053,2119)(1012,2106)
	(968,2091)(920,2076)(871,2061)
	(819,2045)(765,2029)(709,2012)
	(652,1996)(595,1979)(537,1962)
	(479,1945)(422,1929)(367,1913)
	(314,1898)(263,1883)(216,1870)
	(174,1858)(136,1847)(103,1838)
	(76,1830)(54,1824)(37,1819)
	(25,1816)(18,1814)(14,1813)(12,1812)
\put(2862,1662){\makebox(0,0)[lb]{{\SetFigFont{8}{9.6}{\rmdefault}{\mddefault}{\updefault}source}}}
\path(12,612)(2412,12)
\end{picture}
}

\end{center}
\caption{The ``predeformability'' or ``kissing'' condition
on maps of nodes to nodes.  The singularities of the source and target
are nodes (analytically isomorphic to $xy=0$ in $\C^2$),
although it is impossible to depict them as such on the page.
\lremind{kissingfig}}
\label{kissingfig}
\end{figure}

A relative map to $X = \proj^1$ is the following data:
\begin{itemize}
\item
a morphism $f_1$ from a nodal 
$n$-pointed genus $g$ curve $(C, q_1, \dots, q_n)$ (where
as usual the  $q_j$ are distinct non-singular points) to a
chain of $\proj^1$'s, $T = T_0 \cup T_1 \cup \cdots \cup T_t$ (where
$T_i$ and $T_{i+1}$ meet), with a point $\infty \in T_t - T_{t-1}$,
so there are two points named $\infty$.  We will call the one on $X$,  $\infty_X$, and the one on $T$, $\infty_T$, whenever there is any ambiguity.
\item  A projection $f_2: T \rightarrow X$ contracting $T_i$ 
to $\infty_X$ (for $i>0$)
and giving an isomorphism from $(T_0, T_0 \cap T_1)$ 
(resp.\ $(T_0, \infty_T)$) 
to $(X,\infty_X)$ if $t>0$  (resp.\ if $t=0$).  Denote $f_2 \circ f_1$ by $f$.
\item We have an equality of divisors on $C$: $f_1^* \infty_T = \sum \al_i q_i$.
In particular,  $f_1^{-1} \infty_T$ consists 
of non-singular marked points of $C$.
\item The preimage of each node $n$ of $T$ is a union of nodes of $C$.
  At any such node $n'$ of $C$, the two branches map to the two
  branches of $n$, and their orders of branching are the same.
This  is called the {\em predeformability} or {\em kissing} condition.
See Figure~\ref{kissingfig} for a pictorial representation.
Analytically, this map is of the following form. The node $uv=0$ in the
$uv$-plane maps to the node $xy=0$ in the $xy$-plane by $(u,v) \mapsto
(u^m, v^m)=(x,y)$.  (The branching of the $u$-axis over the $x$-axis
is the same as  the branching of the $v$-axis over the $y$-axis.)
\end{itemize}

An {\em isomorphism} of two such maps is a commuting diagram
$$
\xymatrix{
(C, p_1, \dots, p_m, q_1, \dots, q_n) \ar[d]_{f_1} \ar[r]^{\sim} & 
(C', p'_1, \dots, p'_m, q'_1, \dots, q'_n) \ar[d]^{f_1} \\
(T, \infty_T) \ar[r]^{\sim} \ar[d]_{f_2} &  (T, \infty_T) \ar[d]^{f_2}\\
(X,\infty_X) \ar[r]^= & (X,\infty_X)}
$$
where all horizontal morphisms are isomorphisms, the bottom (although
not necessarily the middle!) is an equality, the top horizontal
isomorphism sends $p_i$ to $p'_i$ and $q_j$ to $q'_j$.  Note that the
middle isomorphism must  preserve the isomorphism of $T_0$ with $X$, and is
hence the identity on $T_0$, but for $i>0$, the isomorphism may not be
the identity on $T_i$.

We say that $f$ is {\em stable} if it has finite automorphism group.

Let $r^g_{\al}$ be the ``expected''
number of branch points away from $\infty$ (the number if all the
branching over points other than $\infty$ --- ``away from $\infty$'' --
were {\em simple}, \textit{i.e.}\ corresponding to the partition $(1 ^{d-2} 2) \vdash d$).  
Then
\begin{equation}\label{rgadef}
r^g_{\al} := d+l(\al) + 2g-2,
\end{equation}
by the
Riemann-Hurwitz formula.  Then $\cmbar_{g,\al}(\proj^1)$ supports a
virtual fundamental class $$[\cmbar_{g,\al}(\proj^1)]^{\virt} \in
A_{r^g_{\al}}(\cmbar_{g,\al}(\proj^1))$$ of dimension $r^g_{\al}$.

There is a Fantechi-Pandharipande branch morphism \cite{fanp}\lremind{branchmorphism}
\begin{equation}
\label{branchmorphism}
br:  \cmbar_{g, \al}(\proj^1) \rightarrow \Sym^{r^g_{\al}}(\proj^1)
\end{equation}
sending each map to the set of its branch points in $\proj^1$, excluding
the ones that are ``automatic'' because of the branching over $\infty$.
Note that the total branching $C \rightarrow \proj^1$ above a point in
the target $p \neq \infty \in \proj^1$ is $1$ only if the map is simply branched
above $p$. 

\subsubsection{Relative to two points}\label{sssRTPT}
Similarly, let $\cmbar_{g, \al, \be}(\proj^1)$ be the moduli space of stable
relative maps to $\proj^1$ relative to two points $0$ and $\infty$,
where branching above $0$ is given by a partition $\al \vdash d$ and the
branching above $\infty$ is given by a partition $\be \vdash d$.  By the
Riemann-Hurwitz formula, the 
expected number of branch points away
from $0$ and $\infty$ is 
\begin{equation}\label{rgabdef}
r^g_{\al,\be} := l(\al)+l(\be)+2g-2,
\end{equation}
and the moduli space supports a virtual fundamental class
$$[\cmbar_{g, \al, \be}(\proj^1)]^{\virt} \in A_{r^g_{\al,\be}}
(\cmbar_{g, \al, \be}(\proj^1))$$ of this dimension.
There is also a Fantechi-Pandharipande
branch morphism
$$br:  \cmbar_{g, \al,\be}(\proj^1) \rightarrow \Sym^{r^g_{\al,\be}}(\proj^1).$$

A mild variation of this is the case of relative stable maps to a non-rigidified
target (see \cite[Sec.~2.4]{thmstar}; these are sometimes called ``rubber
maps''), where two relative maps to $(\proj^1, 0, \infty)$ are
considered the same if they differ by an automorphism of the target $\proj^1$
preserving $0$ and $\infty$ (an element of the group $\C^*$).
This moduli space $\cmbar_{g, \al, \be}(\proj^1)_{\sim}$ of such objects supports
a virtual fundamental class of dimension one less than that
of the ``unrigidified'' usual space:
$$[\cmbar_{g, \al, \be}(\proj^1)_{\sim}]^{\virt} \in A_{r^g_{\al,\be}-1}
(\cmbar_{g, \al, \be}(\proj^1)_{\sim}).$$
There is an obvious map $\phi: \cmbar_{g,\al,\be}(\proj^1) \rightarrow
\cmbar_{g,\al,\be}(\proj^1)_\sim$ ``forgetting the rigidification of the target''.
The two fundamental classes are related by the following result.
\label{nonrigiddef}\lremind{nonrigiddef}

\begin{pr} \cite[Lem.~4.6]{thmstar}
\label{bigbusiness}\lremind{bigbusiness}
$$\phi_* \left( br^*(L) \cap [\cmbar_{g, \al, \be}(\proj^1)]^{\virt} \right)
= r^g_{\al, \be} [\cmbar_{g, \al, \be}(\proj^1)_{\sim}]^{\virt}$$
where
$L$ is the class in $\Sym^{r^g_{\al,\be}}(\proj^1)$ 
corresponding to those $r^g_{\al,\be}$-tuples of points containing a given fixed point
$p_0$.
\end{pr}

Intuitively, this is because given any unrigidified map in $\cmbar_{g,
  \al, \be}(\proj^1)_{\sim}$, there should be $r^g_{\al,\be}$ ways to
``rigidify'' it so that a branch point maps to $1 \in \proj^1$, as
there ``should be'' $r^g_{\al,\be}$ branch points away from $0$ and $\infty$.

\subsubsection{Relative stable maps with rational tails} \label{sssSRMRT}
We define the moduli
spaces of relative stable  maps with rational tails 
$\cm_{g,\al}(\proj^1)^{rt}$  as the analogous moduli
spaces of maps where the source curve is a nodal curve with rational
tails.
This is an open substack of the space of relative stable
maps.  We define $\cm_{g,\al,\be}(\proj^1)^{rt}$ similarly.

\subsubsection{Relative stable maps with possibly disconnected source}\label{secpossdis}
The above definitions of relative maps would all make sense
without the requirement that the source curve be connected.  The
resulting moduli space of maps from ``possibly-disconnected'' curves
to $\proj^1$ relative to one point is denoted by
$\cmbar_{g,\al}(\proj^1)_{\bullet}$, and similarly for maps relative
to two points.\label{pds}\lremind{pds}

\subsubsection{Observations using the Riemann-Hurwitz formula}
We make some crucial observations, which are straightforward to verify using
the Riemann-Hurwitz formula.  We shall use them repeatedly.
 A degree $d$ cover is said to be \emph{completely branched} over a point
if the branching data is $(d) \vdash d.$
Define a {\em trivial cover} of $\proj^1$ to be a map $\proj^1
\rightarrow \proj^1$ of the form $[x;y] \mapsto [x^u; y^u]$, 
branched only over $0$ and $\infty$.

\begin{re}
\textrm{
Suppose
we have a map $C \rightarrow \proj^1$ from a nodal (possibly
disconnected) curve, unbranched away from $0$ and $\infty$, where $C$
is non-singular over $0$ and $\infty$.  Then it is a disjoint union of trivial
covers.\label{cr0}\lremind{cr0}
}
\end{re}

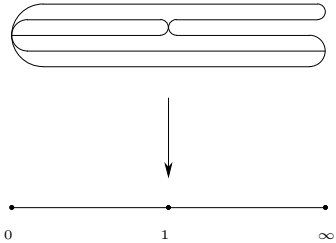
\begin{figure}[ht]
\begin{center}
\setlength{\unitlength}{0.00083333in}
\begingroup\makeatletter\ifx\SetFigFont\undefined%
\gdef\SetFigFont#1#2#3#4#5{%
  \reset@font\fontsize{#1}{#2pt}%
  \fontfamily{#3}\fontseries{#4}\fontshape{#5}%
  \selectfont}%
\fi\endgroup%
{\renewcommand{\dashlinestretch}{30}
\begin{picture}(2025,1517)(0,-10)
\put(0,23){\makebox(0,0)[lb]{{\SetFigFont{5}{6.0}{\rmdefault}{\mddefault}{\updefault}$0$}}}
\put(1076.000,1343.000){\arc{98.000}{1.5708}{4.7124}}
\put(978.000,1343.000){\arc{98.000}{4.7124}{7.8540}}
\put(1955.990,1441.000){\arc{98.021}{4.7330}{7.8334}}
\put(244.500,1294.500){\arc{391.001}{1.5682}{4.7149}}
\put(146.500,1294.500){\arc{195.003}{1.5657}{4.7175}}
\put(49,219){\blacken\ellipse{24}{24}}
\put(49,219){\ellipse{24}{24}}
\put(2005,219){\blacken\ellipse{24}{24}}
\put(2005,219){\ellipse{24}{24}}
\put(1027,219){\blacken\ellipse{24}{24}}
\put(1027,219){\ellipse{24}{24}}
\path(49,219)(2005,219)
\path(1027,903)(1027,414)
\blacken\path(1007.435,492.260)(1027.000,414.000)(1046.565,492.260)(1027.000,468.782)(1007.435,492.260)
\path(245,1490)(1957,1490)
\path(1957,1392)(1076,1392)
\path(1076,1294)(1908,1294)
\path(1908,1099)(245,1099)
\path(49,1294)(978,1294)
\path(978,1392)(147,1392)
\path(147,1197)(2005,1197)
\put(978,23){\makebox(0,0)[lb]{{\SetFigFont{5}{6.0}{\rmdefault}{\mddefault}{\updefault}$1$}}}
\put(1957,23){\makebox(0,0)[lb]{{\SetFigFont{5}{6.0}{\rmdefault}{\mddefault}{\updefault}$\infty$}}}
\put(1907.500,1196.500){\arc{195.003}{4.7175}{7.8489}}
\end{picture}
}

\end{center}
\caption{Almost-trivial covers:  All branched covers from non-singular curves 
with simple branching
above $1$ and no other branching away from $0$ and $\infty$
must look like this. \lremind{almosttrivial}}
\label{almosttrivial}
\end{figure}

\begin{re}
a) Suppose we have a map from a nodal curve $C$ to
$\proj^1$, with no branching away from $0$, $1$, and $\infty$, simple
branching at $1$, and non-singular over $0$ and $\infty$.  Then it is a
union of trivial covers, 
with one further component, that is
non-singular of genus $0$, completely branched over one of $\{ 0,
\infty \}$, and with two preimages over the other (see
Figure~\ref{almosttrivial}).  We call this an {\em almost-trivial cover}.
\label{cr1}\lremind{cr1}

b) More generally, suppose
we have a map from a curve $C$ to a chain of $\proj^1$'s, satisfying
the kissing or predeformability condition, unbranched except
for two non-singular points $0$ and $\infty$ on the ends of the chain, and
simple branching at one more point.  Then the map looks like a number
of trivial covers glued together, 
with one further almost-trivial cover $\proj^1
\rightarrow \proj^1$ of the sort described in Remark~\ref{cr1}(a).
\end{re}

\begin{re}
a) If we have a map from a nodal curve $C$ to
$\proj^1$, with total branching away from $0$ and $\infty$ of degree
less than $2g$, and non-singular over $0$ and $\infty$, then $C$ has no
components of geometric genus $g$.  If instead the total branching
away from $0$ and $\infty$ is precisely $2g$, and $C$ has a component
of geometric genus $g$, then the cover is a disjoint union of trivial
covers, and one connected curve $C'$ of arithmetic genus $g$, where
the map $C' \rightarrow \proj^1$ is completely branched over $0$ and
$\infty$.  \label{cr2}\lremind{cr2}

b) Suppose  more generally we have a map from a curve $C$ to a 
chain of $\proj^1$'s,
satisfying the kissing condition, with branching
of less than $2g$ away from the nodes and two smooth points $0$ and
$\infty$ on the ends of the chain, then $C$ has no component of geometric
genus $g$.
If the branching is precisely $2g$ away from the nodes and $0$ and
$\infty$, then the map is a number of trivial covers glued
together, plus one other cover of the sort described in Remark~\ref{cr2}(a).
\end{re}

The following fact will prove essential.
\begin{tm}\label{tGgd}
 Let $\G_{g,d, \sim} := \mu_*
  \left[ \cm_{g,(d),(d)}(\proj^1,d)_{\sim}^{rt} \right]^{\virt}$
  where $\mu$ is the moduli map to $\cm_{g,2}^{rt}$,
  remembering the curve and the points over  $0$ and $\infty$.  Then $\G_{g,d,
    \sim} = d^{2g} \G_{g,1, \sim}$.  
Similarly, if $\G_{g,d} := \mu_*
  \left[ br^*(L) \cap \cm_{g,(d),(d)}(\proj^1,d)^{rt}
  \right]^{\virt}$, then $\G_{g,d} = d^{2g} \G_{g,1}$.
 \lremind{trick}\label{trick}
\end{tm}

\bpf 
The second statement follows from the first by Proposition~\ref{bigbusiness}.
We now prove the first.

Let $\nu: \operatorname{Jac} \rightarrow \cm_{g,2}^{rt}$ be the universal Jacobian or
Picard stack over $\cm_{g,2}^{rt}$. 
Let $\operatorname{Jac}[d] := \ker( \xymatrix{ \operatorname{Jac} \ar[r]^{\times d} & \operatorname{Jac}})$ be the
$d$-torsion substack.
(Points of $\operatorname{Jac}[d]$ correspond to curves along with a $d$-torsion  point.)

Let $p$ and $q$ be the names of
the points of the genus $g$ curve parametrized by $\cm_{g,2}^{rt}$.
Let $s_{p-q}: \cm_{g,2}^{rt} \rightarrow \operatorname{Jac}$ be
the section corresponding to the line bundle $\oh(p-q)$.  
There is a natural (stack-theoretic) isomorphism \lremind{misosoup}
\begin{equation}\label{misosoup}
\cm_{g,(d),
  (d)}(\proj^1)^{rt}_{\sim} \cong \operatorname{Jac}[d] \cap s_{p-q}(\cm_{g,2}^{rt}),
\end{equation}
 as follows.
Given any family of such  relative   stable  maps (with rational tails), let
$p$ and $q$  be the preimage of $0$ and $\infty$ respectively.
Note that the target $\proj^1$ never degenerates (or sprouts)
in such a family.
Thus $\oh(dp - dq) \cong \oh$, and hence $\oh(p-q)$ is indeed 
$d$-torsion.  Conversely, given any family of curves in $\cm_{g,2}^{rt}$ where
$\oh(p-q)$ is  $d$-torsion,
we obtain a unique family of relative  stable maps with unrigidified target
of the desired sort.  (It is essential to note that this isomorphism holds
over the boundary of $\cm_{g,2}^{rt}$ as well.)
$$
\xymatrix{
\cm_{g,(d),(d)}(\proj^1)^{rt}_{\sim} = \operatorname{Jac}[d] \cap s_{p-q} (\cm_{g,2}^{rt}) \ar[rd]^{\mu} \ar[r] & \operatorname{Jac} \ar[d]_{\nu} \\
& \cm^{rt}_{g,2} \ar@/_/[u]_{s_{p-q}} 
}
$$

Both sides of the isomorphism \eqref{misosoup} have natural
virtual fundamental classes, the latter by intersecting the section
$s_{p-q}(\cm_{g,2}^{rt})$ with the local complete intersection $\operatorname{Jac}[d]
\hookrightarrow \operatorname{Jac}$.  We claim that these virtual fundamental classes
agree.  (It is easy to show that away from the boundary of $\cm_{g,2}^{rt}$, both will
agree with the actual fundamental class, but we shall
not need this fact.)
We show this by showing that they have the same deformation-obstruction
theory over $\cm_{g,2}^{rt}$.

We use the description of the deformation-obstruction theory  
of
$\cmbar_{g,(d),
  (d)}(\proj^1)$ over $\cmbar_{g,2}$ given in \cite[Sec.~2.8]{thmstar}.
On the locus, at a relative stable map $[f: C \rightarrow \proj^1]$, 
the relative deformation space
is identified with  $H^0(C, f^* T_{\proj^1}(-[0] - [\infty])) = H^0(C, f^* \oh_{\proj^1})$ in \cite[equ.\ (1)]{thmstar}.
(The formula given in \cite{thmstar} is more involved, as it applies in  more general
situations.  In particular, in the locus, the target never sprouts, so
the notation $f^{\dagger}$ in \cite{thmstar}
agrees with $f^*$.
Also $T_{\proj^1}(-\log [0] - \log [\infty]) = T_{\proj^1}(-[0] - [\infty])$, as
$\proj^1$ has dimension $1$.)  Thus the relative deformation space is 
canonically identified with
$H^0(C, \oh_C)$, which is one-dimensional.  But such deformations correspond
precisely to the $\C^*$-action induced by its action on the target (keeping the
map otherwise fixed), so the relative deformation space of maps to $\proj^1$ with
{\em unrigidified target} is $0$.

The relative obstruction space $\operatorname{RelOb}(f)$ at a 
relative stable map $[f: C \rightarrow \proj^1]$ has a natural filtration
\cite[equ.\ (2)]{thmstar}
$$
0 \rightarrow H^1(C, f^* T_{\proj^1}(- [0] - [\infty])) \rightarrow
\operatorname{RelOb}(f) \rightarrow H^0(C, f^{-1} \mathcal{E}xt^1(
\Omega_{\proj^1}( -[0] - [\infty]), \oh_{\proj^1})) \rightarrow 0.$$
Now $\mathcal{E}xt^1( \oh_{\proj^1}, \oh_{\proj^1}) = 0$, so this
filtration gives a natural isomorphism $$\xymatrix{ H^1(C, \oh_C)
  \ar[r]^{\sim} & \operatorname{RelOb}(f)}.$$  But this is the relative
obstruction space of the intersection: the normal bundle to the
$d$-torsion locus $\operatorname{Jac}[d]$ is canonically $H^1(C, \oh_C)$.
Thus we have shown that the isomorphism \eqref{misosoup}
indeed extends to an isomorphism of virtual fundamental classes.

We conclude by noting that $[\operatorname{Jac}[d]] = d^{2g} [\operatorname{Jac}[1]]$
in $A_*(\operatorname{Jac})$.  This follows from \cite[Sec.~2.15]{dm}, as observed
in the proof of \cite[Lem.~2.10]{looijenga}.
\epf

\noindent {\em Remark.} 
$\G_{g,1}$ will be our ``natural generator'' of $R_{2g-1}(\cm_{g,1})$.
The generator used by Looijenga in \cite{looijenga} is the
hyperelliptic locus $\Hyp$.  They are related
as follows:  $(2g+2)(2g+1) \Hyp = (2^{2g}-1) \G_{g,1,\sim}$
(both count hyperelliptic covers with two fixed branch points),
and $\G_{g,1} = 2g \G_{g,1,\sim}$ from Proposition~\ref{bigbusiness};
see Corollary~\ref{hyp} and its proof.

\subsection{Faber-Hurwitz classes}
We define Hurwitz classes,
following \cite{thmstar}.
Our motivation is as follows. Suppose we are interested
in dimension $j$ classes on $\cmbar_{g,n}$.  One way
of producing a family of genus $g$ curves
with $n$-marked points is by considering branched covers
of $\proj^1$ with fixed branching over $\infty$
corresponding to a partition $(\al_1, \ldots, \al_n)$ of $d$
with $n$ parts. 
The curve in question will be the source of the map to $\proj^1$, and the $n$ points will
be the points over $\infty$. 
This Hurwitz scheme (the moduli space of such maps) 
will have dimension $r^g_{ \al}$.
In order to get a class of dimension $j$, we fix all but $j$
branch points.  
We can then push forward this class to $\cmbar_{g,n}$.
The definition of ``Hurwitz class'' involves doing this 
``virtually'':\label{hcdef}\lremind{hcdef}
$$
\H^{g,\al}_j := \mu_* \left( br^*(L)^{r^g_{\al}-j} \cap
[ \cmbar_{g, \al}(\proj^1)]^{\virt} \right)$$
where $\mu$ is the forgetful map $\cmbar_{g,\al}(\proj^1) \rightarrow
\cmbar_{g,n}$ and $br$ is the branch morphism \eqref{branchmorphism}.
Here, as in Proposition~\ref{bigbusiness}, $L$ is the class
in $\Sym^{r^g_{\al}}(\proj^1)$ corresponding to unordered tuples of points
containing a given fixed point $p_0$.

Define the \emph{Faber-Hurwitz} class \lremind{fhclass}\label{FHdef}\lremind{FHdef}
$\F^{g,\al}$ by
\begin{equation}\F^{g,\al} := \H^{g,\al}|_{\cm^{rt}_{g,n}} \in
A_{2g-1}(\cm_{g,n}^{rt}).\label{fhclass}
\end{equation}
Equivalently, we can consider the space of relative stable  maps with
rational tails $\cm_{g,\al}(\proj^1)^{rt}$, along with its virtual
fundamental class; fix all but $2g-1$ branch points; and push forward
to the moduli space $\cm_{g,n}^{rt}$.
Let 
\begin{equation}\label{rFgadef}\lremind{rFgadef}
r^{\Faber}_{g,\al} := r^g_{\al} - (2g-1) = d+l(\al)-1
\end{equation}
be the number of {\em fixed} branch points on $\proj^1$ in this construction.

\begin{boxitpara}{box 1.00 setgray fill }
We shall now understand this class in two ways, by degeneration and
localization.  The first will connect us to a Hurwitz-type problem
(and join-cut type recursion) that we can solve.  The second will
connect us to the tautological ring.\end{boxitpara}

Readers less comfortable with these ideas from Gromov-Witten
theory may go  directly to the recursions that are obtained by
these means without losing the thread of the paper. These
recursions are treated quite cleanly through transforms of the corresponding
partial differential equations (see Section~\ref{s:genseries}).

\subsection{Degeneration of Faber-Hurwitz classes}
We now describe a join-cut type recursive formula for the Faber-Hurwitz
classes $\F^{g,\al}$.  This will allow us to compute $\F^{g,\al}$ in terms of
the putative
generator $\G_{g,1}$, first recursively, and later, in closed form.
We use Jun Li's degeneration formula \cite{li1, li2}, which in our case states
that, if we degenerate
$[ \cmbar_{g, \al}(\proj^1)_\bullet]^{\virt}$ by degenerating the target
into two components $\proj^1_L \cup \proj^1_R$, where the
point corresponding to $\al$ is on the right component $\proj^1_R$,
then
\begin{equation}\label{cmvirtsum}
[ \cmbar_{g, \al}(\proj^1)_\bullet]^{\virt} 
= \sum_{g_1,g_2,\ga} \left( \prod_i \ga_i \right) 
[ \cmbar_{g_1, \ga}(\proj^1)_\bullet]^{\virt}  \boxtimes
[ \cmbar_{g_2, \ga, \al}(\proj^1)_\bullet]^{\virt} 
\end{equation}
where the sum is over all ``splitting of the data'':
$|\ga|=|\al|$,
$g_1 + g_2 + l(\ga) -1 = g$, and the genus $g_1$ cover $C_1$
is glued to the genus $g_2$ cover $C_2$ over the node $\proj^1_L \cap \proj^1_R$
by gluing the point corresponding to $\ga_i$ on $C_1$ to the 
point corresponding to $\ga_i$ on $C_2$.
There is the obvious variation for
$[ \cmbar_{g, \al, \be}(\proj^1)_\bullet]^{\virt}$ for
stable maps relative to two points.  If we are interested
in spaces with connected source, \textit{i.e.}\ 
$\cmbar_{g, \al}(\proj^1)$ rather than
$\cmbar_{g, \al}(\proj^1)_\bullet$, then we include only the summands
where the union $C_1 \cup C_2$ is connected.
We shall ``cap'' these equalities with pullbacks under
the branch morphism, and interpret them as considering
maps with given branch points, where we shall specify how
many branch points degenerate to $\proj^1_L$ and $\proj^1_R,$ respectively.

\begin{lm}\label{beginob}
If $j< 2g-1$, then the restriction of $\mu_* ( \H^{g,\al}_j)$ to the
rational-tails locus $\cm^{rt}_{g,n}$ is $0$.  
\end{lm}
\bpf  
Degenerate the target $\proj^1$ into a chain of $(r^g_{\al}-j)$
$\proj^1$'s, where each of the $r^g_{\al}-j$ fixed branch points
lies in a different component of the target.  Then by the degeneration
formula, this class will be obtained as a sum (over all choices of
splittings of the data) of classes glued together from virtual
fundamental classes of various spaces of relative stable maps to each
component.  But any such relative stable map to one component, with
points $0$ and $\infty$ connecting it to adjacent elements of the
chain (with the obvious variation if it is the end of the chain) has at
most $j+1<2g$ branch points away from $0$ and $\infty$.  Hence by
Remark~\ref{cr2}, there is no curve of geometric genus $g$ mapping to this
component.  Thus, after degeneration, the source curve can have no
irreducible component of genus $g$, and hence the class is $0$ in
$\cm_{g,n}^{rt}$.
\epf

\begin{tm}[Degeneration Theorem --- Join-cut recursion] \label{Fcutjoin}
If $r^{\Faber}_{g,\al} > 0$, then\lremind{FcutjoinE}
\begin{equation*}
\F^{g,\al}
= \sum_{i+j=\al_k}  ij H^0_{\al'} \F^{g, \al''}
 \binom {r^{\Faber}_{g,\al}-1} {r^{\Faber}_{g, \al''}}
+
\sum_{\al_i, \al_j} (\al_i + \al_j) \F^{g,\al'}
+
\sum_{i=1}^{l(\al)} \al_i^{2g+1} H^0_{\al} \G_{g,1}.
\end{equation*}
where $\al'$, $\al''$, etc., are as defined in the proof below.
 \lremind{Fcutjoin}
\end{tm}

This recursion inductively determines $\F^{g,\al}$ in terms of $\G_{g,1}$
(Theorem~\ref{connec}).  There is no ``base case'' necessary.
This result is best stated in terms  of generating
series (see Corollary~\ref{JCEforFHS}).

\bpf
We obtain a recursion for Faber-Hurwitz classes by a similar argument
to that used in Lemma~\ref{beginob} above.
We degenerate the target into two pieces $\proj^1_L \cup \proj^1_R$,
where $\infty$ and one of the $r^{\Faber}_{g,\al} = r^g_{\al} - (2g-1)
= d+l(\al)-1$ fixed branch points are on the right component
$\proj^1_R$, and the remaining $r^{\Faber}_{g,\al}-1$ fixed branch
points are on the left component $\proj^1_L$.  By the degeneration
formula~(\ref{cmvirtsum}),
this class will be obtained as a sum (over all choices of
splittings of the data) of classes glued together from virtual
fundamental classes of various spaces of relative stable maps to each
component.  We discard every term that does not have a smooth component of
genus $g$.

There are two cases: the genus $g$ curve maps to the left component
$\proj^1_L$ or the right component $\proj^1_R$.  

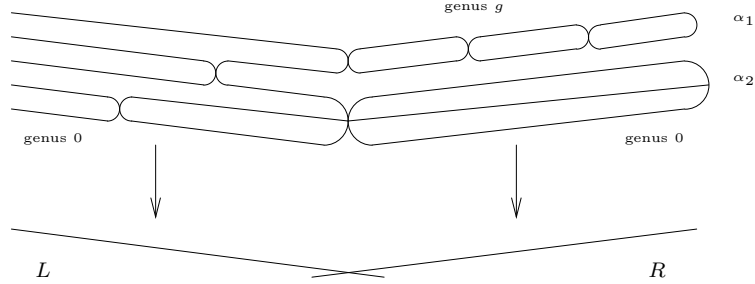
\begin{figure}[ht]
\begin{center}
\setlength{\unitlength}{0.00083333in}
\begingroup\makeatletter\ifx\SetFigFont\undefined%
\gdef\SetFigFont#1#2#3#4#5{%
  \reset@font\fontsize{#1}{#2pt}%
  \fontfamily{#3}\fontseries{#4}\fontshape{#5}%
  \selectfont}%
\fi\endgroup%
{\renewcommand{\dashlinestretch}{30}
\begin{picture}(4512,1756)(0,-10)
\put(2712,1681){\makebox(0,0)[lb]{{\SetFigFont{5}{6.0}{\rmdefault}{\mddefault}{\updefault}genus $g$}}}
\put(2262.000,1006.000){\arc{300.000}{1.5708}{4.7124}}
\put(2187.000,1381.000){\arc{150.000}{1.5708}{4.7124}}
\put(2037.000,1381.000){\arc{150.000}{4.7124}{7.8540}}
\put(762.000,1081.000){\arc{150.000}{1.5708}{4.7124}}
\put(612.000,1081.000){\arc{150.000}{4.7124}{7.8540}}
\put(1362.000,1306.000){\arc{150.000}{1.5708}{4.7124}}
\put(1212.000,1306.000){\arc{150.000}{4.7124}{7.8540}}
\put(2787.000,1456.000){\arc{150.000}{4.7124}{7.8540}}
\put(2937.000,1456.000){\arc{150.000}{1.5708}{4.7124}}
\put(3537.000,1531.000){\arc{150.000}{4.7124}{7.8540}}
\put(3687.000,1531.000){\arc{150.000}{1.5708}{4.7124}}
\put(4212.000,1606.000){\arc{150.000}{4.7124}{7.8540}}
\put(4212.000,1231.000){\arc{300.000}{4.7124}{7.8540}}
\path(2187,1306)(2787,1381)
\path(2187,1456)(2787,1531)
\path(2937,1531)(3537,1606)
\path(2937,1381)(3537,1456)
\path(3687,1456)(4212,1531)
\path(3687,1606)(4212,1681)
\path(1962,856)(762,1006)
\path(762,1156)(2112,1006)
\path(612,1006)(12,1081)
\path(612,1156)(12,1231)
\path(1962,1156)(1362,1231)
\path(2037,1306)(1362,1381)
\path(1212,1231)(12,1381)
\path(1212,1381)(12,1531)
\path(2037,1456)(12,1681)
\path(2262,1156)(4212,1381)
\path(2112,1006)(4362,1231)
\path(2262,856)(4212,1081)
\path(912,856)(912,406)
\path(882.000,526.000)(912.000,406.000)(942.000,526.000)
\path(3162,856)(3162,406)
\path(3132.000,526.000)(3162.000,406.000)(3192.000,526.000)
\path(2337,31)(12,331)
\path(1887,31)(4287,331)
\put(3987,31){\makebox(0,0)[lb]{{\SetFigFont{8}{9.6}{\rmdefault}{\mddefault}{\updefault}$R$}}}
\put(162,31){\makebox(0,0)[lb]{{\SetFigFont{8}{9.6}{\rmdefault}{\mddefault}{\updefault}$L$}}}
\put(4512,1606){\makebox(0,0)[lb]{{\SetFigFont{5}{6.0}{\rmdefault}{\mddefault}{\updefault}$\alpha_1$}}}
\put(4512,1231){\makebox(0,0)[lb]{{\SetFigFont{5}{6.0}{\rmdefault}{\mddefault}{\updefault}$\alpha_2$}}}
\put(3837,856){\makebox(0,0)[lb]{{\SetFigFont{5}{6.0}{\rmdefault}{\mddefault}{\updefault}genus $0$}}}
\put(87,856){\makebox(0,0)[lb]{{\SetFigFont{5}{6.0}{\rmdefault}{\mddefault}{\updefault}genus $0$}}}
\put(1962.000,1006.000){\arc{300.000}{4.7124}{7.8540}}
\end{picture}
}

\end{center}
\caption{First case in degeneration argument in proof of Theorem~\ref{Fcutjoin}\lremind{firstcase}}
\label{firstcase}
\end{figure}

\textbf{Case~1:}
If the genus $g$ curve maps to the
right component $\proj^1_R$, then by Remark~\ref{cr2}, all the moving branch
points must also map to $\proj^1_R$ (we need $2g$ branching
total away from $\proj^1_L \cap \proj^1_R$ and $\infty$ in order
to have a genus $g$ curve in the cover).  Also by Remark~\ref{cr2}(a), the cover of
$\proj^1_R$ must be a union of trivial covers, together with one cover by a
connected arithmetic genus $g$ curve that is completely branched over $\infty$
and $\proj^1_L \cap \proj^1_R$.  In particular, the partition over the
node $\proj^1_L \cap \proj^1_R$ must be the same as the partition
$\al$ over $\infty$, and the genus $g$ component can be on the
component of the source corresponding to any $\al_j$.  Also, the cover
of $\proj^1_L$ must be genus $0$ and connected, and the number of such covers
is the genus
$0$ Hurwitz number
\begin{equation}\label{singHur}
H^0_{\al} = (d-2+l(\al))!\; d^{l(\al)-3}
 \prod_{i=1}^{l(\al )}
\frac {{\al_i}^{\al_i}} {\al_i!} .
\end{equation}
(This celebrated formula was first given by Hurwitz, and has now been 
proved in 
many ways, see for example \cite{gj0}.)
This case is depicted pictorially in Figure~\ref{firstcase};
here the genus $g$ curve corresponds to $\al_1$.
We thus get a contribution to $\F^{g,d}$ of \lremind{cont1}
\begin{equation}\label{cont1}
\sum_{i=1}^{l(\al)}\al_i H^0_{\al} \G_{g,\al_i} = 
\sum_{i=1}^{l(\al)} \al_i^{2g+1} H^0_{\al} \G_{g,1}.
\end{equation}
The equality in \eqref{cont1} is through Theorem~\ref{trick}. 
The contributions to the left side of \eqref{cont1}
are as follows:  we have a contribution of $\G_{g,\al_i}$
from the map to $\proj^1_L$.  We have a contribution of $\G_{g,\al_i}
\prod_{j \neq i} (1/\al_j)$ from the map to $\proj^1_R$ 
(since the trivial cover contributes the order~$\al_k$ of the automorphism group 
and we are therefore counting objects weighted by $1/\al_k.$).
We have a multiplicity
of $\prod \al_i$ from the kissing condition.

\textbf{Case 2a:}
If otherwise the genus $g$ curve maps to the left component $\proj^1_L$, 
then by Lemma~\ref{beginob}, all the non-fixed branch points must also map to
$\proj^1_L$
(as we need at least $2g-1$ moving branch points in order to get a non-zero contribution),
so the only branching over $\proj^1_R$ away from its $0$ and $\infty$
is simple branching over one (fixed) point.
By Remark~\ref{cr1}, the cover of $\proj^1_R$ is a union of trivial covers, with one
almost-trivial cover.  If the almost-trivial cover of $\proj^1_R$  is completely branched
over $\proj^1_L \cap \proj^1_R$, and has two preimages over $\infty$
(\textit{e.g.}\ as shown in Figure~\ref{almosttrivial}),
then it can connect any two of the points over $\infty$ (corresponding to two parts
of the partition $\al$, say $\al_i$ and $\al_j$).  The virtual fundamental class of this space
of relative stable maps to $\proj^1_R$ is then $\al_i \al_j / \prod \al_i$:  any
trivial cover corresponding to $\al_k$ ($k \neq i,j$) 
is weighted by $1/\al_k.$
The 
space of relative maps to $\proj^1_L$ corresponds to maps from a connected
curve of genus $g$ curve, with branching over $\proj^1_L \cap \proj^1_R$ corresponding
to the partition $\al'$  obtained by removing $\al_i$ and $\al_j$ from $\al$, and adding $\al_i+\al_j$,
with all but $2g-1$ branch points fixed.  In other words, the contribution is precisely
$\F^{g,\al'}$.  Finally, the multiplicity in the degeneration formula is the product of
the multiplicities of the kissing over the node $\proj^1_L \cap \proj^1_R$, \textit{i.e.}\ 
$\prod \al'_i = (\prod \al_i)(\al_i+\al_j)/(\al_i \al_j)$.  Thus we obtain a contribution 
to $\F^{g,\al}$ of
\lremind{cont2}
\begin{equation*}
\label{cont2}
\sum_{\al_i, \al_j} (\al_i + \al_j) \F^{g,\al'}.
\end{equation*}

\textbf{Case 2b:}
Finally, if the almost-trivial cover over $\proj^1_R$ has one preimage
over $\infty$ (corresponding to $\al_k$, say), and two preimages over
$\proj^1_L \cap \proj^1_R$ ($q_1$ and $q_2$, say, branching with
multiplicity $i$ and $j$ respectively, where $i+j=\al_k$), then we
have contributions of $\al_k/ \prod \al_i$ from the cover of
$\proj^1_R$, and $\prod \al_i \times ij/ \al_k$ from the kissing
multiplicities.  We next examine the contribution from the stable
relative maps to $\proj^1_L$.  The cover must have two connected
components, one containing the point $q_1$ and one containing the
point $q_2$.  As the map has a component of geometric genus $g$, one
of these two components must be arithmetic genus $g$, and the other
must be arithmetic genus $0$.  
Suppose  the genus $0$ curve contains to
$q_1$, and the genus $g$ curve contains $q_2$.  
Suppose the partition
corresponding to the genus $0$ curve is $\al'$ (\textit{i.e.}\ its ramification
above the node $\proj^1_L \cap \proj^1_R$), and the partition
corresponding to the genus $g$ curve is $\al''$, so $\al' +
\al'' = \al - \al_k + \{ i, j \}$, and $i \in \al'$, $j \in \al''$.
In order for the contribution to be non-zero, by Lemma~\ref{beginob},
all of
the $(2g-1)$ ``moving'' branching must belong to the genus $g$
component, so 
all of the genus $0$ curve's branch points must be fixed. 
There are $\binom {r^{\Faber}_{g,\al}-1} {r^{\Faber}_{g,
    \al''}}$ ways of choosing which fixed branch point on $\proj^1_L$ 
belongs to the genus
$0$ curve, and which belongs to the genus $g$ curve.  Hence the
contribution of the cover of $\proj^1_L$ is $H^0_{\al'} \F^{g,\al}
\binom {r^{\Faber}_{g,\al}-1} {r^{\Faber}_{g, \al''}}$, so the combined
contribution to $\F^{g,\al}$ from this case is
\begin{equation*}
\label{cont3}
\sum_{i+j=\al_k}  ij H^0_{\al'} \F^{g, \al''} \binom {r^{\Faber}_{g,\al}-1} {r^{\Faber}_{g, \al''}}.
\end{equation*}
Summing the three types  of contribution, we obtain the result.
\epf

\begin{tm}\label{connec}
 For each fixed $g$, $\F^{g,\al}$ is a
  rational multiple of $\G_{g,1}$, as determined by the recursion of
  Theorem~\ref{Fcutjoin}.
\end{tm}

\bpf
By Remark~\ref{cr2}(a), if $r^{\Faber}_{g,\al}=0$, then $\F^{g,\al} = 0$;
this is the trivial base case of the recursion in Theorem~\ref{Fcutjoin}. The
result follows by induction.
\epf

For convenience, we define the {\em Faber-Hurwitz number} $F^g_{\al} \in \Q$
\label{FHn}\lremind{FHn}
to be this multiple of $\G_{g,1}$, to remind us that these classes are
all commensurate:\lremind{ianeq} 
\begin{equation}\label{ianeq}
F^g_{\al} \G_{g,1}:=\F^{g,\al}.
\end{equation}
(Corollary~\ref{Fnice} is a sign that $F^g_{\al}$ is a well-behaved quantity to consider.)

Theorem~\ref{Fcutjoin} is best formulated in terms of generating
series.  A natural generating
series for genus $0$ Hurwitz numbers is
\begin{equation}\label{sHser}
\widH^0(z;\bfp):=\sum_{\al \in\sP}
z^{|\al |} \f{p_{\al}}{|\ts{Aut}\al |}
\frac{H^0_{\al }}{r^0_{\al}!},
\end{equation}
where $\bfp =(p_1,p_2,\ld )$,
and a natural generating series for
Faber-Hurwitz numbers, which we shall call the \emph{Faber-Hurwitz series},  is\lremind{gsFnum}
\begin{equation}\label{gsFnum}
F^g(z;\bfp):=\sum_{\al\in\sP}z^{|\al |}
\frac{p_{\al}}{|\ts{Aut}\al |}\frac{F^g_{\al}}{r_{g,\al}^{\Faber}!}.
\end{equation}
Then the recursion of Theorem~\ref{Fcutjoin}, for {\em classes}, becomes the following (linear) partial
differential equation for~$F^g$.\lremind{Fcutjoin2}

\begin{co}[Degeneration Theorem --- Join-cut Equation for Faber-Hurwitz Series]\label{JCEforFHS}
For $g\ge 1$, $F^g(z;\bfp)$ is the unique formal power series solution to
\begin{eqnarray*}
\left( z\frac{\pa}{\pa z}-1+\sum_{i\geq 1} p_i\frac{\pa}{\pa p_i}\right) F^g&=
&\sum_{i,j\geq 1}p_{i+j}\left( i\frac {\pa}{\pa p_i}\widH^0\right)
\left( j\frac {\pa}{\pa p_j}F^g\right)\\
&& + \frac{1}{2}\sum_{i,j\geq 1}p_i p_j(i+j)\frac{\pa}{\pa p_{i+j}}F^g
+ \sum_{i\geq 1}i^{2g+1}p_i\frac{\pa}{\pa p_i}\widH^0,\notag
\end{eqnarray*}
\end{co}

In Corollary~\ref{JCEforFHS}, we have changed from classes to numbers
by erasing $\G_{g,1}$ from Theorem~\ref{Fcutjoin} \textit{via}~(\ref{ianeq}).
To obtain a generating series for Faber-Hurwitz {\em classes},
we would simply multiply this series by $\G_{g,1}$.

\subsection{Localization of Faber-Hurwitz classes}\label{ssLFHC}
This section contains the  Localization Tree Theorem for Faber-Hurwitz 
classes. It is to be thought of in conjunction with the Degeneration Theorem.
The theorem gives the fundamental 
relationship between Faber-Hurwitz classes on the one hand, and intersection numbers and (genus 
$0$) single and double  Hurwitz numbers on the other. Implicitly, it allows us to determine the
Faber-Hurwitz numbers.
It involves a sum over a set~$\sT_{g,m}$, 
which is one of three sets of trees defined as follows.

\begin{df}[Localization trees] \label{Dloctr} 
Consider rooted trees with the following properties.
For each  tree $\rt$ there are three classes of vertices: $\sV_0(\rt)$, the $0$-vertices; $\sV_\infty(\rt)$,
the $\infty$-vertices; $\sV_t(\rt)$, the  $t$-vertices.
All $t$-vertices are monovalent (here $t$ stands for ``tail''). The number
 of non-root $0$-vertices in $\rt$ is denoted by $\eta_0(\rt)$.
The non-root $0$-vertices are labelled (each
receives one of $\eta_0(\rt)$ labels in all possible ways).  The $\infty$-vertices are not labelled.
There are two classes of edges: $\sE_{0\infty}(\rt)$, the $0\infty$-edges, each of
which joins a $0$-vertex to an $\infty$-vertex;
 $\sE_{\infty t}(\rt)$, the $\infty t$-edges, each of
which joins a non-root $\infty$-vertex to a $t$-vertex.

There is at least one edge.  Each edge is assigned a positive integer
weight, and the integer weight on an $0\infty$-edge $e$ is denoted by
$\ep(e)$.  The $t$-vertices incident with an edge of weight $k$ are
labelled among themselves, for each $k\ge 1$.  The list of the weights
on all $0\infty$-edges incident with a $0$-vertex $v$ is specified by
the partition $\de ^v(\rt)$, the list of the weights on all
$0\infty$-edges incident with an $\infty$-vertex $v$ is specified by
the partition $\be ^v(\rt)$, and the list of the weights on all
$\infty t$-edges incident with a non-root $\infty$-vertex $v$ is
specified by the partition $\ga ^v(\rt)$.  For each non-root
$\infty$-vertex $v$, we impose the condition that
\begin{equation}\label{balance}
|\be^v(\rt)|=|\ga^v(\rt)|.
\end{equation}

The root-vertex may be either a $0$-vertex or an $\infty$-vertex,
and is denoted by $\bullet$. 
For $m\ge 1$, we define $\sT_{g,m}$ to be the set of rooted trees above
in which the root-vertex is a $0$-vertex of degree $m$.
For $j\ge 1$, we define $\sT_{0,j}$ to be the subset of $\sT_{g,1}$ in
which the edge incident with the root-vertex has weight $j$.
For $j\ge 1$, we define $\sT_{\infty,j}$ to be the set of rooted trees above
in which the root-vertex is a  monovalent $\infty$-vertex, and
the edge incident with the root-vertex has weight $j$.

We refer to the trees in any of the sets $\sT_{g,m}$, $\sT_{0,j}$ and $\sT_{\infty,j}$ as
``localization trees''.
\end{df}

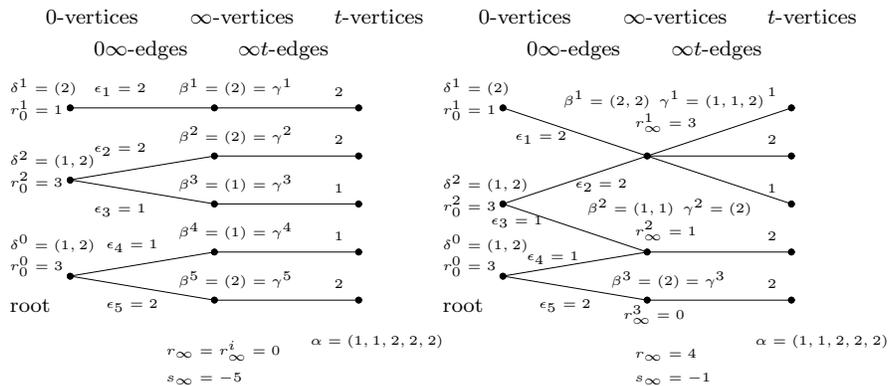
\begin{figure}[h]
\begin{center}
\setlength{\unitlength}{0.00083333in}
\begingroup\makeatletter\ifx\SetFigFont\undefined%
\gdef\SetFigFont#1#2#3#4#5{%
  \reset@font\fontsize{#1}{#2pt}%
  \fontfamily{#3}\fontseries{#4}\fontshape{#5}%
  \selectfont}%
\fi\endgroup%
{\renewcommand{\dashlinestretch}{30}
\begin{picture}(5355,2388)(0,-10)
\put(3000,1002){\makebox(0,0)[lb]{{\SetFigFont{5}{6.0}{\rmdefault}{\mddefault}{\updefault}$\epsilon_3=1$}}}
\put(1275,1752){\blacken\ellipse{36}{36}}
\put(1275,1752){\ellipse{36}{36}}
\put(2175,1752){\blacken\ellipse{36}{36}}
\put(2175,1752){\ellipse{36}{36}}
\put(375,1302){\blacken\ellipse{36}{36}}
\put(375,1302){\ellipse{36}{36}}
\put(1275,1452){\blacken\ellipse{36}{36}}
\put(1275,1452){\ellipse{36}{36}}
\put(2175,1452){\blacken\ellipse{36}{36}}
\put(2175,1452){\ellipse{36}{36}}
\put(2175,1152){\blacken\ellipse{36}{36}}
\put(2175,1152){\ellipse{36}{36}}
\put(1275,1152){\blacken\ellipse{36}{36}}
\put(1275,1152){\ellipse{36}{36}}
\put(1275,852){\blacken\ellipse{36}{36}}
\put(1275,852){\ellipse{36}{36}}
\put(2175,852){\blacken\ellipse{36}{36}}
\put(2175,852){\ellipse{36}{36}}
\put(2175,552){\blacken\ellipse{36}{36}}
\put(2175,552){\ellipse{36}{36}}
\put(1275,552){\blacken\ellipse{36}{36}}
\put(1275,552){\ellipse{36}{36}}
\put(375,702){\blacken\ellipse{36}{36}}
\put(375,702){\ellipse{36}{36}}
\put(3075,1752){\blacken\ellipse{36}{36}}
\put(3075,1752){\ellipse{36}{36}}
\put(3975,1452){\blacken\ellipse{36}{36}}
\put(3975,1452){\ellipse{36}{36}}
\put(4875,1752){\blacken\ellipse{36}{36}}
\put(4875,1752){\ellipse{36}{36}}
\put(4875,1452){\blacken\ellipse{36}{36}}
\put(4875,1452){\ellipse{36}{36}}
\put(4875,1152){\blacken\ellipse{36}{36}}
\put(4875,1152){\ellipse{36}{36}}
\put(3075,1152){\blacken\ellipse{36}{36}}
\put(3075,1152){\ellipse{36}{36}}
\put(3975,852){\blacken\ellipse{36}{36}}
\put(3975,852){\ellipse{36}{36}}
\put(4875,852){\blacken\ellipse{36}{36}}
\put(4875,852){\ellipse{36}{36}}
\put(3075,702){\blacken\ellipse{36}{36}}
\put(3075,702){\ellipse{36}{36}}
\put(3975,552){\blacken\ellipse{36}{36}}
\put(3975,552){\ellipse{36}{36}}
\put(4875,552){\blacken\ellipse{36}{36}}
\put(4875,552){\ellipse{36}{36}}
\path(375,1752)(1275,1752)
\path(1275,1752)(2175,1752)
\path(2175,1452)(1275,1452)(375,1302)
	(1275,1152)(2175,1152)
\path(2175,852)(1275,852)(375,702)
	(1275,552)(2175,552)
\path(4875,852)(3975,852)
\path(4875,1152)(3975,1452)
\path(4875,1752)(3975,1452)(4875,1452)
\path(3075,1752)(3975,1452)(3075,1152)
	(3975,852)(3075,702)(3975,552)(4875,552)
\put(525,1452){\makebox(0,0)[lb]{{\SetFigFont{5}{6.0}{\rmdefault}{\mddefault}{\updefault}$\epsilon_2=2$}}}
\put(525,1077){\makebox(0,0)[lb]{{\SetFigFont{5}{6.0}{\rmdefault}{\mddefault}{\updefault}$\epsilon_3=1$}}}
\put(525,2052){\makebox(0,0)[lb]{{\SetFigFont{8}{9.6}{\rmdefault}{\mddefault}{\updefault}$0\infty$-edges}}}
\put(1425,2052){\makebox(0,0)[lb]{{\SetFigFont{8}{9.6}{\rmdefault}{\mddefault}{\updefault}$\infty t$-edges}}}
\put(4125,2052){\makebox(0,0)[lb]{{\SetFigFont{8}{9.6}{\rmdefault}{\mddefault}{\updefault}$\infty t$-edges}}}
\put(3225,2052){\makebox(0,0)[lb]{{\SetFigFont{8}{9.6}{\rmdefault}{\mddefault}{\updefault}$0\infty$-edges}}}
\put(4650,252){\makebox(0,0)[lb]{{\SetFigFont{5}{6.0}{\rmdefault}{\mddefault}{\updefault}$\al=(1,1,2,2,2)$}}}
\put(2700,477){\makebox(0,0)[lb]{{\SetFigFont{8}{9.6}{\rmdefault}{\mddefault}{\updefault}root}}}
\put(2700,702){\makebox(0,0)[lb]{{\SetFigFont{5}{6.0}{\rmdefault}{\mddefault}{\updefault}$r^0_0=3$}}}
\put(0,702){\makebox(0,0)[lb]{{\SetFigFont{5}{6.0}{\rmdefault}{\mddefault}{\updefault}$r^0_0=3$}}}
\put(0,477){\makebox(0,0)[lb]{{\SetFigFont{8}{9.6}{\rmdefault}{\mddefault}{\updefault}root}}}
\put(3300,477){\makebox(0,0)[lb]{{\SetFigFont{5}{6.0}{\rmdefault}{\mddefault}{\updefault}$\epsilon_5=2$}}}
\put(600,852){\makebox(0,0)[lb]{{\SetFigFont{5}{6.0}{\rmdefault}{\mddefault}{\updefault}$\epsilon_4=1$}}}
\put(600,477){\makebox(0,0)[lb]{{\SetFigFont{5}{6.0}{\rmdefault}{\mddefault}{\updefault}$\epsilon_5=2$}}}
\put(3900,1602){\makebox(0,0)[lb]{{\SetFigFont{5}{6.0}{\rmdefault}{\mddefault}{\updefault}$r^1_{\infty}=3$}}}
\put(3825,402){\makebox(0,0)[lb]{{\SetFigFont{5}{6.0}{\rmdefault}{\mddefault}{\updefault}$r^3_{\infty}=0$}}}
\put(3225,777){\makebox(0,0)[lb]{{\SetFigFont{5}{6.0}{\rmdefault}{\mddefault}{\updefault}$\epsilon_4=1$}}}
\put(3900,927){\makebox(0,0)[lb]{{\SetFigFont{5}{6.0}{\rmdefault}{\mddefault}{\updefault}$r^2_{\infty}=1$}}}
\put(3900,177){\makebox(0,0)[lb]{{\SetFigFont{5}{6.0}{\rmdefault}{\mddefault}{\updefault}$r_{\infty}=4$}}}
\put(225,2277){\makebox(0,0)[lb]{{\SetFigFont{8}{9.6}{\rmdefault}{\mddefault}{\updefault}$0$-vertices}}}
\put(1125,2277){\makebox(0,0)[lb]{{\SetFigFont{8}{9.6}{\rmdefault}{\mddefault}{\updefault}$\infty$-vertices}}}
\put(2025,2277){\makebox(0,0)[lb]{{\SetFigFont{8}{9.6}{\rmdefault}{\mddefault}{\updefault}$t$-vertices}}}
\put(2925,2277){\makebox(0,0)[lb]{{\SetFigFont{8}{9.6}{\rmdefault}{\mddefault}{\updefault}$0$-vertices}}}
\put(3825,2277){\makebox(0,0)[lb]{{\SetFigFont{8}{9.6}{\rmdefault}{\mddefault}{\updefault}$\infty$-vertices}}}
\put(4725,2277){\makebox(0,0)[lb]{{\SetFigFont{8}{9.6}{\rmdefault}{\mddefault}{\updefault}$t$-vertices}}}
\put(0,1827){\makebox(0,0)[lb]{{\SetFigFont{5}{6.0}{\rmdefault}{\mddefault}{\updefault}$\delta^1=(2)$}}}
\put(0,1677){\makebox(0,0)[lb]{{\SetFigFont{5}{6.0}{\rmdefault}{\mddefault}{\updefault}$r^1_0=1$}}}
\put(0,1377){\makebox(0,0)[lb]{{\SetFigFont{5}{6.0}{\rmdefault}{\mddefault}{\updefault}$\delta^2=(1,2)$}}}
\put(0,1227){\makebox(0,0)[lb]{{\SetFigFont{5}{6.0}{\rmdefault}{\mddefault}{\updefault}$r^2_0=3$}}}
\put(0,852){\makebox(0,0)[lb]{{\SetFigFont{5}{6.0}{\rmdefault}{\mddefault}{\updefault}$\delta^0=(1,2)$}}}
\put(2700,852){\makebox(0,0)[lb]{{\SetFigFont{5}{6.0}{\rmdefault}{\mddefault}{\updefault}$\delta^0=(1,2)$}}}
\put(2700,1077){\makebox(0,0)[lb]{{\SetFigFont{5}{6.0}{\rmdefault}{\mddefault}{\updefault}$r^2_0=3$}}}
\put(2700,1227){\makebox(0,0)[lb]{{\SetFigFont{5}{6.0}{\rmdefault}{\mddefault}{\updefault}$\delta^2=(1,2)$}}}
\put(2700,1677){\makebox(0,0)[lb]{{\SetFigFont{5}{6.0}{\rmdefault}{\mddefault}{\updefault}$r^1_0=1$}}}
\put(2700,1827){\makebox(0,0)[lb]{{\SetFigFont{5}{6.0}{\rmdefault}{\mddefault}{\updefault}$\delta^1=(2)$}}}
\put(975,27){\makebox(0,0)[lb]{{\SetFigFont{5}{6.0}{\rmdefault}{\mddefault}{\updefault}$s_{\infty}=-5$}}}
\put(975,177){\makebox(0,0)[lb]{{\SetFigFont{5}{6.0}{\rmdefault}{\mddefault}{\updefault}$r_{\infty}=r^i_{\infty}=0$}}}
\put(1875,252){\makebox(0,0)[lb]{{\SetFigFont{5}{6.0}{\rmdefault}{\mddefault}{\updefault}$\al=(1,1,2,2,2)$}}}
\put(525,1827){\makebox(0,0)[lb]{{\SetFigFont{5}{6.0}{\rmdefault}{\mddefault}{\updefault}$\epsilon_1=2$}}}
\put(3450,1752){\makebox(0,0)[lb]{{\SetFigFont{5}{6.0}{\rmdefault}{\mddefault}{\updefault}$\be^1=(2,2)$}}}
\put(4050,1752){\makebox(0,0)[lb]{{\SetFigFont{5}{6.0}{\rmdefault}{\mddefault}{\updefault}$\ga^1=(1,1,2)$}}}
\put(4725,1227){\makebox(0,0)[lb]{{\SetFigFont{5}{6.0}{\rmdefault}{\mddefault}{\updefault}$1$}}}
\put(4200,1077){\makebox(0,0)[lb]{{\SetFigFont{5}{6.0}{\rmdefault}{\mddefault}{\updefault}$\ga^2=(2)$}}}
\put(4725,927){\makebox(0,0)[lb]{{\SetFigFont{5}{6.0}{\rmdefault}{\mddefault}{\updefault}$2$}}}
\put(4725,627){\makebox(0,0)[lb]{{\SetFigFont{5}{6.0}{\rmdefault}{\mddefault}{\updefault}$2$}}}
\put(2025,1527){\makebox(0,0)[lb]{{\SetFigFont{5}{6.0}{\rmdefault}{\mddefault}{\updefault}$2$}}}
\put(2025,1227){\makebox(0,0)[lb]{{\SetFigFont{5}{6.0}{\rmdefault}{\mddefault}{\updefault}$1$}}}
\put(2025,927){\makebox(0,0)[lb]{{\SetFigFont{5}{6.0}{\rmdefault}{\mddefault}{\updefault}$1$}}}
\put(2025,627){\makebox(0,0)[lb]{{\SetFigFont{5}{6.0}{\rmdefault}{\mddefault}{\updefault}$2$}}}
\put(2025,1827){\makebox(0,0)[lb]{{\SetFigFont{5}{6.0}{\rmdefault}{\mddefault}{\updefault}$2$}}}
\put(3600,1077){\makebox(0,0)[lb]{{\SetFigFont{5}{6.0}{\rmdefault}{\mddefault}{\updefault}$\be^2=(1,1)$}}}
\put(4725,1527){\makebox(0,0)[lb]{{\SetFigFont{5}{6.0}{\rmdefault}{\mddefault}{\updefault}$2$}}}
\put(4725,1827){\makebox(0,0)[lb]{{\SetFigFont{5}{6.0}{\rmdefault}{\mddefault}{\updefault}$1$}}}
\put(1050,1827){\makebox(0,0)[lb]{{\SetFigFont{5}{6.0}{\rmdefault}{\mddefault}{\updefault}$\be^1=(2)=\ga^1$}}}
\put(1050,1527){\makebox(0,0)[lb]{{\SetFigFont{5}{6.0}{\rmdefault}{\mddefault}{\updefault}$\be^2=(2)=\ga^2$}}}
\put(1050,1227){\makebox(0,0)[lb]{{\SetFigFont{5}{6.0}{\rmdefault}{\mddefault}{\updefault}$\be^3=(1)=\ga^3$}}}
\put(1050,927){\makebox(0,0)[lb]{{\SetFigFont{5}{6.0}{\rmdefault}{\mddefault}{\updefault}$\be^4=(1)=\ga^4$}}}
\put(1050,627){\makebox(0,0)[lb]{{\SetFigFont{5}{6.0}{\rmdefault}{\mddefault}{\updefault}$\be^5=(2)=\ga^5$}}}
\put(3750,627){\makebox(0,0)[lb]{{\SetFigFont{5}{6.0}{\rmdefault}{\mddefault}{\updefault}$\be^3=(2)=\ga^3$}}}
\put(3900,27){\makebox(0,0)[lb]{{\SetFigFont{5}{6.0}{\rmdefault}{\mddefault}{\updefault}$s_{\infty}=-1$}}}
\put(3150,1527){\makebox(0,0)[lb]{{\SetFigFont{5}{6.0}{\rmdefault}{\mddefault}{\updefault}$\epsilon_1=2$}}}
\put(3525,1227){\makebox(0,0)[lb]{{\SetFigFont{5}{6.0}{\rmdefault}{\mddefault}{\updefault}$\epsilon_2=2$}}}
\put(375,1752){\blacken\ellipse{36}{36}}
\put(375,1752){\ellipse{36}{36}}
\end{picture}
}

\end{center}
\caption{Sample decorated trees (Def.~\ref{Dloctr}) corresponding roughly to the examples
  of Figure~\ref{twoexamples}.
Omitted are further combinatorial
  labelling of vertices.
  \lremind{treehouse}}
\label{treehouse}
\end{figure}

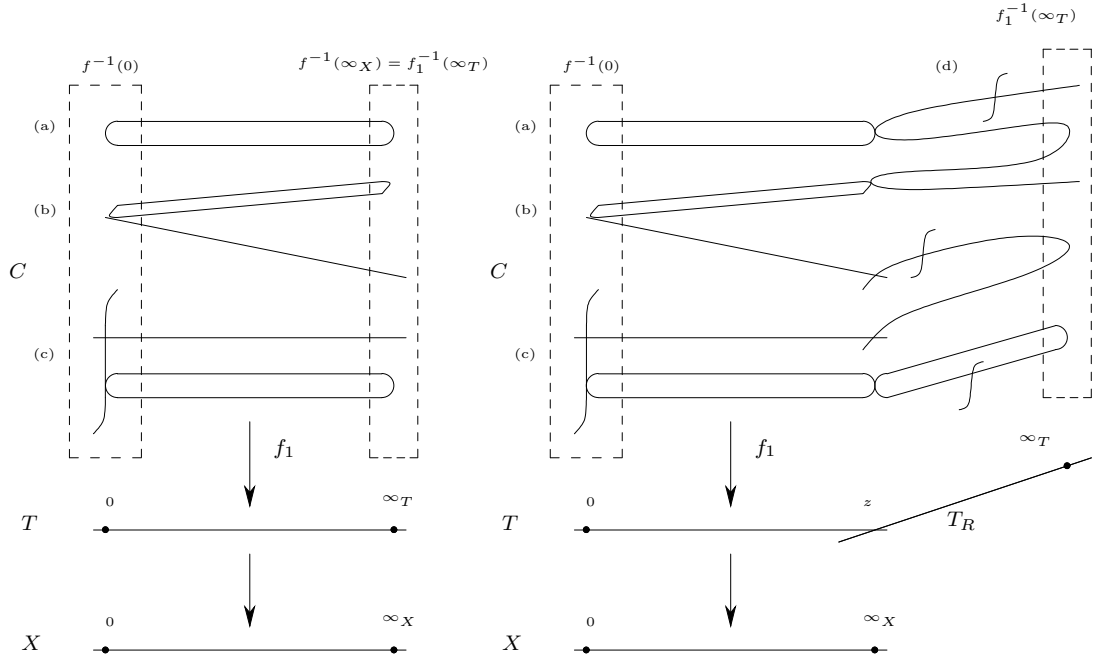
\begin{figure}[ht]
\begin{center}
\setlength{\unitlength}{0.00083333in}
\begingroup\makeatletter\ifx\SetFigFont\undefined%
\gdef\SetFigFont#1#2#3#4#5{%
  \reset@font\fontsize{#1}{#2pt}%
  \fontfamily{#3}\fontseries{#4}\fontshape{#5}%
  \selectfont}%
\fi\endgroup%
{\renewcommand{\dashlinestretch}{30}
\begin{picture}(6762,4018)(0,-10)
\put(4650,1231){\makebox(0,0)[lb]{{\SetFigFont{8}{9.6}{\rmdefault}{\mddefault}{\updefault}$f_1$}}}
\put(5325.000,3256.000){\arc{150.000}{4.7124}{7.8540}}
\put(3675.000,1681.000){\arc{150.000}{1.5708}{4.7124}}
\put(5325.000,1681.000){\arc{150.000}{4.7124}{7.8540}}
\put(5475.000,1681.000){\arc{150.000}{1.5708}{4.7124}}
\put(6525.000,1981.000){\arc{150.000}{4.7124}{7.8540}}
\put(675.000,3256.000){\arc{150.000}{1.5708}{4.7124}}
\put(2325.000,3256.000){\arc{150.000}{4.7124}{7.8540}}
\put(675.000,1681.000){\arc{150.000}{1.5708}{4.7124}}
\put(2325.000,1681.000){\arc{150.000}{4.7124}{7.8540}}
\put(3600,31){\blacken\ellipse{36}{36}}
\put(3600,31){\ellipse{36}{36}}
\put(5400,31){\blacken\ellipse{36}{36}}
\put(5400,31){\ellipse{36}{36}}
\put(3600,781){\blacken\ellipse{36}{36}}
\put(3600,781){\ellipse{36}{36}}
\put(6600,1181){\blacken\ellipse{36}{36}}
\put(6600,1181){\ellipse{36}{36}}
\put(600,31){\blacken\ellipse{36}{36}}
\put(600,31){\ellipse{36}{36}}
\put(2400,31){\blacken\ellipse{36}{36}}
\put(2400,31){\ellipse{36}{36}}
\put(600,781){\blacken\ellipse{36}{36}}
\put(600,781){\ellipse{36}{36}}
\put(2400,781){\blacken\ellipse{36}{36}}
\put(2400,781){\ellipse{36}{36}}
\path(3675,2806)(5325,2956)
\path(3675,2731)(5325,2881)
\path(3600,2731)(5475,2356)
\path(3675,3181)(5325,3181)
\path(3675,3331)(5325,3331)
\path(5475,1981)(3525,1981)
\path(3675,1756)(5325,1756)
\path(5325,1606)(3675,1606)
\path(3600,1606)(3600,2056)
\path(4500,1456)(4500,931)
\blacken\path(4470.000,1051.000)(4500.000,931.000)(4530.000,1051.000)(4500.000,1015.000)(4470.000,1051.000)
\path(4500,631)(4500,181)
\blacken\path(4470.000,301.000)(4500.000,181.000)(4530.000,301.000)(4500.000,265.000)(4470.000,301.000)
\path(3525,781)(5475,781)
\path(3525,31)(5475,31)
\dashline{60.000}(3375,3556)(3825,3556)(3825,1231)
	(3375,1231)(3375,3556)
\path(5175,706)(6750,1231)
\path(5175,706)(6750,1231)
\path(5475,1756)(6525,2056)
\path(5475,1606)(6525,1906)
\dashline{60.000}(6450,3781)(6750,3781)(6750,1606)
	(6450,1606)(6450,3781)
\path(675,2806)(2325,2956)
\path(675,2731)(2325,2881)
\path(600,2731)(2475,2356)
\path(675,3181)(2325,3181)
\path(675,3331)(2325,3331)
\path(2475,1981)(525,1981)
\path(675,1756)(2325,1756)
\path(2325,1606)(675,1606)
\path(600,1606)(600,2056)
\path(1500,1456)(1500,931)
\blacken\path(1470.000,1051.000)(1500.000,931.000)(1530.000,1051.000)(1500.000,1015.000)(1470.000,1051.000)
\path(1500,631)(1500,181)
\blacken\path(1470.000,301.000)(1500.000,181.000)(1530.000,301.000)(1500.000,265.000)(1470.000,301.000)
\path(525,781)(2475,781)
\path(525,31)(2475,31)
\dashline{60.000}(375,3556)(825,3556)(825,1231)
	(375,1231)(375,3556)
\dashline{60.000}(2250,3556)(2550,3556)(2550,1231)
	(2250,1231)(2250,3556)
\path(3675,2806)(3674,2805)(3668,2799)
	(3658,2788)(3647,2775)(3637,2764)
	(3630,2755)(3626,2748)(3625,2743)
	(3627,2739)(3634,2735)(3647,2733)
	(3662,2732)(3673,2731)(3675,2731)
\path(5325,2956)(5327,2956)(5338,2955)
	(5353,2954)(5366,2952)(5373,2948)
	(5375,2943)(5374,2939)(5370,2932)
	(5363,2923)(5353,2912)(5342,2899)
	(5332,2888)(5326,2882)(5325,2881)
\path(3600,2056)(3600,2057)(3600,2064)
	(3601,2078)(3601,2098)(3602,2119)
	(3603,2139)(3605,2156)(3607,2171)
	(3609,2183)(3613,2193)(3616,2204)
	(3621,2213)(3627,2223)(3635,2234)
	(3644,2246)(3655,2258)(3664,2269)
	(3671,2277)(3674,2280)(3675,2281)
\path(3600,1606)(3600,1605)(3600,1598)
	(3599,1584)(3599,1564)(3598,1543)
	(3597,1523)(3595,1506)(3593,1491)
	(3591,1479)(3588,1468)(3584,1458)
	(3579,1449)(3573,1439)(3565,1428)
	(3556,1416)(3545,1404)(3536,1393)
	(3529,1385)(3526,1382)(3525,1381)
\path(6675,3556)(6672,3556)(6667,3555)
	(6656,3553)(6640,3551)(6618,3548)
	(6589,3544)(6553,3539)(6513,3534)
	(6467,3527)(6419,3520)(6367,3513)
	(6315,3506)(6262,3498)(6211,3491)
	(6160,3483)(6112,3476)(6066,3469)
	(6023,3463)(5983,3457)(5945,3451)
	(5910,3445)(5877,3440)(5846,3434)
	(5818,3429)(5792,3425)(5767,3420)
	(5743,3415)(5721,3411)(5700,3406)
	(5667,3398)(5636,3391)(5606,3383)
	(5579,3375)(5553,3366)(5529,3358)
	(5507,3350)(5487,3341)(5469,3333)
	(5453,3324)(5439,3316)(5428,3308)
	(5419,3301)(5412,3294)(5406,3287)
	(5403,3280)(5401,3274)(5400,3268)
	(5401,3263)(5403,3257)(5406,3252)
	(5411,3247)(5418,3242)(5426,3237)
	(5437,3233)(5449,3229)(5463,3225)
	(5480,3222)(5497,3220)(5517,3218)
	(5538,3216)(5560,3216)(5584,3215)
	(5608,3216)(5635,3217)(5663,3218)
	(5683,3220)(5704,3222)(5726,3224)
	(5750,3226)(5774,3229)(5800,3232)
	(5826,3235)(5854,3238)(5882,3242)
	(5911,3246)(5940,3249)(5970,3254)
	(6000,3258)(6030,3262)(6060,3266)
	(6089,3270)(6118,3275)(6146,3279)
	(6174,3283)(6200,3287)(6226,3290)
	(6250,3294)(6274,3297)(6296,3300)
	(6317,3303)(6338,3306)(6365,3310)
	(6392,3313)(6416,3316)(6440,3318)
	(6463,3320)(6484,3321)(6503,3321)
	(6522,3321)(6538,3321)(6553,3319)
	(6566,3317)(6578,3314)(6587,3310)
	(6595,3306)(6602,3301)(6606,3295)
	(6610,3288)(6613,3281)(6614,3274)
	(6614,3266)(6613,3257)(6612,3248)
	(6609,3238)(6604,3227)(6598,3216)
	(6591,3205)(6581,3193)(6571,3181)
	(6558,3169)(6544,3157)(6528,3146)
	(6511,3135)(6492,3124)(6471,3114)
	(6450,3105)(6426,3096)(6402,3088)
	(6375,3081)(6353,3076)(6330,3071)
	(6306,3067)(6280,3063)(6253,3059)
	(6225,3055)(6195,3052)(6165,3049)
	(6133,3046)(6101,3043)(6067,3041)
	(6034,3038)(6000,3036)(5966,3034)
	(5933,3033)(5899,3031)(5867,3029)
	(5835,3028)(5805,3027)(5775,3026)
	(5747,3024)(5720,3023)(5694,3022)
	(5670,3021)(5647,3020)(5625,3018)
	(5596,3016)(5568,3014)(5543,3012)
	(5519,3009)(5497,3006)(5477,3003)
	(5459,2999)(5442,2996)(5428,2992)
	(5415,2988)(5404,2984)(5396,2980)
	(5389,2976)(5383,2972)(5379,2968)
	(5377,2964)(5375,2960)(5375,2956)
	(5375,2952)(5377,2948)(5380,2944)
	(5384,2941)(5389,2937)(5397,2933)
	(5406,2929)(5418,2925)(5431,2922)
	(5447,2919)(5465,2916)(5485,2913)
	(5507,2911)(5532,2909)(5558,2908)
	(5586,2907)(5617,2906)(5650,2906)
	(5670,2906)(5691,2906)(5713,2907)
	(5736,2907)(5761,2908)(5788,2909)
	(5817,2910)(5848,2911)(5881,2912)
	(5917,2914)(5955,2916)(5995,2918)
	(6039,2920)(6084,2922)(6132,2925)
	(6182,2928)(6234,2930)(6286,2933)
	(6338,2936)(6389,2939)(6439,2942)
	(6485,2945)(6528,2947)(6565,2949)
	(6597,2951)(6624,2953)(6644,2954)
	(6658,2955)(6668,2956)(6673,2956)(6675,2956)
\path(5325,2281)(5326,2283)(5330,2288)
	(5336,2296)(5344,2306)(5355,2319)
	(5367,2334)(5381,2349)(5396,2363)
	(5412,2376)(5429,2389)(5447,2401)
	(5466,2412)(5487,2422)(5511,2433)
	(5538,2443)(5555,2450)(5574,2457)
	(5594,2464)(5616,2471)(5638,2478)
	(5662,2486)(5687,2493)(5714,2501)
	(5741,2508)(5769,2516)(5799,2523)
	(5828,2531)(5859,2538)(5889,2546)
	(5920,2553)(5950,2559)(5981,2566)
	(6011,2572)(6040,2577)(6069,2583)
	(6097,2588)(6124,2592)(6151,2596)
	(6176,2600)(6201,2603)(6225,2606)
	(6259,2609)(6292,2612)(6324,2614)
	(6355,2615)(6386,2615)(6415,2615)
	(6444,2613)(6471,2611)(6496,2608)
	(6518,2604)(6539,2600)(6557,2595)
	(6573,2590)(6586,2584)(6596,2577)
	(6604,2571)(6609,2563)(6613,2556)
	(6613,2548)(6612,2540)(6608,2531)
	(6601,2521)(6592,2511)(6581,2500)
	(6567,2489)(6550,2477)(6531,2465)
	(6510,2453)(6487,2440)(6462,2428)
	(6436,2415)(6408,2403)(6380,2391)
	(6350,2379)(6319,2367)(6288,2356)
	(6265,2348)(6241,2340)(6217,2332)
	(6191,2324)(6165,2315)(6138,2307)
	(6110,2298)(6081,2289)(6051,2280)
	(6022,2271)(5991,2262)(5961,2253)
	(5930,2243)(5900,2234)(5870,2225)
	(5840,2216)(5812,2207)(5784,2198)
	(5757,2189)(5731,2180)(5706,2172)
	(5683,2163)(5660,2155)(5639,2147)
	(5619,2139)(5600,2131)(5574,2120)
	(5551,2108)(5528,2096)(5508,2083)
	(5488,2070)(5468,2055)(5449,2039)
	(5431,2022)(5412,2004)(5394,1985)
	(5377,1967)(5361,1949)(5348,1934)
	(5338,1922)(5331,1913)(5327,1908)(5325,1906)
\path(6075,1831)(6074,1831)(6071,1830)
	(6064,1829)(6055,1826)(6044,1822)
	(6035,1817)(6027,1811)(6021,1803)
	(6016,1793)(6013,1781)(6010,1770)
	(6008,1759)(6006,1745)(6005,1730)
	(6003,1715)(6001,1698)(6000,1681)
	(5999,1664)(5997,1647)(5995,1632)
	(5994,1617)(5992,1603)(5990,1592)
	(5988,1581)(5984,1569)(5979,1559)
	(5973,1551)(5965,1545)(5956,1540)
	(5945,1536)(5936,1533)(5929,1532)
	(5926,1531)(5925,1531)
\path(5775,2656)(5774,2656)(5771,2655)
	(5764,2654)(5755,2651)(5744,2647)
	(5735,2642)(5727,2636)(5721,2628)
	(5716,2618)(5713,2606)(5710,2595)
	(5708,2584)(5706,2570)(5705,2555)
	(5703,2540)(5701,2523)(5700,2506)
	(5699,2489)(5697,2472)(5695,2457)
	(5694,2442)(5692,2428)(5690,2417)
	(5688,2406)(5684,2394)(5679,2384)
	(5673,2376)(5665,2370)(5656,2365)
	(5645,2361)(5636,2358)(5629,2357)
	(5626,2356)(5625,2356)
\path(6225,3631)(6224,3631)(6221,3630)
	(6214,3629)(6205,3626)(6194,3622)
	(6185,3617)(6177,3611)(6171,3603)
	(6166,3593)(6163,3581)(6160,3570)
	(6158,3559)(6156,3545)(6155,3530)
	(6153,3515)(6151,3498)(6150,3481)
	(6149,3464)(6147,3447)(6145,3432)
	(6144,3417)(6142,3403)(6140,3392)
	(6138,3381)(6134,3369)(6129,3359)
	(6123,3351)(6115,3345)(6106,3340)
	(6095,3336)(6086,3333)(6079,3332)
	(6076,3331)(6075,3331)
\path(675,2806)(674,2805)(668,2799)
	(658,2788)(647,2775)(637,2764)
	(630,2755)(626,2748)(625,2743)
	(627,2739)(634,2735)(647,2733)
	(662,2732)(673,2731)(675,2731)
\path(2325,2956)(2327,2956)(2338,2955)
	(2353,2954)(2366,2952)(2373,2948)
	(2375,2943)(2374,2939)(2370,2932)
	(2363,2923)(2353,2912)(2342,2899)
	(2332,2888)(2326,2882)(2325,2881)
\path(600,2056)(600,2057)(600,2064)
	(601,2078)(601,2098)(602,2119)
	(603,2139)(605,2156)(607,2171)
	(609,2183)(613,2193)(616,2204)
	(621,2213)(627,2223)(635,2234)
	(644,2246)(655,2258)(664,2269)
	(671,2277)(674,2280)(675,2281)
\path(600,1606)(600,1605)(600,1598)
	(599,1584)(599,1564)(598,1543)
	(597,1523)(595,1506)(593,1491)
	(591,1479)(588,1468)(584,1458)
	(579,1449)(573,1439)(565,1428)
	(556,1416)(545,1404)(536,1393)
	(529,1385)(526,1382)(525,1381)
\put(3600,931){\makebox(0,0)[lb]{{\SetFigFont{5}{6.0}{\rmdefault}{\mddefault}{\updefault}$0$}}}
\put(3600,181){\makebox(0,0)[lb]{{\SetFigFont{5}{6.0}{\rmdefault}{\mddefault}{\updefault}$0$}}}
\put(5325,181){\makebox(0,0)[lb]{{\SetFigFont{5}{6.0}{\rmdefault}{\mddefault}{\updefault}$\infty_X$}}}
\put(3150,3256){\makebox(0,0)[lb]{{\SetFigFont{5}{6.0}{\rmdefault}{\mddefault}{\updefault}(a)}}}
\put(3150,2731){\makebox(0,0)[lb]{{\SetFigFont{5}{6.0}{\rmdefault}{\mddefault}{\updefault}(b)}}}
\put(3150,1831){\makebox(0,0)[lb]{{\SetFigFont{5}{6.0}{\rmdefault}{\mddefault}{\updefault}(c)}}}
\put(3450,3631){\makebox(0,0)[lb]{{\SetFigFont{5}{6.0}{\rmdefault}{\mddefault}{\updefault}$f^{-1}(0)$}}}
\put(3000,2356){\makebox(0,0)[lb]{{\SetFigFont{8}{9.6}{\rmdefault}{\mddefault}{\updefault}$C$}}}
\put(3075,781){\makebox(0,0)[lb]{{\SetFigFont{8}{9.6}{\rmdefault}{\mddefault}{\updefault}$T$}}}
\put(3075,31){\makebox(0,0)[lb]{{\SetFigFont{8}{9.6}{\rmdefault}{\mddefault}{\updefault}$X$}}}
\put(6300,1306){\makebox(0,0)[lb]{{\SetFigFont{5}{6.0}{\rmdefault}{\mddefault}{\updefault}$\infty_T$}}}
\put(6150,3931){\makebox(0,0)[lb]{{\SetFigFont{5}{6.0}{\rmdefault}{\mddefault}{\updefault}$f_1^{-1}(\infty_T)$}}}
\put(5775,3631){\makebox(0,0)[lb]{{\SetFigFont{5}{6.0}{\rmdefault}{\mddefault}{\updefault}(d)}}}
\put(5850,781){\makebox(0,0)[lb]{{\SetFigFont{8}{9.6}{\rmdefault}{\mddefault}{\updefault}$T_R$}}}
\put(5325,931){\makebox(0,0)[lb]{{\SetFigFont{5}{6.0}{\rmdefault}{\mddefault}{\updefault}$z$}}}
\put(1800,3631){\makebox(0,0)[lb]{{\SetFigFont{5}{6.0}{\rmdefault}{\mddefault}{\updefault}$f^{-1}(\infty_X) = f_1^{-1}(\infty_T)$}}}
\put(600,931){\makebox(0,0)[lb]{{\SetFigFont{5}{6.0}{\rmdefault}{\mddefault}{\updefault}$0$}}}
\put(2325,931){\makebox(0,0)[lb]{{\SetFigFont{5}{6.0}{\rmdefault}{\mddefault}{\updefault}$\infty_T$}}}
\put(600,181){\makebox(0,0)[lb]{{\SetFigFont{5}{6.0}{\rmdefault}{\mddefault}{\updefault}$0$}}}
\put(2325,181){\makebox(0,0)[lb]{{\SetFigFont{5}{6.0}{\rmdefault}{\mddefault}{\updefault}$\infty_X$}}}
\put(150,3256){\makebox(0,0)[lb]{{\SetFigFont{5}{6.0}{\rmdefault}{\mddefault}{\updefault}(a)}}}
\put(150,2731){\makebox(0,0)[lb]{{\SetFigFont{5}{6.0}{\rmdefault}{\mddefault}{\updefault}(b)}}}
\put(150,1831){\makebox(0,0)[lb]{{\SetFigFont{5}{6.0}{\rmdefault}{\mddefault}{\updefault}(c)}}}
\put(450,3631){\makebox(0,0)[lb]{{\SetFigFont{5}{6.0}{\rmdefault}{\mddefault}{\updefault}$f^{-1}(0)$}}}
\put(0,2356){\makebox(0,0)[lb]{{\SetFigFont{8}{9.6}{\rmdefault}{\mddefault}{\updefault}$C$}}}
\put(75,781){\makebox(0,0)[lb]{{\SetFigFont{8}{9.6}{\rmdefault}{\mddefault}{\updefault}$T$}}}
\put(75,31){\makebox(0,0)[lb]{{\SetFigFont{8}{9.6}{\rmdefault}{\mddefault}{\updefault}$X$}}}
\put(1650,1231){\makebox(0,0)[lb]{{\SetFigFont{8}{9.6}{\rmdefault}{\mddefault}{\updefault}$f_1$}}}
\put(3675.000,3256.000){\arc{150.000}{1.5708}{4.7124}}
\end{picture}
}

\end{center}
\caption{Sketches of two examples of torus-fixed relative stable maps to $X=(\proj^1, \infty)$. 
$T$ is the ``sprouted target''.
 \lremind{twoexamples}}
\label{twoexamples}
\end{figure}

Hereinafter for any localization tree $\rt$ we shall
subsume the dependence on $\rt$ of  the above sets of vertices and edges, and 
the partitions, by suppressing the occurrence of $\rt$ as an argument. 
We also define some more notation. Let
\begin{equation}\label{rrr}
r_{\infty}:=\sum_{v \in\sV_\infty}r^0_{\ga^v,\be^v} ,
\end{equation}
and let
\begin{equation}\label{spal}
\al := \coprod_{v\in\sV_\infty}\gamma^v
\end{equation} 
\dmjdel{
If $v$  is an $\infty$-vertex, define
\begin{equation}\label{rrr}
r^v_{\infty} := l(\be^v)+l(\ga^v)-2 \;\text{and}\; r_{\infty}:=\sum_{v \in \sV_{\infty}} r^v_{\infty}.
\end{equation}
If $v$ is a $0$-vertex, define
\begin{equation}
\label{rrr2}
r^v_0:=l(\de^v)+|\de^v|-2.
\end{equation} 
We also let 
\begin{equation}\label{spal}
\al:=\coprod_{v\in\sV_\infty}\gamma^v,
\end{equation} 
} 
be the partition formed by all the parts of each of the
$\ga^v$'s for non-root $\infty$-vertices $v$, and let $d=|\al |$.  Examples of such trees
{\em without a specified root} are given in Figure~\ref{treehouse} (these
examples correspond to the geometric picture of
Figure~\ref{twoexamples}).

\newcommand{\nsp}{\!\!\!}

We now define some terminology that will allow us
to write the relative virtual localization calculation cleanly.
Let
$\P^g_m(\al_1, \dots, \al_m)$ be the
dimension $2g-1$
portion of\label{fpoly}\lremind{fpoly}
\begin{equation*}
\frac { 1 - \la_1 + \cdots + (-1)^g \la_g}
{ (1- \al_1 \psi_1) \cdots (1 - \al_m \psi_m)}
\end{equation*}
on $\cm_{g,m}^{rt}$, \emph{viz.},
\begin{equation}\label{Fabpoly2}
\P^g_m(\al_1,\ld ,\al_m ) := 
\sum_{\substack{a_1,\ldots,a_m, k\ge0, \\ a_1+\cdots+a_m+k=g-2+m}}
(-1)^k \laFab\tau_{a_1}\cd\tau_{a_m} \la_k \raFab \al_1^{a_1}\cdots\al_m^{a_m}.
\end{equation}
Note that $\P^g_m$ is a
(Chow-valued) polynomial in the numbers $\al_1$, \dots, $\al_m$,
symmetric of degree between $m-2$ and $g-2+m$, and its leading
coefficients (the portion of homogeneous degree $g-2+m$) are precisely
the subject of Faber's Intersection Number Conjecture.  (It is easy to
show that the homogeneous degree $m-2$ and $m-1$ portions of this
polynomial vanish, but we shall not need this fact.)
We shall refer to~$\P^g_m$ as the \emph{Faber polynomial}.

For any localization tree $\rt$, let
\begin{equation}\label{ABCD}
B(\rt) :=\nsp\prod_{e\in\sE_{0\infty}}\nsp \ep( e),\,\,\,
C(\rt) :=\nsp\prod_{v\in\sV_0} \nsp \frac{H^0_{\de^v}}{r^0_{\de^v}!},\,\,\,
D(\rt) :=\prod_{v\in\sV_\infty}\frac{H^0_{\ga^v,\be^v}}{r^0_{\ga^v,\be^v}!}.
\end{equation}
and, for $\rt\in\sT_{g,m}$, let
\begin{equation}\label{Aonly}
A(\rt) :=\P^g_m(\de^\bullet )\prod
\frac{\ep(e)^{\ep(e)}}{\ep(e)!},
\end{equation}
where the product is over all edges $e$ incident with the root-vertex $\bullet$ of $\rt$.
Let ``$\,\,\mbox{}^\ddagger\,\,$'' as a superscript on a product denote the removal of
the contribution of the root-vertex from that product.

\begin{tm}[Localization Tree Theorem --- tree summation]\label{bigyuck}
For $g\ge 1$ and $\al=(\al_1,\ldots,\al_m)\vdash d,$
\begin{equation}\label{RaviR}
\F^{g,\al}=\sum_{m\ge 1}\sum_{\rt\in\sT_{g,m}} (-1)^{r_{\infty}}
r^{\Faber}_{g,\al}!
\binom {  r^g_{\al} - r_{\infty}}   {r^{\Faber}_{g,\al}}
\frac{1}{\eta_0(\rt)!} 
A(\rt) B(\rt) C^\ddagger(\rt) D(\rt)
\end{equation}
where the sum is subject to~(\ref{spal}). \end{tm}

\noindent Note that the sum in~\eqref{RaviR} is finite, and that
the binomial coefficient in it is
zero unless $r_{\infty}$ is small (at most
$r^g_{\al} -r^{\Faber}_{g,\al} = 2g-1$).  

Before proving Theorem~\ref{bigyuck}, we digress to observe that just the ``shape'' of the formula
\eqref{RaviR} quickly yields the result   {\bf (I)} promised  in Section~\ref{sor}.

\begin{tm} [Socle statement for $\cm^{rt}_{g,n}$]  $R_{2g-1}(\cm^{rt}_{g,n}) \cong \Q$.
\label{onedim}\lremind{onedim}
\end{tm}

\bpf
 We shall show that any monomial in the $\psi$-classes (pushed
forward to $\cm_{g,1}$) is a multiple of $\G_{g,1}$.  This (with 
Faber's Non-vanishing Theorem~\ref{nv} and Remark~\ref{summary} (iii)) implies that
$R_{2g-1}(\cm^{rt}_{g,n})$ is generated by a single element.
{\em In particular, $\G_{g,1} \neq 0$ and
\begin{equation}\label{divide}
\xymatrix{  \Q \ar[rr]^{\times \G_{g,1}} & & R_{2g-1}(\cm_{g,1})
}
\end{equation}
is an isomorphism.
}

We show first that $\P^g_n(\al)$ is a multiple of $\G_{g,1}$ for each
partition $\al$, by induction on $| \al | = \sum_i \al_i$.  Now
$\F^{g,\al}$ is a multiple of $\G_{g,1}$ (Thm.~\ref{connec}), and the
contribution of the unique graph with $r_{\infty}=0$ (``the simplest
graph in $\sT_{g,m}$'') to Theorem~\ref{bigyuck} is a non-zero multiple
of $\P^g_n(\al)$.  The contribution of any other graph is a multiple
of $\P^g_{n'}(\al')$ for some smaller $| \al' |$.  Thus, by induction,
$\P^g_n(\al)$ is a multiple of $\G_{g,1}$ as desired.

Next, we apply the ``polynomiality trick'' used in \cite{socle}
and \cite{thmstar}.  We fix $g$ and $n$, and hence the polynomial
$\P^g_n(\cdot)$.
For arbitrary choices of $\al_1$, \dots, $\al_n$, $\P^g_n(\al)$ is a multiple
of $\G_{g,1}$.    But by knowing enough values of  a polynomial
of known degree,
we can determine  its coefficients as linear combinations of these values.
Hence all coefficients of $\P^g_n(\al)$ are multiples
of $\G_{g,1}$, and in particular, the monomials in $\psi$-classes
(of  degree $g-2+n$) are multiples of $\G_{g,1}$. 
\epf

\noindent {\em Proof of Theorem~\ref{bigyuck}.}
We apply relative virtual localization as
developed in \cite{thmstar} (based on the 
foundational \cite{vl}).   
As with many virtual localization calculations, Faber classes will be
expressed as sums over certain graphs.
We show that the set of trees $\cT_{g,m}$ is precisely the set of graphs that
is required for this purpose by supplying a geometric meaning to
the vertices, edges, partitions, weights and constants associated with $\rt\in\sT_{g,m}$
through a geometry-combinatorics lexicon.

\noindent\underline{\emph{Classification of torus-fixed loci:}}
We shall have a contribution from each torus-fixed locus of stable
relative maps.  Recall that Faber-Hurwitz classes are defined by
considering the pullback of a linear space under the branch map,
applied to the virtual fundamental class of the moduli space of stable
relative maps (equ.\ \eqref{fhclass}).  As in the proof of the 
``tautological vanishing theorem'' of \cite{thmstar}, we choose a
linearization on the branch-class that corresponds to requiring the
$r^{\Faber}_{g,\al} = r^g_{\al}-2g+1$ fixed branch points to go to
$0$.  Hence in any contributing torus-fixed locus, 
the amount of branching over $\infty$ is at most $2g-1$.

We now classify the fixed loci which can appear in the ``rational
tails'' case, where we must have a smooth irreducible component of
genus $g,$ and then consider the evaluation of their contributions
(\textbf{L1} to \textbf{L5} below), although some require further
elaboration.  Fixed loci correspond to
maps of the following sort (see Figure~\ref{twoexamples}), and in each
case we shall see that the same statements  hold for both the
left hand side (the simple case) and the right hand side (the
composite case) of Figure~\ref{twoexamples}, although the arguments
differ slightly.  The components mapping surjectively onto $\proj^1$ are
trivial covers (see \textbf{L1}).  Over $0$, there can be smooth
points (as in Figure~\ref{twoexamples}(a)) (see \textbf{L3}), nodes
(Figure~\ref{twoexamples}(b)), or contracted components
(Figure~\ref{twoexamples}(c)) (see \textbf{L2}).  Over $\infty$,
either there is no ``sprouting'', and the preimage of $\infty$
consists of smooth points (the ``\emph{simple}'' case, see the left
side of Figure~\ref{twoexamples}), or the target ``sprouts'', and we
obtain a relative stable map to an unrigidified target (that we denote
by $T_R$), with possibly-disconnected source (the
``\emph{composite}'' case, see the right side of
Figure~\ref{twoexamples}) (see \textbf{L4}).  In the composite case,
let $z$ be the node where the ``sprouted'' portion $T_R$ of
the target meets the ``unsprouted'' portion.  If $T_L$ denotes the
``unsprouted'' portion of $T,$ then $\infty_{T_L}$ and
$0_{T_R}$ are identified.  Note that in the composite case,
the sprouted portion of $T$ need not be a single $\proj^1;$ it could
be a chain of $\proj^1$'s.

Over $\infty$, we can have no genus $g$ components
by Remark~\ref{cr2} --- the only possible branching is over the two
endpoints ($z$ and $\infty_T$), together with at most $2g-1$ more.  Thus the
genus $g$ component must lie over $0$, and all irreducible curves over
$\infty$ must be genus $0$.  (In particular, in the second figure in
Figure~\ref{twoexamples}, there cannot be any contracted components mapping to $T_R$, so the picture is misleading.)

\noindent\underline{\emph{Geometry-combinatorics lexicon:}} 
With the following combinatorial elements, we associate the following  geometrical information. 
\newenvironment{Boxedminipage}%
{\begin{Sbox}\begin{minipage}}%
{\end{minipage}\end{Sbox}\fbox{\TheSbox}}
\begin{center}
\begin{Boxedminipage}{6.25in}
\noindent\emph{\underline{Vertices and edges of $\rt\in\cT_{g,m}$}}: \\
\noindent-- $0$-\emph{vertex}: $\leftrightarrow$ connected component of the preimage of $0$ 
              (contracted curve, node or  smooth point).  \\
\noindent-- $\infty$-\emph{vertex}:   $\leftrightarrow$ connected component of pre-images of $\infty_X.$  \\  
\noindent-- $t$-\emph{vertex}:    $\leftrightarrow$   preimage of $\infty_T.$    \\
\noindent--  $0\infty$-\emph{edge}:  $\leftrightarrow$ trivial cover of $\P^1;$ edge joins vertices
                corresponding to the  loci it meets.  \\
\noindent--  root $\bullet$: $\leftrightarrow$ genus $g$ (contracted) curve.

\noindent\emph{\underline{Edge weights}}: \\               
\noindent--  \emph{weight $\epsilon(e)$ on the  $0\infty$-edge $e$}:  $\leftrightarrow$ degree of 
               corresponding trivial cover. \\
\noindent--  \emph{weight on an $\infty t$-edge}: $\leftrightarrow$ contribution to ramification 
               over $\infty_X$ on component specified by the $\infty$-vertex.
\end{Boxedminipage}
\end{center}
\begin{center}
\begin{Boxedminipage}{6.25in}
\noindent\emph{\underline{Partitions}}:  \\
\noindent-- $\delta^v$ for a $0$-vertex  $v$: 
    \emph{Comb.}  formed by weights on $0\infty$-edges incident with  $v$.
    \emph{Geom.}  formed by the degrees of the trivial covers meeting that component. \\
\noindent-- $\beta^v$ for an $\infty$-vertex $v$:
     \emph{Comb.} formed by weights on $0\infty$-edges incident with $v$.
     \emph{Geom.}  formed by the degrees of the trivial covers meeting that component. \\
\noindent-- $\gamma^v$  for an $\infty$-vertex $v$: 
     \emph{Comb.} formed by weights on $\infty t$-edges incident with $v$.
     \emph{Geom.}  specifies ramification over $\infty_X$ on that component.      \\                            
\noindent-- $\alpha$: 
     \emph{Comb.}   see equation~(\ref{spal}).
     \emph{Geom.}  ramification over $\infty$.
\end{Boxedminipage}
\end{center}
\begin{center}
\begin{Boxedminipage}{6.25in}
\noindent\emph{\underline{Conditions}}: \\
\noindent-- $|\be^v|=|\ga^v|$:
       \emph{Comb.}  see equation~(\ref{balance}).
       \emph{Geom.}  the  degree of the map
from this component to $T_R$ is $| \be^v |$ (by examining the preimage
of $z$) and $| \ga^v |$ (from the preimage of $\infty_T$).
\end{Boxedminipage}
\end{center}
\begin{center}
\begin{Boxedminipage}{6.25in}
\noindent\emph{\underline{Constants}}: \\
\noindent-- $\eta_0(\rt)$:
       \emph{Comb.}  number of non-root $0$-vertices.
       \emph{Geom.} number of connected components of the preimage of $0$ excluding
                                    the component containing the contracted genus $g$ curve. \\
\noindent-- $n$: 
       \emph{Comb.}  $l(\al).$
       \emph{Geom.}  number of pre-images of $\infty_T.$ \\
\noindent-- $m$:  
       \emph{Comb.} degree of the root-vertex $\bullet$ (note $m \leq n$).
       \emph{Geom.}  number of trivial covers of $\proj^1$ meeting the 
genus $g$ contracted curve.
\\
\noindent-- $r_\infty$:
       \emph{Comb.}  see equation~(\ref{rrr}).
       \emph{Geom.}  the total branching over~$\infty_X.$ \\
\noindent-- $r_{\de^v}^0$ for a $0$-vertex  $v$: 
       \emph{Comb.}  see equation~(\ref{rgadef}).
       \emph{Geom.}  the total branching contributed by 
                                     the component corresponding to $v.$ \\
\noindent-- $r_{\ga^v,\be^v}^0$ for an $\infty$-vertex  $v$: 
       \emph{Comb.}  see equation~(\ref{rgabdef}).
       \emph{Geom.}  the total branching contributed by the component
corresponding to $v$  \\
\noindent-- $d$:   
       \emph{Comb.}  $d=|\al|.$
       \emph{Geom.}  degree of the relative stable map
\end{Boxedminipage}
\end{center}

\noindent\underline{\emph{Relative virtual localization:}}
Relative virtual localization (\cite[Thm.~3.6]{thmstar}, see \cite[Sec.~3.7]{thmstar}
for the special case of target $\proj^1$, and \cite[Sec.~4]{vl}
for the non-relative case) tells us that the contribution of a
fixed locus can be deduced by looking at the various parts of
Figure~\ref{twoexamples}.  In what follows, $t$ is the generator of
the equivariant cohomology (or Chow) ring of a point, $t \in
A^1_{\C^*}(pt) = \Q[t]$, although we will quickly forget the $t$.  The
relative virtual localization formula (abbreviated to \textbf{L} below) gives the
following contributions associated with the salient ``parts'' of the
graph.   Each of the factors in the summand of \eqref{RaviR} will
be readily derivable, with the exception of the sign and $D(\rt)$
which will require more attention.

Items \textbf{L1} to \textbf{L4} below come from the description of the 
fixed loci, and \textbf{L5} arises from 
the  cohomology class corresponding to fixing some branch points.
Moreover, \textbf{L1} to \textbf{L3} hold for both the simple and the composite case.

The results for the simple case and the composite case are the same, but
the arguments are slightly different.

\textbf{L1:}
For each trivial cover of the target $\proj^1$ of
degree $a$, we have  contribution  $a^a/ (a! t^a)$.
Hence we obtain the product  
$\boxed{ \prod \frac {q^q } {q!}}$ 
 in~\eqref{RaviR} by collecting those contributions associated
with the (genus $g$) root $0$-vertex.  

\textbf{L2:}
For each contracted curve above $0$ of genus $h$ (Figure~\ref{twoexamples}(c)),
meeting trivial covers (components mapping surjectively onto $\proj^1$) 
of degree $\al_1$, \dots,
$\al_m$ respectively,
we have a contribution\lremind{genush}
\begin{equation}\label{genush}
 t^{-1}(t^h - \la_1 t^{h-1} + \cdots + (-1)^h \la_h)  \prod_{i=1}^m 
\frac t {t / \al_i- \psi_i}.
\end{equation}
This contribution is on the factor $\cmbar_{h,m}$ corresponding
to the contracted curve.    We will have $h=0$ or $h=g$, as all components
have one of these two genera.

\textbf{L3:} 
For each node above $0$ (Figure~\ref{twoexamples}(b)) joining
  trivial covers of degrees $\al_1$ and $\al_2$, we get a contribution
  of $t^{-1} \al_1 \al_2 / (\al_1+\al_2)$.  For each smooth point
  above $0$ (Figure~\ref{twoexamples}(a)), on a trivial cover of
  degree $\al_1$, we get $t^{-1} / \al_1$.  Using the contribution
 $\prod (\al_i^{\al_i} / \al_i!)$
  from L1, we obtain
  the factor $\boxed{ H^0_{\de^j} / r^0_{\de^j}!}$ in \eqref{RaviR} for each (non-root)
  $0$-vertex of degree $1$ or $2$, using the formula \eqref{singHur}
  for genus $0$ single Hurwitz numbers.  We also obtain
the product of the $\boxed{\ep(e)}$  in \eqref{RaviR} corresponding to $0 \infty$-edges $e$
meeting a degree $1$ or $2$ (genus $0,$ non-root) $0$-vertex.

\textbf{L4:} In the composite case (Figure~\ref{twoexamples}(d)), we
have a contribution of $1 / (-t-\psi_z)$, where $\psi_z$ is the first
Chern class of the line bundle corresponding to the cotangent space of
$T_R$ at $z$.

\textbf{L5:} 
{}From the pullback of the linear space by the branch morphism
(informally, requiring $r^{\Faber}_{g,\al}$ branch points to map to $0$),
we have a contribution of 
$\left(   (2g-r_\infty) t \right)
\left(   (2g+1-r_\infty) t \right)
\cdots 
\left(   (r^g_{\al} - r_\infty) t \right)$
(there are $r^{\Faber}_{g,\al}$ factors), from which we obtain 
$\boxed{ 
r^{\Faber}_{g,\al}!
\binom {  r^g_{\al} - r_{\infty}}   {
r^{\Faber}_{g,\al}}
}$ in \eqref{RaviR}.

So {\bf L1--L4} arise from the fixed loci, and {\bf L5} arises from
the cohomology class corresponding to fixing some branch points.  We
take the product of these contributions, and read off the constant
($t^0$) term to obtain the contribution of this fixed locus.

The result is a $(2g-1)$-dimensional class on $\cm_{g,n}^{rt}$.  One
of the ingredients (from \textbf{L3}) is a (tautological) class on
$\cmbar_{g,m}$ corresponding to the contracted genus $g$ curve mapping
to $0$.  Now all tautological classes of dimension {\em less} than
$2g-1$ vanish on $\cmbar_{g,m}^{rt}$ by Remark~\ref{summary}(i).  Thus a
non-zero contribution is possible only by taking the contribution of a
class of dimension {\em precisely} $2g-1$ on $\cm_{g,m}^{rt}$, and
thus the contributions from every other ingredient must have dimension
$0$.  In light of this observation, we list the contributions from
each of the parts of Figure~\ref{twoexamples}, ignoring the
equivariant parameter $t$.
Also, from \textbf{L3}, the contribution by the contracted genus $g$
curve is $\boxed{ \P^g_m(\de^0)}$ so, with the contribution from \textbf{L1},
we obtain the term $\boxed{A(\rt)},$ a term in~\eqref{RaviR}. In addition, 
we obtain the product of the $\boxed{\ep(e)}$ over those edges $e$ meeting the (genus $g$) root  $0$-vertex.

{}From \eqref{genush} (using $h=0$),
any contracted genus  $0$ component
over $0 \in \proj^1$ meeting trivial covers of degrees
$\de^j_1$, \dots, $\de^j_p$ ($p \geq 3$ by the stability condition)
gives the dimension $0$ 
contribution
\begin{equation*}\label{g0}
\left( \prod_{i=1}^p \de^j_i \right) \int_{\cmbar_{0,p}} \frac 1 { (1 - \de^j_1 \psi_1) \cdots (1 - \de^j_p \psi_p)}.
\end{equation*}
Combining this with $\prod_{i=1}^p \frac { {\de^j_i}^{\de^j_i}} {\de^j_i!}$
from \textbf{L1}, using the ELSV formula \cite{elsv1, elsv2, gvelsv}
we obtain $(\prod \de^j_i) H^0_{\de^j} / r^0_{\de^j}!$.
(We could have bypassed the ELSV formula, using instead the
formula \eqref{singHur} for genus $0$ Hurwitz numbers and the string
equation, given in Proposition~\ref{sd}.)
Thus 
we obtain the factor $\boxed{H^0_{\de^j} / r^0_{\de^j}!}$ for each non-root $0$-vertex of degree
at least three. 
 We also obtain the product of the $\boxed{\ep(e)}$
corresponding to $0 \infty$-edges $e$ meeting all (genus $0$, non-root) $0$-vertices
of degree at least $3$.  \textbf{L3} gives the same values for degree $1$ and $2$ so,
combining the three sources for the $\epsilon$'s, we now obtain the entire product 
$\boxed{B(\rt)},$
and combining the two sources for the  $H^0_{\de^j} / r^0_{\de^j}!$, we obtain the entire product  
$\boxed{C(\rt)}.$

The contributions from \textbf{L1--L3}, and \textbf{L5} are now exhausted. 

\newcommand{\rmb}{\scriptstyle{\scriptscriptstyle{\cmbar}}}  

It remains to obtain the sign $(-1)^{r_\infty}$, the term $D(\rt)$, and
the division by $\eta_0(\rt)!$.
To do so, we now appeal to \textbf{L4}. In the composite case (if the target ``sprouts'', 
Figure~\ref{twoexamples}(d)), suppose $\cmbar := \cmbar_{\be, \al}$ is the moduli space 
of relative stable maps to the unrigidified $T_R$,  where $\be$ is the kissing partition 
above $z$.  It is the moduli space of  relative stable maps, where the source is a disjoint 
union of genus $0$ curves.   As always for maps to unrigidified targets, the virtual
dimension of this space is one less than the number of ``moving branch points'' $r_{\cmbar}$. 
Then the contribution is the dimension $0$ portion of
$\frac 1 {-1-\psi_z} [ \cmbar ]^{\virt}$
which is\lremind{overinfinity}
\begin{equation}\label{overinfinity}
(-1) ^{r_{\cmbar}} \psi_z^{r_{\cmbar}-1}  
[ \cmbar ]^{\virt}.
\end{equation}

The sign gives us the factor $\boxed{  \prod  (-1)^{r^0_{\ga^i,\be^i}}}$ which, from~\eqref{rrr}, is equal
to  $(-1)^{r_\infty},$ the sign in \eqref{RaviR}.

The proof of \cite[Lem.~4.8]{thmstar} (see also
\cite[Fig.~2]{thmstar}) shows that $\psi^a$ applied to $\cmbar$ can be
interpreted as requiring that the target break into $a+1$ components.
More precisely, $\psi^a$ applied to $[\cmbar]^{\virt}$ is the same as
gluing virtual fundamental classes of relative stable maps to the
$a+1$ components of $T_R$ (in the same sense as the
degeneration formula, with kissing multiplicities arising for each
node of the target $T_R$), divided by $(a+1)!$.
This latter term gives a factor of $\boxed{1/r_{\infty}!}$.

In particular, $\psi^{r_{\cmbar}-1}$ corresponds to the target
breaking into $r_{\cmbar}$ components.  Because the resulting map must
be stable, there must be some branching on each of these components
(away from the nodes of $T_R$, and $z$ and
$\infty_T$).  Thus as the total amount of branching is $r_{\cmbar}$
away from the nodes of $T_R$, and there is precisely
this number of components, we must have branching number~$1$ on each
irreducible component of $T_R$.  By 
Remark~\ref{cr1}(b), above each component of the
$T_R$, we must have precisely one almost-trivial
cover, along with some trivial covers. 

Thus the contribution of \eqref{overinfinity} is the size of a
discrete set (counted modulo automorphisms).  This set counts the
number of branched covers of $T_R$ (a chain of
$r_{\cmbar}$ $\proj^1$'s), with one simple branching on each
component, and given branching $\ep$ over $z$ (a point at the end of
the chain) and $\al$ over $\infty_T$ (a point at the other end of the
chain), satisfying the kissing condition over each node of
$T_R$.  By the gluing formula (or indeed, the much
older technique of just studying the degeneration), this is the number
of branched covers of $\proj^1$ by a union of genus $0$ curves with
branching given by $\ep$ and~$\al$ over two points $z$ and $\infty_T$,
and simple branching over $r_{\cmbar}$ other given fixed points.

We shall now see that this is (up to a combinatorial factor) a product
of genus $0$ double Hurwitz numbers.  Recall that a genus $0$ double
Hurwitz number 
\begin{equation}\label{eHalbe}
H^0_{\al,\vep}
\end{equation}
 (where $\al$ and $\vep$ are partitions
of some number $e$) counts the number of degree $e$ covers of
$\proj^1$ by $\proj^1$, with branching $\al$ at one fixed point,~$\vep$
at another, and simple branching at $r^0_{\al,\vep}$ other
fixed points.  Suppose we are considering covers by $N$
$\proj^1$'s, where component $i$ corresponds to the subpartition
$\be^i$ of $\vep$ (over $z$) and the subpartition $\ga^i$ of $\al$
(over $\infty_T$).  (Thus $| \be^i | = | \ga^i |$ is the degree of
that subcover, $\vep = \coprod_i \be^i$, and $\al = \coprod_i \ga^i$.)
Then component $i$ has simple branching over $r^i_{\infty} :=
r^0_{\be^i,\ga^i}$ of the fixed simple branch
points.  There are $\binom{ \sum_i r^i_{\infty} } { r^1_{\infty},
  \dots, r^N_{\infty}}$ ways of partitioning the branch points into
these $N$ sets.  Once this partition is chosen, there are $\prod_{i=1}^N
H^0_{\be^i,\ga^i}$ such branched covers.  (One caution: we have
cavalierly described $H^0_{\be^i,\ga^i}$ as enumerating a set.  In
reality, each cover is counted with multiplicity equal to the inverse
of the size of its automorphism group, so $H^0_{\be^i,\ga^i}$ need not
be integral, and in fact is not precisely for $\be^i=\ga^i = (d_i)$;
trivial covers of degree $d_i$ ``count for'' $1/d_i$.)  We have
obtained the factors $\boxed{ H^0_{\be^k, \ga^k}}$ and
$\boxed{\binom{ r_{\infty} } { r^1_{\infty}, \dots, r^N_{\infty}}}$ in
\eqref{RaviR}. (The numerator $r_{\infty}!$ in the multinomial
coefficient cancels the $1/r_{\infty}!$ from earlier.) This is the term $D(\rt).$

Finally, the $\boxed{ 1/\eta_0(\rt)!}$ in \eqref{RaviR} is present because
the trees have labelled non-root $0$-vertices (Def.~\ref{Dloctr}).
\epf


\part{ALGEBRAIC COMBINATORICS}\label{partC}
At this point, we have defined the Faber-Hurwitz classes $\F^{g,\al},$
which ``virtually'' correspond to ``rational tail'' curves admitting
a branched cover of $\proj^1$ with branching at $\infty$ corresponding to $\al,$
and ``all but $2g-1$ branching fixed''. Such classes are a multiple of 
a basic class $\G_{g,1};$ this multiple is the Faber-Hurwitz \emph{number} 
$F^{g,\al}$. By degeneration, we have obtained
Corollary~\ref{JCEforFHS}, the Degeneration Theorem for
the generating series $F^g$ for these numbers,
involving the genus $0$ Hurwitz series $\widH^0$.
By localization, we have also obtained Theorem~\ref{bigyuck}, the
Localization Tree Theorem, which describes these classes
(or numbers) as a sum over certain rooted, labelled trees, 
involving genus $0$ Hurwitz, double Hurwitz numbers and the desired
intersection numbers (of $\psi$-classes).
Theorem~\ref{onedim} shows us that we can ``invert'' this expression, to determine
intersection numbers in terms of Hurwitz numbers and double Hurwitz numbers. Our goal is to
formalize this.
The strategy is to show that Localization Tree Theorem and the Degeneration Theorem,
taken together, give a non-singular system of linear equations for 
the top Faber intersection numbers, so it has a unique solution, and that the conjectural values satisfy it.

We accomplish this by a sequence of transformations, which yield a number of refined
versions of the Localization Tree Theorem and the Degeneration Theorem.
These versions of the Localization Tree Theorem (we say that these are results for
the ``\textsf{localization side}'') are given by the sequence
\begin{equation*}
\boxed{\mbox{Thm.~\ref{bigyuck} $\leadsto$ Thm.~\ref{branchfg}
$\leadsto$ Cor.~\ref{tSFHT4}  $\leadsto$ Cor.~\ref{PsiPhieq} } }
\end{equation*}
These versions of the Degeneration Theorem (we say that these are results for
the ``\textsf{degeneration side}'') are given by the sequence
\begin{equation*}
\boxed{\mbox{Thm.~\ref{Fcutjoin} $\leadsto$ Cor.~\ref{JCEforFHS} $\leadsto$ Cor.~\ref{tSJCE}
$\leadsto$ Lem.~\ref{symFjcutm} }}
\end{equation*}

\section{Exponential generating series for localization trees}\label{s:genseries}
The purpose of this section is to ``evaluate'' the sum over localization trees that arises from the 
localization arguments in Theorem~\ref{bigyuck}. Localization trees are a class of rooted, labelled
trees, and we use the standard multivariate exponential generating series for combinatorial
structures with many sets of labels, as well as variants of the standard branch decomposition
for rooted trees.

\subsection{Exponential generating series and the $\star$-product}\label{preliminaries}

For a localization tree $\rt$, let $\eta_k(\rt)$ denote the number of $t$-vertices
in $\rt$ that are incident with an edge of weight $k$, $k\ge 1$.

\begin{df}\label{expgf}
Let $\sA$ be a set of localization trees, with \emph{weight function}  $\wt$.
Then the \emph{exponential generating series for $\sA$ with
respect to $\wt$} is
\begin{equation*}
\left[\sA,\wt\right]_{\bf\eta}:=\sum_{\rt\in\sA}
\f{p_1^{\eta_1(\rt)}p_2^{\eta_2(\rt)}\cd}{\eta_0(\rt)!\eta_1(\rt)!\eta_2(\rt)!\cd}
\wt(\rt).
\end{equation*}
\end{df}

For now, we shall allow the range of the weight function to be any ring, or even
a vector space, and in particular we allow geometric classes.
Note that, as a formal power series in $p_1,p_2,\ld$, $\left[\sA,\wt\right]_{\bf\eta}$ is
always well-formed 
because of the balance condition~(\ref{balance}),
which ensures that there is only a finite number of localization trees $\rt$ with
 $\eta_k(\rt)=i_k$, $k\ge 1$, for each $i_1,i_2,\ld$ (so the coefficients are finite
sums of weight function values).

Since localization trees are labelled objects, we consider a particular version
of the standard $\star$-product for them, which is the Cartesian product together
with a ``label-distribution'' operation.
We define this $\star$-product as
follows. Consider two localization trees $\rt_1$ and
 $\rt_2$.
Suppose that $\eta_k(\rt_1)=i_k$, and $\eta_k(\rt_2)=j_k$, for $k\ge 0$.
Now choose subsets $\al_k\subseteq [i_k+j_k]$, with $|\al_k|=i_k$, for $k\ge 0$,
and let $\be_k=[i_k+j_k]\setminus\al_k$ (so $|\be_k|=j+k$)
(we use the notation $[n]=\{ 1,\ld,n\}$).
Let $\al(\rt_1)$ be the tree obtained from $\rt_1$ by relabelling the
labelled vertices as follows: replace the
label $m$ on a non-root $0$-vertex by the $m$th smallest element of $\al_0$,
for $m=1,\ld, i_0$; replace the label $m$ on a $t$-vertex incident with an edge
of weight $k$ by the $m$th smallest element of $\al_k$, for $m=1,\ld ,i_k$, $k\ge 1$.
Let $\be(\rt_2)$ be the tree obtained from $\rt_2$ by relabelling the
labelled vertices with the elements of $\be_0,\be_1,\ld$, in the analogous manner.
We call $(\al,\be)$ a \emph{compatible relabelling} (of $(\rt_1,\rt_2)$).
Where convenient, we also refer to $\rt_1$ and $\rt_2$ as \emph{canonically
labelled}.

\begin{df}\label{starprod}
Let $\sA,\sB$ be sets of localization trees. Then
\begin{equation*}
\sA\star\sB:=
\{(\al(\rt_1),\be(\rt_2)):(\rt_1,\rt_2)\in\sA\times\sB,(\al,\be)\;\rm{a}\;
\rm{compatible}\;\rm{relabelling}\;\rm{of}
\; (\rt_1,\rt_2)\}.
\end{equation*}
\end{df}

The reason for using exponential generating series for labelled combinatorial
objects is the Product Lemma, given in the following result. 
For a proof, see, \textit{e.g.}, Goulden and Jackson~\cite[Lem.~3.2.11]{gjce}.

\begin{lm}[Product Lemma (for localization tree generating series)]\label{prodlem}
Let $\cS,\sA,\sB$ be sets of localization trees, with weight functions $\wt,\wt_1,\wt_2$,
respectively. Suppose there is a bijection
\begin{equation*}
\cS\bij\sA\star\sB:\rt\mapsto (\al(\rt_1),\be(\rt_2)),
\end{equation*}
subject to $\wt(\rt)=\wt_1(\rt_1)\wt_2(\rt_2)$ and $\eta_k(\rt)=\eta_k(\rt_1)+\eta_k(\rt_2)$,
$k\ge 0$. Then
\begin{equation*}
[\cS,\wt]_{\eta}=[\sA,\wt_1]_{\eta}[\sB,\wt_2]_{\eta}.
\end{equation*}
\end{lm}

For any set $\cS$ of localization trees $\cS$,
let $\cS^{\star k}:=\cS\star\cd \star\cS$, where there
are $k$ $\cS$'s in this $k$-fold $\star$-product. We
define $\cU_k(\cS):=\cS^{\star k}/\mathfrak{S}_k$,
under the natural action of $\mathfrak{S}_k$, for $k\ge 0$.

\subsection{Branch decompositions for localization trees}\label{branchdec}

In order to decompose localization trees, we require three variants of the
standard branch decomposition, described below. These are formalized
as combinatorial mappings, called $\ombr$, $\ombr_0$ and $\ombr_{\infty}$.

First, we give $\ombr$:
If the root-vertex $\bullet$ and incident edges are deleted from a localization tree $\rt$, then
we obtain a list $(\rt_1,\ld ,\rt_m)$ of rooted trees on mutually
distinct sets of vertices, each inheriting as root vertex
its unique vertex that was adjacent to $\bullet$ in $\rt$, where $m$ is the degree
of the root-vertex of $\rt$. Now
let $\rt_i^{\prime}$ be obtained from $\rt_i$ by
joining the root vertex of $\rt_i$ to a new copy of $\bullet$ (which
becomes the new root-vertex of $\rt_i^{\prime}$), joined by
an edge whose weight is equal to the weight of the
edge joining $\bullet$ to the root-vertex of $\rt_i$ in $\rt$.
Since $\bullet$ is unlabelled in $\rt$, then $(\rt_1',\ld ,\rt_m')$ is
equal to $(\al_1(\rt_1''),\ld ,\al_m(\rt_m''))$, a compatible
relabelling of canonical localization trees $(\rt_1'',\ld ,\rt_m'')$.
Then we define $\ombr(\rt)=(\al_1(\rt_1''),\ld ,\al_m(\rt_m''))$.
This corresponds to removing the root from the tree, and describing
the remainder as a list of trees.

Second, we give $\ombr_0$:
Suppose that $\rt$ is a localization tree whose root-vertex $\bullet$ is
a monovalent $0$-vertex.  Let the $\infty$-vertex adjacent
to $\bullet$ be $u$, and
let $\rt^{\prime}$ be the tree obtained  by deleting $\bullet$ and incident edge
from $\rt$, and deleting all $t$-vertices adjacent to $u$, together with their incident edges,
and rooting the resulting tree
at $u$. Now form a new graph $G$, containing these deleted $t$-vertices (labelled as in $\rt$), joined to
a new $\infty$-vertex $w$ by an edge of the same weight as the deleted incident edge in $\rt$.
Then we define $\ombr_0(\rt)=(G,\ombr(\rt^{\prime}))$.
Now consider partition $\ga =(1^{a_1}2^{a_2}\ld)$, and let $\rw_{\ga}$ be the 
graph consisting of a single $\infty$-vertex, joined
by an edge of weight $k$ to $a_k$ canonically labelled monovalent $t$-vertices,
 $k\ge 1$.
Then note
that $\ombr_0(\rt)=(\al_0(\rw_{\ga}),\al_1(\rt_1''),\ld ,\al_m(\rt_m''))$, a
compatible relabelling of the canonical $(\rw_{\ga},\rt_1'',\ld ,\rt_m'')$,
where $\ga=\ga^u(\rt)$, $|\ga|=|\be^u(\rt)|$ (because of the balance
condition~(\ref{balance}) at $u$), and $\be^u(\rt)$ has $m+1$ parts, $m\ge 0$.

Third, we give $\ombr_{\infty}$:
Suppose that $\rt$ is a localization tree whose root-vertex $\bullet$ is
a monovalent $\infty$-vertex. Let the $0$-vertex adjacent
to $\bullet$ be $v$, with label $i$,
then let $\rt^{\prime}$ be the tree obtained  by removing $\bullet$
 from $\rt$, and removing the label from the vertex $v$, and
rooting the resulting tree
at $v$. Let $\rv_i$ be the graph consisting of a single, $0$-vertex, labelled $i$.
Then we define $\ombr_{\infty}(\rt)=(\rv_i,\ombr(\rt^{\prime}))$,
and note that $\ombr_{\infty}(\rt)=(\al_0(\rv_1),\al_1(\rt_1''),\ld ,\al_m(\rt_m''))$,
a compatible relabelling of the canonical $(\rv_1,\rt_1'',\ld ,\rt_m'')$,
where $\de^v(\rt)$ has $m+1$ parts, $m\ge 0$.

\subsection{Generating series form for the Localization Tree Theorem}\label{genserFHsT}

In the next result, Theorem~\ref{branchfg}, which is the second form of the
Localization Tree Theorem, we
introduce a generating series, $\zeta^g$, for the tree summation that
arose in Theorem~\ref{bigyuck}, which was the first form of the
Localization Tree Theorem. We also introduce two ancillary classes
of generating series $f_j,g_j,j\ge 1$. These are all
exponential generating series with respect to particular weight functions,
and the proofs of the equations that relate them
are combinatorial, applying the variants of the branch decomposition,
together with the $\star$-product and the Product Lemma.
However, the statement of Theorem~\ref{branchfg} is purely
algebraic (without reference to trees), and the series $\zeta^g$, $f_j$, $g_j$ are
all uniquely defined by the equations in the statement.
Note that the coefficients in $\zeta^g$ are geometric \emph{classes}, while the
coefficients in $f_j$, $g_j$, $F^g$ are \emph{rationals}.

The statement of the theorem involves the generating series $\widH^0$ for
genus $0$ Hurwitz numbers, defined in~(\ref{sHser}), and 
the generating series for genus $0$ double Hurwitz numbers,
given by
\begin{equation}\label{dHser}
H^0(z,u;\bfp;\bfq):=\sum_{\substack{\al ,\be\in\sP, \\ |\al |=|\be |}}
z^{|\be |} p_{\al}q_{\be}u^{l(\be)}
\frac{H^0_{\al ,\be}}{r^0_{\al ,\be}!|\ts{Aut}\al |\, |\ts{Aut}\be |},
\end{equation}
where $\bfp =(p_1,p_2,\ld )$, $\bfq =(q_1,q_2,\ld )$, and $r^0_{\al ,\be}$ is
defined in~(\ref{rgabdef}).

\begin{tm}[Localization Tree Theorem --- functional equations]\label{branchfg}
For $g\ge 1$, we have
\begin{equation*}
a)\quad F^g(z;\bfp )\; \G_{g,1}=[u^{2g-1}]\;\zeta^g\left(\f{z}{1-u},-u;\f{-u}{1-u}\bfp\right),
\end{equation*}
\newline
where 
$$b)\quad \zeta^g(z,u;\bfp )=\sum_{n\ge 1}\frac{1}{n!}
\!\!\!\!
\sum_{\substack{a_1,\ldots,a_n, k\ge0, \\ a_1+\cdots+a_n+k=g-2+n}}
\!\!\!\!\!\!\!\!
(-1)^k \laFab\tau_{a_1}\cd\tau_{a_n} \la_k \raFab
\prod_{j=1}^n \xi^{(a_j)}(z,u;\bfp),$$
\newline
with
\begin{equation}\label{Fseries}
\xi^{(i)}(z,u;\bfp ):=\sum_{j\geq 1}\frac{j^{j+i}}{j!}f_j(z,u;\bfp ),\qquad i\ge 0,
\end{equation}
\newline
and, for $j\geq 1,$
$\mbox{c)}\quad 
\quad f_j=u^{-2}\left.\left( j\frac{\pa}{\pa q_j}H^0(z,u;\bfp ;\bfq)\right)\right|_{q_i=g_i,i\geq 1},\qquad
\mbox{d)}\quad g_j=\left.\left( j\frac{\pa}{\pa q_j}\widH^0(1;\bfq)\right)\right|_{q_i=f_i,i\geq 1},$
with $f_j=f_j(z,u;\bfp )$, $g_j=g_j(z,u;\bfp ).$
\end{tm}

\bpf
As special cases of weight functions for localization trees $\rt$,
consider $\wt_{g}(\rt):=A(\rt)\wt_{0}(\rt)$,
 $\wt_{0}(\rt):=B(\rt)C^{\ddagger}(\rt)D(\rt)$
$\prod_{v\in\sV_\infty}u^{l(\be^v)-2}z^{|\be^v|}$,
and $\wt_{\infty}(\rt):=$
$B(\rt)C(\rt)D^{\ddagger}(\rt)$ $\prod^{\ddagger}_{v\in\sV_\infty}u^{l(\be^v)-2}z^{|\be^v|}$,
where $A,B,C,D$ were defined in~(\ref{Aonly}),~(\ref{ABCD}).
Now define generating series
\begin{equation*}\label{combgsTg}
\zeta^g(z,u;\bfp ):=[\cup_{m\ge 1}\sT_{g,m},\wt_g]_{\bf \eta},\qquad g\ge 1,
\qquad f_j:=[\sT_{0,j},\wt_0]_{\eta},
\quad g_j:=[\sT_{\infty,j},\wt_{\infty}]_{\eta},\quad j\ge 1.
\end{equation*}

For part (a), in $\zeta^g$ we
have $p_1^{\eta_1(\rt)}p_2^{\eta_2(\rt)}\cd/\eta_1(\rt)!\eta_2(\rt)!\cd =p_{\al}/|\ts{Aut}\al |$,
from~(\ref{spal}). Thus, for $g\geq 1$ and $|\al |=d\geq 1$, we
have $\F^{g,\al }=r^{\Faber}_{g,\al}!|\ts{Aut}\al |\sum_{r_{\infty}}
\binom{r^g_{\al}-r_{\infty}}{r^{\Faber}_{g,\al}}(-1)^{r_{\infty}}
\left[p_{\al}u^{r_{\infty} - l(\al)}z^{d}\right]$
$\zeta^g(z,u;\bfp) $, from Theorem~\ref{bigyuck}.
But
$\binom{r^g_{\al}-r_{\infty}}{r^{\Faber}_{g,\al}}=
[u^{2g-1-r_{\infty}}]{(1-u)^{-(d+l(\al ))}},$
so
$\F^{g,\al }=r^{\Faber}_{g,\al}!|\ts{Aut}\al |$ $[p_{\al}u^{2g-1}z^d]$
$\zeta^g\left(\frac{z}{1-u},-u;\frac{-u}{1-u}\bfp\right) $.
Then part (a) follows from~(\ref{gsFnum}) and~(\ref{ianeq}).

For part (b), define $\sT_{g,\de}$ to be the set of localization trees $\rt$
 in $\sT_{g,m}$ with $\de^{\bullet}(\rt)=\de$, for any partition $\de$ with $m$ parts,
 $m\ge 1$.
Then variant $\ombr$ of the branch decomposition gives a bijection
\begin{equation*}
\ombr:\sT_{g,\de}\bij \cU_{a_1}(\sT_{0,1})\star\cU_{a_2}(\sT_{0,2})\star\cd
:\rt\mapsto (\al_1(\rt_1''),\ld ,\al_m(\rt_m'')),
\end{equation*}
where $\de=(1^{a_1}2^{a_2}\cd)$ with $a_1+a_2+\cd =m$. It is straightforward
to check in this bijection that $\eta_k(\rt)=\eta_k(\rt_1'')+\cd +\eta_k(\rt_m'')$, for $k\ge 0$.
Moreover, $\wt_g(\rt)={\P}_m^g(1^{a_1}2^{a_2}\cd)$ $\left(\prod_{j\ge 1}\left(j^j/j!\right)^{a_j}\right)$
$\wt_0(\rt_1'')\cd\wt_0(\rt_m'')$, so
from the Product Lemma we obtain
\begin{equation*}
[\sT_{g,\de},\wt_g]_{\eta}={\P}_m^g(1^{a_1}2^{a_2}\cd)\prod_{j\ge 1}\f{\left(j^jf_j/j!\right)^{a_j}}{a_j!},
\end{equation*}
and part (b) follows by summing this over all $a_1,a_2,\ld \ge 0$ with $a_1+a_2+\cd =m$, and $m\ge 1$,
and applying~(\ref{Fabpoly2}).

For part (c), for $j\ge 1$, variant $\ombr_0$ of the branch decomposition gives a bijection
\begin{equation*}
\ombr_{0}:\sT_{0,j}\bij\bigcup\{\rw_{\ga}\}\star\cU_{b_1}(\sT_{\infty,1})
\star\cU_{b_2}(\sT_{\infty,2})\star\cd:
\rt\mapsto (\al_0(\rw_{\ga}),\al_1(\rt_1''),\ld ,\al_m(\rt_m'')),
\end{equation*}
where the union is over all $\ga=(1^{a_1}2^{a_2}\cd)$,
 $\be=(1^{b_1}\cd (j-1)^{b_{j-1}}j^{b_j+1}(j+1)^{b_{j+1}}\cd)$
with $|\ga|=|\be|\ge j$, and $l(\be)=m+1$, $m\ge 0$ (here, $\ga=\ga^u(\rt)$ and
 $\be=\be^u(\rt)$, where $u$ is the $\infty$-vertex adjacent to $\bullet$ in $\rt$,
as given in the description of $\ombr_0$ above).
It is straightforward to check in this bijection
that $\eta_k(\rt)=\eta_k(\rw_{\ga})+\eta_k(\rt_1'')+\cd +\eta_k(\rt_m'')$, for $k\ge 0$.
Moreover, $\wt_{0}(\rt)=j
\left(H^0_{\ga,\be}u^{l(\be)-2}z^{|\be|}/r^0_{\ga,\be}!\right)$
 $\wt_{\infty}(\rt_1'')\cd\wt_{\infty}(\rt_m'')$, so
from the Product Lemma we obtain
\begin{equation*}
f_j=u^{-2}j\sum\f{H^0_{\ga,\be}u^{l(\be)}z^{|\be|}}{r^0_{\ga,\be}!}
\prod_{i\ge 1}\f{p_i^{a_i}}{a_i!}
\f{g_i^{b_i}}{b_i!},
\end{equation*}
since $[\rw_{\ga},1]_{\eta}=\prod_{i\ge 1}\f{p_i^{a_i}}{a_i!}$,
and part (c) follows immediately from~(\ref{dHser}).

For part (d), variant $\ombr_{\infty}$ of the branch decomposition gives a bijection
\begin{equation*}
\ombr_{\infty}:\sT_{\infty,j}\bij\bigcup\{\rv_1\}\star\cU_{a_1}(\sT_{0,1})
\star\cU_{a_2}(\sT_{0,2})\star\cd:\rt\mapsto (\al_0(\rv_1),\al_1(\rt_1''),\ld ,\al_m(\rt_m'')),
\end{equation*}
where the union is over all $\de=(1^{a_1}\cd (j-1)^{a_{j-1}}j^{a_j+1}(j+1)^{a_{j+1}}\cd)$ with
 $|\de|\ge j$, and $l(\de)=m+1$, $m\ge 0$. (Here, $\de=\de^v(\rt)$,
where $v$ is the $0$-vertex adjacent to $\bullet$ in $\rt$,
as given in the description of $\ombr_{\infty}$ above.) 
 It is straightforward to check in this bijection
that $\eta_k(\rt)=\eta_k(\rv_1)+\eta_k(\rt_1'')+\cd +\eta_k(\rt_m'')$, for $k\ge 0$,
and that $\wt_{\infty}(\rt)=j\left(H^0_{\de}/r^0_{\de}!\right)\wt_0(\rt_1'')\cd\wt_0(\rt_m'')$. 
Then from the Product Lemma we obtain
\begin{equation*}
g_j=j\f{1}{1!}\sum\f{H^0_{\de}}{r^0_{\de}!}\prod_{i\ge 1}
\f{f_i^{a_i}}{a_i!},
\end{equation*}
since $[\{\rv_1\},1]_{\eta}=1/1!$,
and part (d) follows immediately from~(\ref{sHser}).
\epf


\dmjdel{
\begin{minipage}{2.5in}
\newcommand{\e}{\varepsilon}
\newcommand{\bu}{$\bullet$}

\textbf{Edge weights}: integers~$>0$ \\
\quad--~$\e_i:$~on $i$-th {$0\infty$-edges} \\
\quad\bu~$c_i^j:$~on  {$\infty t$-edges} \\

\medskip

\textbf{Partitions:} indicated by a dark sector at a vertex. \\
\quad\bu~${\beta^1=(\e_1,\e_2,\e_4,\e_5)}$  etc. \\
\quad\bu~${\gamma_i=(c_i^1, \dots, c_i^{k_i})}$  \\
\quad\bu~${\delta^0=(\e_1,\e_6)}$ etc. \\
\quad\bu~Balance: ${\left |\beta^i\right |} ={ \left |\gamma^i\right |}$

\medskip
\textbf{Vertex weights: } \\
\quad\bu~{$H^0_{{\beta^v},{\gamma^v}}$}:  
{Hurwitz Number};
ramification {$\beta^v$} at $0$ \& {$\gamma^v$} at $\infty$
for $\infty$-vertex $v.$ \\
\quad\bu~${H^0_{{\delta^v}}}$: {Hurwitz Number};
ramification ${\delta^v}$ at $\infty$ for $0$-vertex $v.$ \\
\quad\bu~{$P_n^g({\delta^0})$}: {Top potential} 
at \bu-vertex \\[-4pt]

\end{minipage}
} 


\section{Symmetrization and a polynomial transformation}\label{sssSOD}

At this stage, we have established the two generating series results
that we need to prove Faber's Conjecture. The first of these results,
on the \textsf{degeneration side}, is Corollary~\ref{JCEforFHS}, the second
form of the Degeneration Theorem, which gives a
linear partial differential equation for the Faber-Hurwitz series $F^g$, in
terms of the genus $0$ Hurwitz series $\widH^0$. The
second of these results, on the \textsf{localization side}, is Theorem~\ref{branchfg},
the second form of the Localization Tree Theorem,
which expresses Faber-Hurwitz classes as a linear combination
of the intersections $ \laFab\tau_{a_1}\cd\tau_{a_n} \la_k \raFab$,
in terms of $\widH^0$ and the genus $0$ double Hurwitz series $H^0$.

In this section, in order to
prove Faber's Intersection Number Conjecture, we introduce three operators,
giving us a three-step transformation that will enable us to apply
Corollary~\ref{JCEforFHS} and Theorem~\ref{branchfg} conveniently,
so that we can extract the intersection numbers
of top degree, that are the subject of Faber's Conjecture.
The first step is a symmetrization operator $\Xi_m$,
and the second step is a change of variables $\sC$. The composition $\sC\Xi_m$ yields
polynomials in the new variables, and the third step is the
operator $\Ups'$, that restricts to terms of maximum total degree.
We shall refer to $\Ups'\sC\Xi_m=\Ups\Xi_m$ as the \emph{fundamental transformation}.

Note that our fundamental transformation is only a slight modification of
the three-step transformation that was used in~\cite{gjvlambdag} to
give another proof of Getzler and Pandharipande's  $\lambda_g$-Conjecture. There the
transformation was applied to the Hurwitz number generating series in arbitrary genus,
and the
first two steps, of symmetrizing and changing variables, were identical. In that
Hurwitz case, polynomiality also held, but the third step was to
restrict to (full) terms of \emph{minimum} total degree, so
that the intersection numbers of bottom degree, that are
the subject of the $\lambda_g$-Conjecture, could be extracted.
Changes of variables for similar purposes arise in work of
Kazarian and Lando~\cite{kl} and Shadrin and Zvonkine~\cite{sz}, as
well as in~\cite{gj0,gjvai,gjv1,gjv2}.
 
\subsection{Symmetrization} 

The first step of the fundamental transformation is the linear symmetrization operator $\Xi_m$.

\noindent 
\underline{\emph{Definition of} $\Xi_m$ (Step 1)}:
Following~\cite{gjv2}, we define
\begin{equation}\label{eXi}
\Xi_m(p_{\al}z^{|\al |}):=\sum_{\si\in \mathfrak{S}_m}x_{\si(1)}^{\al_1}\cd x_{\si(m)}^{\al_m},
\qquad m\ge 1,
\end{equation}
if $l(\al )=m$ (with $\al=(\al_1,\ld ,\al_m)$), and zero otherwise.

In order to provide more compact expressions when applying $\Xi_m$, some
notation is required.
For $\rho=\{\rho_1,\ld ,\rho_l\}\subseteq\{1,\ld ,m\}$,
 let $x_{\rho}=\{x_{\rho_1},\ld ,x_{\rho_l}\}$.
Given $i_1,\ld ,i_k\geq 1$, with $i_1+\cd +i_k=m$,
let $s_j=i_1+\cd +i_j$, for $j=1,\ld ,k$, and $s_0=0$. Also
let
\begin{equation}\label{eRHOJ}
\rho_j=\{s_{j-1}+1,\ld ,s_j\},
\end{equation}
for $j=1,\ld ,k$.
Then $\SYM_{i_1,\ld,i_k}$ is a summation operator, over
the set of ordered set partitions $(\nu_1,\ld ,\nu_k)$ of $\{1,\ld ,m\}$,
in which, in the summand, $x_{\rho_j}$ is replaced by $x_{\nu_j}$,
for $j=1,\ld ,k$.

We require the symmetrized series
\begin{equation}\label{symrest1}
\left\{
\begin{array}{llllll}
F^g_m(x_1,\ld ,x_m)&:=&\Xi_m F^g(z;\bfp),\qquad
 &\widbfH^0_m(x_1,\ld ,x_m)&:=&\Xi_m\widH^0(z;\bfp ),\\
\zeta^g_m(x_1,\ld ,x_m,u)&:=&\Xi_m \zeta^g(z,u;\bfp),\qquad
&\xi^{(i)}_m(x_1,\ld ,x_m,u)&:=&\Xi_m \xi^{(i)}(z,u;\bfp)),
\end{array}
\right.
\end{equation}
for $m\ge 1$. When the arguments of these series are suppressed, they are the ones stated above.

\subsection{Symmetrizing on the localization and degeneration sides}

On the \textsf{localization side}, 
the following result is the third form of the Localization Tree Theorem. It
is the symmetrized form of Theorem~\ref{branchfg},
and uses the substitution operator
\begin{equation}\label{eLSO}
\La\colon u\mapsto -u,\, x_i\mapsto(1-u)^{-1} x_i,\qquad i=1,\ld,m.
\end{equation}

\begin{co}[Localization Tree Theorem --- symmetrized functional equations]\label{tSFHT4} For $m,g\ge1$,
\begin{equation*}
a)\quad F^g_m \G_{g,1} =[u^{2g-1}]\left(\f{-u}{1-u}\right)^m\!\!\!\!\La
 \zeta^g_m,
\end{equation*}
\newline
where
$$b)\quad \zeta^g_m=
\sum_{n\ge 1}\frac{1}{n!}\!\!\!\!\!\!\!\!
\sum_{\substack{a_1,\ld ,a_n,k\geq 0, \\ a_1+\cdots+a_n+k=g-2+n, \\  i_1,\ld ,i_n\geq 1, \\  i_1+\cdots +i_n=m}}    
\!\!\!\!\!\!\!\!\!\!\!\!\!\!\!\! (-1)^k \laFab\tau_{a_1}\cd\tau_{a_n}\la_k\raFab
\SYM_{i_1,\ld,i_n}
\prod_{j=1}^n \xi^{(a_j)}_{i_j}(x_{\rho_j}).$$
\end{co}

\bpf
The result follows by applying $\Xi_m$ to Theorem~\ref{branchfg}(a) and~(b), for $m\ge 1$, and
using the properties
of $\Xi_m$ given in Lemmas 4.1--4.3 of~\cite{gjvai}.
\epf

On the \textsf{degeneration side}, the following result is
the symmetrized form of Corollary~\ref{JCEforFHS}.
For each $1\le i <j\le m$, note that there are two terms in the
result with denominator $x_i-x_j$, and that the numerator is antisymmetric in $x_i,x_j$,
so these terms combine to give a formal power series. 
To account for each individual term $(x_i-x_j)^{-1}$ that arises in the statement, we
adopt  the total ordering  $x_i \prec x_j$ if $i<j,$ and then define $(x_i - x_j)^{-1}$ by
$(x_i-x_j) ^{-1}:=x_j^{-1} (1-x_i x_j^{-1})^{-1}.$
This gives us an ordered Laurent series ring in $x_1,\ldots,x_m$.

\begin{co}[Degeneration Theorem --- symmetrized Join-cut Equation]\label{tSJCE}
For $m,g\ge 1$,
\begin{eqnarray*}
\left(\sum_{i=1}^m x_i\f{\pa}{\pa x_i}+m-1\right) F^g_m&=&
\!\!\!\!\!\!\sum_{\substack{i_2,i_3\ge 0, \\ i_2+i_3=m-1}}\!\!\!\!\SYM_{1,i_2,i_3}
\left( x_1\f{\pa}{\pa x_1}\widbfH^0_{1+i_2}(x_1,x_{\rho_2})\right)
\left( x_1\f{\pa}{\pa x_1}F^g_{1+i_3}(x_1,x_{\rho_3})\right) \nonumber \\
&{}&+\SYM_{1,1,m-2}\f{x_2}{x_1-x_2}x_1\f{\pa}{\pa x_1}F^g_{m-1}(x_1,x_{\rho_3})
+\sum_{i=1}^m\left( x_i\f{\pa}{\pa x_i}\right)^{2g+1}\widbfH^0_m.
\end{eqnarray*}
\end{co}
\bpf
The result follows by applying $\Xi_m$ to Corollary~\ref{JCEforFHS} for $m\ge 1$, and
using the properties
of $\Xi_m$ given in Lemmas 4.1--4.3 of~\cite{gjvai}.
\epf

Note that the second term on the right hand side of Corollary~\ref{tSJCE} vanishes
in the case $m=1$.

\subsection{Polynomiality and terms of top degree}\label{Spft}

The second step of the fundamental transformation is the change of variables $\sC$.

\noindent  
\underline{\emph{Definition of} $\sC$ (Step 2)}:
Let
\begin{equation}\label{Ty}
a)\quad w(x):=\sum_{n\geq 1}\f{n^{n-1}}{n!}x^n,\qquad\qquad
b)\quad y(x):=\frac{1}{1-w(x)}=\sum_{n\geq 0}\f{n^n}{n!}x^n,
\end{equation} 
where the second equality in the expression for $y(x)$ is
well-known (see, \textit{e.g.},~\cite[Prop.~3.2.1]{gj0}).
Let $w_i:=w(x_i)$ and $y_i=y(x_i)$, for $i=1,\ld ,m$,
\dmjdel{
We shall see that it suffices to work in a particular subring of the formal Laurent series ring 
$\Q[[x_1,\ldots,x_m]](u),$ Let $\sC$ be 
the ring homomorphism defined by
\begin{equation}\label{eOPC}
\sC\colon \Q\left[x_1\f{dw_1}{dx_1},\ldots,x_m\f{dw_m}{dx_m}\right](x)
\rightarrow \Q[y_1,\ldots,y_m](u)
\colon x_i \f{dw_i}{dx_i}\mapsto y_i-1, \quad i=1,\ldots,m.
\end{equation}
where $y_j:=y(x_j)$, $j=1,\ldots ,m.$
} 
and let $\sC$ be an operator, applied to a formal power series in $x_1,\ld ,x_m$,
that changes variables from the indeterminates $x_1,\ld ,x_m$ to $y_1,\ld ,y_m$.
Thus, a direct way to apply $\sC$ is to substitute $x_i=G(y_i-1)$, where $G(z)$ is
the compositional inverse of the formal power series $\sum_{n\ge 1}n^nz^n/n!$. In
this paper, we do \emph{not} apply this substitution directly to the
symmetrized series $F_m^g$, $\La\;\zeta_m^g$,  $\La\;\xi_m^{(i)}$, but instead
quite indirectly. The details are intricate, but the key to our method
is that in all cases the result
is a \emph{polynomial} in $y_1,\ld ,y_m$  (we
say that a Laurent series in
another indeterminate, either $u$ or $t$ in this paper, is
\emph{polynomial} if each of its coefficients is polynomial).

This key fact is recorded in the
following result. The proof involves a number of transformations for
implicitly defined series, and is deferred until Appendix A to avoid
interrupting the present development. Note that 
part~(b) follows from part~(a) by Theorem~\ref{branchfg}(a),(b), and by~(\ref{ianeq}) and
Theorem~\ref{onedim} (which allows us to ``divide'' $\laFab \tau_{a_1}$ $\cd\tau_{a_n}\la_k\raFab$
 by $\G_{g,1}$).
Consequently, in the Appendix it suffices to prove part~(a) of Theorem~\ref{polyfact} only.

\begin{tm}\label{polyfact}
For $m=1,2,3$, we have
$$ a)\quad\sC\;\La\;\xi^{(i)}_m\in  \Q[y_1,\ldots,y_m]((u)),
\quad i\ge 1, \qquad\qquad
 b)\quad\sC\; F_m^g\in \Q[y_1,\ldots,y_m],\quad g\ge 1.$$
\end{tm}

Theorem~\ref{polyfact} is essential to our proof of Faber's Intersection  Conjecture
for at most $3$ parts. In order for us to extend our proof to, say, $n$ parts, we
would first need to prove Theorem~\ref{polyfact} for all $m\le n$. We are presently
able to prove this result for $m\le 5$, but not for larger values of $m$, since
we require the symmetrized double Hurwitz series in our method of proof, and
we only have explicit expressions for this when $m\le 5$, as given in~\cite{gjv2}.
We conjecture that Theorem~\ref{polyfact} holds for all positive integers $m$.

The third and final step of the fundamental transformation is the operator $\Ups'$, that restricts
a polynomial to terms of maximum total degree (these are referred to as the ``top'' terms).

\noindent
\underline{\emph{Definition of} $\Ups'$ (Step 3)}:
Let $\sS_i$ be the operator that restricts a polynomial in $y_1,\ld ,y_m$
to the terms of total degree $i$. Then we define $\Ups'$, when applied to $\sC\; F^g_m$,
to denote $\sS_{4g+3m-5}$. It turns out that there are no terms of higher
degree in $\sC\; F^g_m$, so we say that $\Ups'$ restricts to the
terms of top degree, though we understand that this is informal, since
it assumes that the terms of total degree ${4g+3m-5}$ are not all zero.

If $\sC\;\La\;\xi^{(i)}_m=\sum_k a_k u^k$, we define $\Ups'\; \sC\;\La\;\xi^{(i)}_m$
$:=\sum_k \left( \sS_{2g+2m-2+k}a_k\right) u^k$. Again, it turns out
that there are no terms is any of the $a_k$ of higher total degree, so
in this case also $\Ups'$ restricts to terms of top degree.
We define $\Ups'$, when applied to $\sC\; {\widbfH}^0_m$,
to denote $\sS_{3m-6}$. Again there are no terms of higher
total degree, so in this case also $\Ups'$ restricts to terms of top degree.

Finally, in all other cases, we define $\Ups'$ as a homomorphism, and we
define $\Ups:=\Ups'\;\sC$.

\dmjdel{Let $\Ups'$ be the operator 
$$\Ups'\colon \Q[y_1,\ldots,y_m](u) \rightarrow \Q[y_1,\ldots,y_m](u)$$
such that, for $f=\sum_i a_i(\bfy_m)u^i\in \Q[y_1,\ldots,y_m](u),$
then i)~$\Ups' f=\sum_i\Ups'(a_i(\bfy_m)u^i,$ and 
ii)~$\Ups' a_i(\bfy_m)=[\la^p]a_i(\la\bfy_m),$ where $p=\deg_\la a_i(\la\bfy_m)$
and $\bfy_m:=(y_1,\ldots,y_m).$
Thus $\Ups' a_i(\bfy_m)$ is the sum of terms of maximum (total) degree in $a_i(\bfy_m).$
Note, however, that $\Ups' $ is \emph{not} linear but it is multiplicative on $\Q[y_1,\ldots,y_m].$

Let $\Ups:=\Ups'\sC.$ 

Since $\Ups'$ is not linear, we first need to determine the maximum degree of its operand and then
apply the operator.

We shall  apply $\Ups$ the following operands  by extracting the coefficient of the degree indicated
below. The terms thus obtained will have \emph{maximum} total degree, although we have
to verify that the result is non-zero.  This will be done in each instance of its use.

On the \emph{degeneration side}: \textbf{i)}  for $F_m^g,$ the total degree is 
\underline{$4g+3m-5$}, $m,g\ge 1;$
\textbf{ii)}  for the symmetrized single Hurwitz series, genus $0$,
the total degree is \underline{$3m-6$}, for $m\ge 1$.

On the \emph{localization side}: \textbf{iii)}  for $a_l$ where $\La\zeta_m^g:=\sum_{\ell} a_{\ell}u^{\ell}$, 
the total degree is $\underline{2g+2m-4+\ell};$

\textbf{iv)}  for $a_l,$ where $\xi^{(g)}_m:=\sum_{\ell} a_{\ell}u^{\ell}$,
the \emph{symmetrized localization tree series},  genus $0$,
the  total
degree is \underline{$2g+2m-2+\ell$}.
} 

\subsection{Polynomiality and terms of top degree on the localization and degeneration sides}

On the \textsf{localization side}, we can apply the fundamental transformation
\textit{via} Theorem~\ref{branchfg}, because of the
polynomiality of $\La\;\xi^{(i)}_m$ that was established
in Theorem~\ref{polyfact}(a) for $m=1,2,3$. This requires the top terms, $\Ups\;\La\; \xi^{(i)}_m$,
which are given in the following result for $m=1,2,3$. Again, the proof is intricate,
and is deferred until Appendix A to avoid interrupting the present development.
The result uses the notation 
\begin{equation}\label{Yiu}
Y_i(u):=\f{y_i}{1-uy_i},\;\;\;\;\;\; i\ge 1.
\end{equation}

\begin{tm}\label{xi123}
For $i\ge 0$,
$$a)\quad u\Ups\;\La\; \xi^{(i)}_1=-(2i-1)!!Y_1(u)^{2i+1}, \qquad
b) \quad u^2\Ups\;\La\; \xi^{(i)}_2
=-(2i+1)!!u\,\SYM_{1,1}\f{y_1^2y_2}{y_1-y_2} Y_1(u)^{2i+3},$$
\begin{eqnarray*}
c)\quad\f{(-u)^3\Ups\;\La\; \xi^{(i)}_3}{(2i+1)!!}
&=&u^2\SYM_{1,1,1}\f{y_1^3 y_2^4 y_3}{(y_2-y_3)(y_1-y_2)^2}
Y_1(u)^{2i+3}Y_2(u)\\
&{}&-\SYM_{1,2}\,
\left( uy_1Y_1(u)^{2i+5}Y_2(u)Y_3(u)-u^2y_1^3\f{\pa}{\pa y_1}\f{y_1^3 y_2 y_3}{(y_2-y_1)(y_3-y_1)}Y_1(u)^{2i+4}\right) .
\end{eqnarray*}
\end{tm}

On the \textsf{degeneration side}, we can apply the fundamental transformation
\textit{via} Corollary~\ref{tSJCE} if we can apply the change of variables $\sC$ to
the partial differential operator $x_i\f{\pa}{\pa x_i}$ and to the symmetrized
generating series ${\widbfH}^0_m$ for Hurwitz numbers, genus $0$. 
This is straightforward in both cases.

For the partial differential operator, from ~(\ref{Ty}a,b) we have 
the functional equations $y_i=1/(1-w_i)$, $y_i=1+x_i\f{dw_i}{dx_i}$, and
together these imply $x_i\f{dy_i}{dx_i}=y_i^2(y_i-1)$. This immediately
implies the operator identities
\begin{equation}\label{Upident}
a)\quad\sC\, x_i\f{\pa}{\pa x_i}=y_i^2(y_i-1)\f{\pa}{\pa y_i}\sC,\qquad\qquad
b)\quad\Ups\, x_i\f{\pa}{\pa x_i}=y_i^3\f{\pa}{\pa y_i}\Ups.
\end{equation}

For the symmetrized generating series ${\widbfH}^0_m$, 
~(\ref{singHur}) and~(\ref{Ty}b) give
\begin{equation}\label{singexp}
 {\widbfH}^0_m=\left(\sum_{j=1}^m x_j\f{\pa}{\pa x_j}\right)^{m-3}
\prod_{i=1}^m (y_i-1),
\end{equation}
for $m\ge 1$ (note that~(\ref{Upident}) applied to ~(\ref{singexp}) immediately
identifies $\sC\; {\widbfH}^0_m$ as having degree $3m-6$, consistent with
the definition of $\Ups'$ above). Also, we have
\begin{equation}\label{H0singone}
a)\quad x_1\f{\pa}{\pa x_1}{\widbfH}^0_1 = 1-y_1^{-1},\qquad\qquad b)\quad
{\widbfH}^0_2=\log\left(\f{y_1-y_2}{y_1y_2(x_1-x_2)}\right) +y_1^{-1}+y_2^{-1}-2,
\end{equation}
where the first of these expressions arises by
applying $x_1\f{d}{dx_1}$ to~(\ref{singexp}) when $m=1$, and using~(\ref{Upident}a),
and the second is given in~\cite[p.~38]{gjvai}.

\section{Intersection numbers with one part}\label{Sinop}
In this section, we use our strategy for the first time. This
involves applying the fundamental transformation to give results
for the symmetrized Faber-Hurwitz series $F_m^g$ on both
the \textsf{degeneration} and \textsf{localization sides},
to obtain results
for Faber's intersection numbers. Here, these results are given in Section~\ref{in1p} for the
case of a single ($m=1$) part. Faber's Intersection Number Conjecture,
giving $\laFab\tau_{g-1}\raFab$ $ =\psi_1^{g-1}$, is
immediate in this case, but it gives us an equation relating the two generators
 $\psi_1^{g-1}$ and $\G_{g,1}$ of $A_{2g-1}(\cm_{g,n}^{rt})$. There
is non-trivial geometric information to be gained by relating these generators, which is
described in Section~\ref{sigc}.

\subsection{Intersection numbers with one part}\label{in1p}
We begin with a result for the \textsf{localization side}.
\begin{tm}\label{LHS1}
For $g\ge 1$, we have
$\Ups\, F^g_1(x_1)\G_{g,1}=
(2g-3)!!{4g-3\choose 2g-1}y_1^{4g-2} \psi_1^{g-1}.$
\end{tm}

\bpf
{}From Corollary~\ref{tSFHT4}(b) we have
$\zeta_1^g= \sum_{k=0}^{g-1} (-1)^k \laFab \tau_{g-1-k} \la_k \raFab\xi_1^{(g-1-k)}.$
Together with Theorem~\ref{xi123}(a), this gives
\begin{eqnarray*}
[u^{2g-1}](-u)^1\Ups\,\La\,
\zeta^g_1(x_1,u) = 
[u^{2g-1}]\laFab\tau_{g-1}\raFab
(2g-3)!!Y_1(u)^{2g-1}
=[u^{2g-1}](2g-3)!!Y_1(u)^{2g-1} \psi_1^{g-1}
\end{eqnarray*}
where, for the second equality, we have used
the immediate 
fact that
 $\laFab\tau_{g-1}\raFab =\psi_1^{g-1}$ for  $g\ge 1.$
But
$Y_1(u)^{2g-1}=\sum_{i\ge 0}{-(2g-1)\choose i}(-1)^iy_1^{2g-1+i}u^i
=\sum_{i\ge 0}{2g-2+i\choose i}y_1^{2g-1+i}u^i,$
and the result follows immediately.
\epf

Now we turn to the corresponding result for the \textsf{degeneration side}.

\begin{tm}\label{onesymexp}
For $g\ge 1$, we have
$\Ups\, F^g_1(x_1)=
(4g-3)!!\f{y_1^{4g-2}}{4g-2}.$
\end{tm}

\bpf
Let $m=1$ in Corollary~\ref{tSJCE}, to
obtain $x_1\f{\pa}{\pa x_1}F^g_1=\left( x_1\f{\pa}{\pa x_1}\widbfH^0_1\right)
\left( x_1\f{\pa}{\pa x_1}F^g_1\right)$
$+\left(x_1\f{\pa}{\pa x_1}\right)^{2g+1}\widbfH^0_1$.
Solving for  $x_1\f{\pa}{\pa x_1}F^g_1$, and using~(\ref{H0singone}a), we obtain
\begin{equation}\label{F1x1}
x_1\f{\pa}{\pa x_1}F^g_1=y_1
\left(x_1\f{\pa}{\pa x_1}\right)^{2g-1}\!\!\!\!\!\!\!\!(y_1-1).
\end{equation}
The polynomiality of $F^g_1$ follows immediately from~(\ref{Upident}a) (or we
can use Theorem~\ref{polyfact}(b), with $m=1$), so
from~(\ref{Upident}b) we have
$y_1^3\f{\pa}{\pa y_1}\Ups\, F^g_1$ $=y_1\left(y_1^3\f{\pa}{\pa y_1}\right)^{2g-1}\!\!\!\!\!\!y_1$
$=(4g-3)!!\,y_1^{4g}$,
where the second equality follows by induction. The result follows by integrating in $y_1$.
\epf

\subsection{Some immediate geometric consequences}
\label{sigc}\lremind{sigc}
Comparing the above two results,  
we obtain the following result, which relates the generator $\G_{g,1}$ of $A_{2g-1}(\cm_{g,n}^{rt})$ with 
the generator $\psi_1^{g-1}$ suggested by Faber's Intersection Number Conjecture.
\begin{co}\label{cgformula}
For $g\geq 1$,
$\;\G_{g,1} = \f{2^g}{(g-1)!} \psi_1^{g-1}$
on $A_{2g-1}(\cm_{g,1})$.
\end{co}

In particular, $\psi_1^{g-1} = \laFab \tau_{g-1} \raFab$
is another non-zero element (=basis) of the one-dimensional vector 
space $A_{2g-1}(\cm_{g,1})$, and this describes the change of basis
from $\G_{g,1}$ to $\psi_1^{g-1}$.

\bpf
Comparing Theorems~\ref{LHS1} and~\ref{onesymexp}, we obtain
$(2g-3)!!{4g-3\choose 2g-1} \psi_1^{g-1} =
\f{(4g-3)!!}{4g-2} \G_{g,1},$
and the result follows immediately.
\epf

\begin{co}[{\cite[Thm.~4]{phi}}]
$\qquad\sum_{i=0}^{g-2} (-1)^i \la_i \psi_1^{g-1-i} = 2^{g-1} \psi_1^{g-1}/g!$
\end{co}

\bpf
Applying the usual (non-relative) virtual localization formula to
$\G_{g,1,\sim}$, we obtain
$\G_{g,1,\sim} = \psi_1^{2g-1} - \psi_1^{2g-2} \la_1
+ \cdots + (-1)^g \psi_1^{g-1} \la_g$.  
The result then follows from 
Corollary~\ref{cgformula} and the fact that
$\G_{g,1}= 2g \G_{g,1,\sim}$ (Proposition~\ref{bigbusiness}).
\epf

\begin{co}
The class of the hyperelliptic curves in $\cm_{g}$ is\label{hyp}\lremind{hyp}
$\;\frac { (2^{2g}-1) 2^{g-1}}   {g! (2g+1)(2g+2)} \ka_{g-2}.$
\end{co}

(This was stated, for example, in the concluding remarks of
\cite{fconj}.)
The argument carries through for the locus of degree $d$ covers
fully ramified over two points.

\noindent {\em Sketch of proof.}
The class $\G_{g,2,\sim}$ corresponds to the points $[C,p,q] \in \cm_{g,2}^{rt}$
where $\oh_C(p-q)$ is a $2$-torsion point (along with a 
virtual class).    The locus where $\oh_C(p-q)$ is trivial corresponds
to $\G_{g,1}$.  Thus the locus where $\oh_C(p-q)$ is $2$-torsion but non-trivial
has virtual class  $(2^{2g}-1) \G_{g,1}$ by Theorem~\ref{trick}.
It is straightforward to check the virtual fundamental class on this locus is the actual
fundamental class (as mentioned in the proof of Theorem~\ref{trick}).
This locus corresponds  to hyperelliptic curves with two marked Weierstrass points.
Thus the
locus in $\cm_{g,1}$ of hyperelliptic curves with choice of {\em one} Weierstrass
point has class 
$\frac { 2^{2g}-1 } {2g+1} \G_{g,1}
= \frac { (2^{2g}-1) 2^{g-1}}   {g! (2g+1)} \psi_1^{g-1}.$
Pushing forward to $\cm_g$, and forgetting the choice
of one of the $2g+2$ Weierstrass points, we get the desired result. \epf

\noindent
\begin{co}\label{Fnice}\lremind{Fnice}
For $d,g\ge 1$, we have
$F^g_{(d)}=\f{1}{d}\sum_{i=1}^d {d \choose i} i^{2g+i-1}(d-i)^{d-i}.$
\end{co}
This is a corollary of the proof
of Theorem~\ref{onesymexp}, not of the result itself. 

\bpf
{}From~(\ref{F1x1}) and~(\ref{Ty}b), we have
$x_1\f{\pa}{\pa x_1}F^g_1=\sum_{n\ge 0}\f{n^n}{n!}x_1^n
\sum_{i\ge 1}\f{i^{2g+i-1}}{i!}x_1^i,$
and the result follows from~(\ref{symrest1}) and~(\ref{gsFnum}).
\epf

Thus, for example, the class of non-singular genus $g$ curves $C$ admitting 
a degree $d$ cover of $\proj^1$ \textit{via} $\oh(d p)$ for some $p \in C$,
ramified over an appropriate number of points is
$$
\f{1}{d}\sum_{i=1}^d {d \choose i} i^{2g+i-1}(d-i)^{d-i} \G_{g,1} = 
\f{1}{d}\sum_{i=1}^d {d \choose i} i^{2g+i-1}(d-i)^{d-i} 
\frac {2^g \ka_{g-2}}
{ (g-1)!} \in A_{2g-1}(\cm_g).
$$

\section{A strategy  for an inductive proof of Faber's Intersection Number Conjecture}
\label{Safhst}

In this section, we consider an inductive strategy for proving the Conjecture.

\subsection{Final form of the Localization Tree Theorem}
First, we define two generating series, in which the
genus~$g$ is marked by the indeterminate~$t$.

\noindent\underline{\em The localization tree generating series} is
\begin{eqnarray}\label{Psigf}
\Psi_m(y_1,\ld ,y_m,t) &:=&\sum_{g\ge 1}\f{t^{2g}}{(2g-1)!!}
[u^{2g-1}](-u)^m\Ups\,\La\,\zeta_m^g(x_1,\ld ,x_m,u)/\psi_1^{g-1}, \;\;\;\;\;\; 
m\ge 1.
\end{eqnarray}

\noindent\underline{\em The Faber generating series} is
\begin{eqnarray}\label{Phigf}
\Phi_m(y_1,\ld ,y_m,t)&:=&\sum_{g\ge 1}2^{2g-1}\f{t^{2g}}{(2g-1)!}
\Ups\,F_m^g(x_1,\ld ,x_m),\;\;\;\;\;\; m\ge 1.  
\end{eqnarray}

For $\Psi_m$, 
thanks to the fact that $R_{2g-1}(\cm_{g,1})$ $ \cong \Q$ (Thm.~\ref{onedim}),
we ``divide by $\psi_1^{g-1}$'' using the isomorphism
$$
\xymatrix{ \Q \ar[rr]^{\times \psi_1^{g-1}} & & R_{2g-1}(\cm_{g,1})}
$$
(combining \eqref{divide} and Cor.~ \ref{cgformula}), and thus
we define the {\em Faber intersection number}
\begin{equation}\label{Fabrational}
\langle \cdots \rangle^{\Faber}_g := \laFab \cdots \raFab / \psi_1^{g-1}
=  \laFab \cdots \raFab / \laFab \tau_{g-1} \raFab .
\end{equation}

Note that $\Psi_m$ and $\Phi_m$, since they require the operator $\Ups$, are
only defined when the operator $\sC$ yields a polynomial. So far, we have
proved this polynomiality only for $m=1,2,3$, in Theorem~\ref{polyfact}.  
Subject to
this, $\Psi_m$ and $\Phi_m$  are both
generating series with rational coefficients, and 
the next result shows that they are equal.
\begin{co}[Localization Tree Theorem --- for genus generating series] \label{PsiPhieq}
For $m\ge 1$,
$$\Psi_m(y_1,\ld ,y_m,t)=\Phi_m(y_1,\ld ,y_m,t).$$
\end{co}
\bpf
The result follows  from Corollary~\ref{tSFHT4}(a) by noting that 
$\f{\G_{g,1}}{(2g-1)!!}=\f{2^g \psi_1^{g-1}}{(g-1)!(2g-1)!!}=\f{2^{2g-1} \psi_1^{g-1}}{(2g-1)!}$
for $g\ge 1,$ from Corollary~\ref{cgformula}.
\epf

Corollary~\ref{PsiPhieq} is the final form of Theorem~\ref{bigyuck}.
It summarizes  the relationship between
the \textsf{localization side} (dealing with $\Psi_m$) and 
the \textsf{degeneration side} (dealing with $\Phi_m$) 
as the equality of~$\Psi_m$ and~$\Phi_m$.
Our strategy for exploiting this relationship in order to prove
the Faber Conjecture for a fixed number of parts is indirect. We prove
in the following lemma that Corollary~\ref{PsiPhieq} implies
sufficiently many linearly independent linear equations for Faber's
intersection numbers to uniquely identify them. This means that
we can prove the Faber Conjecture iteratively on the number of parts,
by simply verifying the equality $\Psi_m=\Phi_m$ of Corollary~\ref{PsiPhieq}
with Faber's conjectured values substituted for the intersection numbers.
In fact,
as we shall see, with these substituted values, we are able to verify the
equality $\Psi_m=\Phi_m$ for particular $m$ simply by equating
polynomials.

\begin{lm}\label{nonsing}
For any $n\ge 2$, if the Faber Conjecture is true for up to $n-1$ parts, and 
$$\Psi_n(y_1,\ld ,y_n,t)=\Phi_n(y_1,\ld ,y_n,t)$$
with Faber's conjectured values substituted for the intersection numbers with up to $n$ parts,
then the Faber Conjecture is true for up to $n$ parts. 
\end{lm}

\bpf
Consider the equality
$\Psi_n(y_1,\ld ,y_n,t)=\Phi_n(y_1,\ld ,y_n,t),$
and equate coefficients of $t^{2g}$ for each fixed $g\ge 2$.
This gives an equation involving Faber's intersection numbers
with at most $n$ parts, and from Corollary~\ref{tSFHT4}(b),~(\ref{Psigf}) and
Theorem~\ref{xi123}(a), we deduce that the only terms involving
the intersection numbers
with exactly $n$, positive, parts are given by
\begin{equation}\label{linthr}
\left\{
\begin{array}{cl}
{}&\sum_{\substack{a_1,\ld ,a_n\ge 1, \\ a_1+\cd +a_n=g+n-2}}
\!\!\!\!\langle \tau_{a_1}\cd \tau_{a_n}\rangle^\Faber_g
\f{1}{(2g-1)!!}[u^{2g-1}]
\prod_{j=1}^n(2a_j-1)!!Y_j(u)^{2a_j+1}\\
=&\sum_{\substack{a_1,\ld ,a_n\ge 1, \\ a_1+\cd +a_n=g+n-2}}
\;\;
\sum_{\substack{m_1,\ld ,m_n\ge 0, \\ m_1+\cd +m_n=2g-1}}
\!\!\!\!\langle \tau_{a_1}\cd \tau_{a_n}\rangle^\Faber_g
\f{\prod_{j=1}^n(2a_j-1)!!{2a_j+m_j \choose m_j}y_j^{2a_j+m_j+1} }{(2g-1)!!}.
\end{array}
\right.
\end{equation}
Now under the hypothesis that the Faber Conjecture has been proved for up to  $n-1$ parts,
we are able to evaluate all Faber intersection numbers in the above equation
with up to $n-1$ parts, as well as those with $n$ parts, at least
one of which is zero (by use of the string equation), and we consider
Faber's intersection numbers with exactly $n$, positive, parts
as unknowns.
Now equate coefficients of $y_1^{i_1}\cd y_n^{i_n}$ in this equation
for each $i_1,\ld ,i_n \ge 3$ with $i_1+\cdots +i_n=4g-5+3n$, to obtain
a linear equation that we shall refer to as $B_{i_1,\ld ,i_n}$. 
Then we thus obtain a system
of ${4g-6+n\choose n-1}$ equations $B_{i_1,\ld ,i_n}$ in
the ${g-3+n\choose n-1}$ ``unknowns'' $\langle \tau_{a_1}\cdots \tau_{a_n}\rangle^\Faber_g$,
where $a_1, \ld ,a_n\ge 1$, with $a_1+ \cdots +a_n=g-2+n$. Let the 
coefficient matrix for this system be denoted by $\bfM$,
and note that ${4g-6+n\choose n-1}\ge {g-3+n\choose n-1}$, since $g,n\ge 2$.
(We will not use the symmetry of these  unknowns.)

Now we prove that this
system has at most one solution.
To do so,
suppose that we consider the unknowns in lexicographic order of
the $a_1\ld a_n$, and the equations in lexicographic order of
the $i_1\ld i_n$ (so we have ordered the rows and columns of $\bfM$.)
Note that, from~(\ref{linthr}), the entry in row $i_1\ld i_n$ and
column $a_1\ld a_n$ is non-zero if and only if $i_j\ge 2a_j+1$ for
all $j=1,\ld ,n$.
Then it is straightforward to check that the submatrix of $\bfM$ consisting
of all columns, and rows $i_1,\ld ,i_n$ with $i_1,\ld ,i_{n-1}\ge 3$ and odd, and
 $i_1+\cd +i_{n-1}\le 2g-7+3n$, is lower triangular, with non-zero coefficients on
the diagonal. Thus these rows of $\bfM$ are linearly independent, so
the system has at most one solution. The result follows since, if Faber's
conjectured values satisfy this system, then they are uniquely the
correct values for these intersection numbers.
\epf

\subsection{A lemma for the localization tree generating series}
We consider first the \textsf{localization side}.
In order to determine a compact rational expression for the localization tree generating
series $\Psi_m$, we shall need the following lemma, where
$\oE f$ denotes the {\em even} subseries
of the formal power series $f$ in the indeterminate $t$.
Also, for $i\ge 1$ we let
\begin{equation}\label{AiBi}
B_i:=\sqrt{1-4y_i^2t},\qquad\qquad\qquad A_i:=\f{1-B_i}{2y_i}.
\end{equation}

\begin{lm}\label{lagY1u}
For a formal power series $f$,
$\sum_{g\ge 1} t^{2g}\left[ u^{2g-1}\right] f(u) Y_1(u)^{2g-1}
=\tfrac{1}{2}\oE\, t\left( 1+B_1^{-1}\right) f(A_1).$
\end{lm}
\bpf      
Now
$\left[ u^{2g-1}\right] f(u) Y_1(u)^{2g-1}=
\left[ t^{2g-1}\right] \f{f(s)}{1-t\f{d}{ds}Y_1(s)},$
from Lagrange's Implicit Function Theorem (see, \textit{e.g.},~\cite[Thm.~1.2.4]{gjce}), where~$s$ is
the unique formal power series  solution 
(in~$t$) of  $s=tY_1(s)$,  which is the equation $y_1s^2-s+y_1t=0$, and so~$s=A_1$.
But $\f{d}{ds}Y_1(s)=Y_1(s)^2$, so with $s=A_1$,  we have
$1-t\f{d}{ds}Y_1(s)$ $=\f{2B_1}{1+B_1}.$
Then
$\left[ u^{2g-1}\right] f(u) Y_1(u)^{2g-1}= \left[ t^{2g-1}\right]\tfrac{1}{2}\left( 1+B_1^{-1}\right) f(A_1).$
The result follows.
\epf

\subsection{Final form of the Degeneration Theorem}\label{ssFGSJCE}
For the \textsf{degeneration side}, we describe the general results
that will allow us to obtain an explicit rational expression
for the Faber generating series $\Phi_m$.
First we give a compact expression for $\Phi_1$.

\begin{co}\label{Phione}
$\qquad y_1\f{\pa}{\pa y_1}\Phi_1(y_1,t)=\oE\, tB_1^{-1}.$
\end{co}

\bpf
{}From~(\ref{Phigf}) with $m=1$, and Theorem~\ref{onesymexp}, we obtain
\begin{equation*}
y_1\f{\pa}{\pa y_1}\Phi_1(y_1,t)
=\sum_{g\ge 1}2^{2g-1}\f{(4g-3)!!}{(2g-1)!}y_1^{4g-2}t^{2g}
=t\sum_{g\ge 1}{-\tfrac{1}{2}\choose 2g-1}\left(-4y_1^2t\right) ^{2g-1},
\end{equation*}
and the result follows immediately.
\epf

In the next result, we consider the top terms in Corollary~\ref{tSJCE},
to  obtain a simple iterative equation for the
Faber generating series $\Phi_m$, $m\ge 2$.
Similarly to Corollary~\ref{tSJCE}, there are two terms in this
result with denominator $y_i-y_j$, and the numerator is antisymmetric in $y_i,y_j$,
so these terms combine to give a formal power series.
Also, to account for each individual term $(y_i-y_j)^{-1}$, we
adopt the total ordering  $y_i \prec y_j$ if $i<j$, to obtain
an ordered Laurent series ring in $y_1,\ldots,y_m$.

The statement of the result uses the operator
\begin{equation}\label{eDOp}
\De_{k}:=\sum_{i\ge 1} y_i^k\f{\pa}{\pa y_i},\;\;\;\;\;\;\;\;
k\ge 2,
\end{equation}
and the generating series
\begin{equation}\label{ser2g}
\Upsilon_m(y_1,\ld ,y_m,t):=\sum_{g\ge 1}2^{2g-1}\f{t^{2g}}{(2g-1)!}\Ups\sum_{i=1}^m
\left(x_i\f{\pa}{\pa x_i}\right)^{2g+1}\widbfH^0_m(x_1,\ld ,x_m),\;\;\;\;\;\; m\ge 1.
\end{equation}

\begin{lm}[Degeneration Theorem --- Join-cut equation for $\Phi_m$]\label{symFjcutm}
For $m\ge 2$ we have
\begin{eqnarray*}
\De_{2}\Phi_m(y_1,\ld ,y_m,t)&=&\sum_{\substack{i_2\ge 2,i_3\ge 0, \\ i_2+i_3=m-1}}
\SYM_{1,i_2,i_3}
\left( y_1^3\f{\pa}{\pa y_1}\Ups\,\widbfH^0_{1+i_2}(x_1,x_{\rho_2})\right)
\left( y_1^3\f{\pa}{\pa y_1}\Phi_{1+i_3}(y_1,y_{\rho_3},t)\right)\\
&{}&+\SYM_{1,1,m-2}\f{y_1^5y_2}{y_1-y_2}\f{\pa}{\pa y_1}\Phi_{m-1}(y_1,y_{\rho_3},t)
+\Upsilon_m(y_1,\ld ,y_m,t).
\end{eqnarray*}
\end{lm}

\bpf
{}From the symmetrized Join-cut equation~(Cor.~\ref{tSJCE}), and 
moving the contribution for $i_2=0$ from the summation on the right hand side to
the left hand side, and using~(\ref{H0singone}a) and~(\ref{Upident}a), we have
$\left(\sum_{i=1}^m  \f{x_i}{y_i}\f{\pa}{\pa x_i}+m-1\right) F^g_m$
on the left hand side.
But, from~(\ref{H0singone}b) and~(\ref{Upident}a), we have
\begin{equation*}
x_1\f{\pa}{\pa x_1}\widbfH^0_{2}=
\f{y_1^2(y_1-1)}{y_1-y_2}-y_1(y_1-1)-(y_1-1)-\f{x_1}{x_1-x_2}
=  \f{y_1^2(y_2-1)}{y_1-y_2}-\f{x_2}{x_1-x_2},
\end{equation*}
and we apply this result to combine the contribution for $i_2=1$ in
the summation on the right hand side with the $\SYM_{1,1,m-2}$ term on the right hand side,
to give
\begin{eqnarray*}
&{}&\sum_{\substack{i_2\ge 2,i_3\ge 0, \\ i_2+i_3=m-1}}\SYM_{1,i_2,i_3}
\left( x_1\f{\pa}{\pa x_1}\widbfH^0_{1+i_2}(x_1,x_{\rho_2})\right)
\left( x_1\f{\pa}{\pa x_1}F^g_{1+i_3}(x_1,x_{\rho_3})\right)\\
&{}&+\SYM_{1,1,m-2}\f{y_1^2(y_2-1)}{y_1-y_2}x_1\f{\pa}{\pa x_1}F^g_{m-1}(x_1,x_{\rho_3})
+\sum_{i=1}^m\left( x_i\f{\pa}{\pa x_i}\right)^{2g+1}\widbfH^0_m,
\end{eqnarray*}
on the right hand side. Now, multiply on both sides by $2^{2g-1}\f{t^{2g}}{(2g-1)!}$, sum
over $g\ge 1$, and apply $\Ups$, to obtain the result,
from~(\ref{Phigf}) and~(\ref{Upident}b).
\epf

Lemma~\ref{symFjcutm} is the  final form of Theorem~\ref{Fcutjoin}.
In order to apply this result, we require a technical result
about the invertibility of the partial differential operator $\De_{k}$.

\begin{pr}\label{invDe}
For $m\ge 1,k\ge 2$, the operator $\De_{k}$ is invertible
for formal power series in $y_1,\ld ,y_m$ in which every
monomial has positive exponent for each of $y_1,\ld ,y_m$.
\end{pr}
\bpf
Let $A:=\sum_{n_1,\ld ,n_m\ge 1}a_{(n_1,\ld ,n_m)}y_1^{n_1}\cd y_m^{n_m}$
 and $B:=\sum_{n_1,\ld ,n_m\ge 1}b_{(n_1,\ld ,n_m)}y_1^{n_1}\cd y_m^{n_m}$,
and suppose that $\De_{k}A=\De_{k}B$.
We prove that $a_{(n_1,\ld ,n_m)}$ $=b_{(n_1,\ld ,n_m)}$ for all $n_1,\ld ,n_m\ge 1$ by
putting a partial order $\prec_{\ell}$ on the $(n_1,\ld ,n_m)$: we let $(n_1,\ld ,n_m)_{\le}$ be
the vector of length $m$ containing the $n_i$ in weakly increasing order,
and we define $(n_1,\ld ,n_m)\prec_{\ell}\,(i_1,\ld ,i_m)$  if and only
if $(n_1,\ld ,n_m)_{\le}$ precedes $(i_1,\ld ,i_m)_{\le}$ in
lexicographic order.

To prove the result for $(n_1,\ld ,n_m)$, suppose
that $n_i=\max\{n_1,\ld ,n_m\}$ (make an arbitrary choice of $i$ if
there is more than one, equal, maximum element), and equate coefficients
of $y_1^{n_1}\cd y_i^{n_i+k-1}\cd y_m^{n_m}$ in $\De_{k}A=\De_{k}B$,
to obtain
\begin{equation}\label{indrec}
\left\{
\begin{array}{cl}
{}&n_ia_{(n_1,\ld ,n_m)}+\sum_{\substack{j\ne i, \\ n_j\ge k}}
\left(n_j-k+1\right)a_{(n_1,\ld ,n_j-k+1,\ld ,n_i+k-1,\ld ,n_m)}\\
=&n_ib_{(n_1,\ld ,n_m)}+\sum_{\substack{j\ne i, \\ n_j\ge k}}
\left(n_j-k+1\right)b_{(n_1,\ld ,n_j-k+1,\ld ,n_i+k-1,\ld ,n_m)}.
\end{array}
\right.
\end{equation}
Now consider any fixed linear extension of $\prec_{\ell}$,
and note that, for any term in the summations in~(\ref{indrec}), we have
$(n_1,\ld ,n_j-k+1,\ld ,n_i+k-1,\ld ,n_m)\prec_{\ell}\, (n_1,\ld ,n_m).$
The result
follows from~(\ref{indrec}), by induction  on the position in this linear extension.
\epf

\section{Proof of Faber's Intersection Number Conjecture for 2 and 3 parts, with comments
about the general case of  $m$ parts}\label{sFTINC23Pm}

\subsection{Faber's Intersection Number Conjecture for two parts}\label{Sin2p}
In this section, we prove Faber's Conjecture for the case of two parts.  
It will be convenient to ``smooth'' the result to allow zeros, by
applying the string equation.  For cases of~(\ref{Fabconj}) when some
$d_i=0$, we can use the string equation to deduce a consistent
value. For example, when $n=2$, $d_1=g$, $d_2=0$, the string equation
gives
$\langle \tau_g\tau_0\rangle^\Faber_g=\langle \tau_{g-1}\rangle^\Faber_g=1$,
which agrees with~(\ref{Fabconj}), so Faber's Conjecture
for $n=2$ becomes
\begin{equation}\label{fabtwo}
\langle \tau_{d_1}\tau_{d_2}\rangle^\Faber_g=\f{(2g-1)!!}{(2d_1-1)!!(2d_2-1)!!},
\;\;\;\;\;\; d_1+d_2=g,\;\; d_1,d_2\ge 0,\;\; g\ge 1.
\end{equation}

For the \textsf{localization side}, we have the following result. In the
proof, we use the notation
\begin{equation}\label{whYA}
\whY_m:=Y_m(A_1),\;\;\;\;\;\;\;\; m\ge 1,
\end{equation}
where $Y_m$, $A_1$ are defined in~(\ref{AiBi}), (\ref{Yiu}), respectively.
\begin{tm}\label{LHS2}
If the Faber Conjecture is true for $n=2$, then
\begin{equation*}
\Psi_2(y_1,y_2,t)=\oE\, \SYM_{1,1}
\left(\f{y_1^3y_2t}{y_1-y_2}-\f{y_1^4y_2^2t}{(y_1-y_2)(y_1^2-y_2^2+4y_1^2y_2^2t)}\right)B_1^{-1}.
\end{equation*}
\end{tm}

\bpf
{}From Corollary~\ref{tSFHT4}(b) we have
$ \zeta_2^g=
\sum_{\substack{a_1,a_2\ge 0, \\ a_1+a_2\le g}}
 (-1)^{g-a_1-a_2} 
\langle \tau_{a_1}\tau_{a_2} \la_{g-a_1-a_2} \raFab \xi_1^{(a_1)}(x_1)$ $\xi_1^{(a_2)}(x_2)$
$+\sum_{k=0}^{g-1}$ $ (-1)^k\langle \tau_{g-1-k} \la_k \raFab \xi_2^{(g-1-k)}$, 
where we have used the symmetry of the Faber symbol in the first summation.
Together with Theorem~\ref{xi123}(a),(b) and~(\ref{Psigf}), this gives
\begin{eqnarray*}
\Psi_2(y_1,y_2,t)&=&\sum_{g\ge 1}\f{t^{2g}}{(2g-1)!!}
[u^{2g-1}]\sum_{\substack{a_1,a_2\ge 0, \\ a_1+a_2=g}}
\langle\tau_{a_1}\tau_{a_2}\rangle^\Faber_g
\prod_{j=1}^2 (2a_j-1)!!
Y_j(u)^{2a_j+1}\\
&{}&-\sum_{g\ge 1}\f{t^{2g}}{(2g-1)!!}
[u^{2g-1}]\langle \tau_{g-1}\rangle^\Faber_g (2g-1)!!
u\SYM_{1,1}\f{y_1^2y_2}{y_1-y_2}
Y_1(u)^{2g+1}.
\end{eqnarray*}
Now, if the Faber Conjecture is true for $n=2$, then
from~(\ref{fabtwo}), we have
\dmjdel{
\begin{eqnarray*}
\Psi_2(y_1,y_2,t)&=&\sum_{g\ge 1}t^{2g}
[u^{2g-1}]\sum_{\substack{a_1,a_2\ge 0, \\ a_1+a_2=g}}
Y_1(u)^{2a_1+1}Y_2(u)^{2a_2+1}\\
&{}&-\sum_{g\ge 1}t^{2g}
[u^{2g-1}]\SYM_{1,1}\f{y_1^2y_2}{y_1-y_2}
uY_1(u)^{2g+1},
\end{eqnarray*}
} 
\begin{eqnarray*}
\Psi_2(y_1,y_2,t) =
\sum_{g\ge 1}t^{2g}[u^{2g-1}] \left(
\sum_{\substack{a_1,a_2\ge 0, \\ a_1+a_2=g}} Y_1(u)^{2a_1+1}Y_2(u)^{2a_2+1}
-  
\SYM_{1,1}\f{y_1^2y_2}{y_1-y_2} uY_1(u)^{2g+1}\right),
\end{eqnarray*}
But
$\sum_{\substack{a_1,a_2\ge 0, \\ a_1+a_2=g}}
Y_1(u)^{2a_1+1}Y_2(u)^{2a_2+1}=
\SYM_{1,1}\f{Y_1(u)^{2g+3}Y_2(u)}{Y_1(u)^2-Y_2(u)^2},$
and from Lemma~\ref{lagY1u} we conclude that
\begin{equation*}
\Psi_2(y_1,y_2,t)=\tfrac{1}{2}\oE\, \SYM_{1,1}
t\left( 1+B_1^{-1}\right)
\left(\f{\whY_1^4\whY_2}{\whY_1^2-\whY_2^2}\right.
-\left. \f{y_1^2y_2}{y_1-y_2}A_1\whY_1^2\right).
\end{equation*}
{}From the proof of Lemma~\ref{lagY1u}, we have $A_1=ty_1/(1-y_1A_1)$, so
routinely simplifying, we obtain
\begin{eqnarray*}
\Psi_2(y_1,y_2,t)&=&\oE\, \SYM_{1,1}
\left(\f{y_1^3y_2tA_1(1-y_2A_1)B_1^{-1}}{(y_1-y_2)(y_1+y_2B_1)}\right.
-\left.\f{y_1^3y_2A_1^2B_1^{-1}}{y_1-y_2}\right)\\
&=&\oE\, \SYM_{1,1}\left(\f{y_1^3y_2tB_1^{-1}}{y_1-y_2}\right.
+\left.\f{y_1^3y_2^3t-y_1^4y_2^2tB_1^{-1}}{(y_1-y_2)(y_1^2-y_2^2B_1^2)}\right).
\end{eqnarray*}
Now it is straightforward to verify that
$\oE\, \SYM_{1,1}\f{y_1^3y_2^3t}{(y_1-y_2)(y_1^2-y_2^2B_1^2)}=0,$
and the result follows with $y_1^2-y_2^2B_1^2=y_1^2-y_2^2+4y_1^2y_2^2t$.
\epf

For the \textsf{degeneration side}, we have the following result.
\begin{tm}\label{RHS2}
$\De_{3}\De_{2}\Phi_2(y_1,y_2,t)=
\oE\,\SYM_{1,1}
\left(\De_{3}\f{y_1^4y_2}{y_1-y_2}tB_1^{-1}+3y_1^5y_2tB_1^{-5}\right).$
\end{tm}

\bpf
We begin by using~(\ref{singexp}) to obtain
$\De_{3}\Ups\,{\widbfH}^0_2=y_1y_2,$
so, from~(\ref{Upident}b), we have
$ \Ups\left( x_1\f{\pa}{\pa x_1}\right)^{2g+1}\!\!\!\!\De_{3}{\widbfH}^0_2$
$=\left( y_1^3\f{\pa}{\pa y_1}\right)^{2g+1}\!\!\!\!y_1y_2$
$=(4g+1)!!y_1^{4g+3}y_2,$ 
where the second equality follows by induction. This gives
\begin{eqnarray*}
\De_{3}\Upsilon_2&=&\SYM_{1,1}
\sum_{g\ge 1}2^{2g-1}\f{(4g+1)!!}{(2g-1)!}y_1^{4g+3}y_2t^{2g}
=\SYM_{1,1}y_1^5y_2t\sum_{g\ge 1}3{-\tfrac{5}{2}\choose 2g-1}
\left(-4y_1^2t\right)^{2g-1}\\
&=& 3\oE\,\SYM_{1,1}y_1^5y_2tB_1^{-5}.
\end{eqnarray*}
Now, consider Lemma~\ref{symFjcutm} with $m=2$, and apply $\De_{3}$ to both sides,
to obtain
$ \De_{3}\De_{2}\Phi_2=
\De_{3}\SYM_{1,1}\f{y_1^5y_2}{y_1-y_2}\f{\pa}{\pa y_1}\Phi_1(y_1,t)
+\De_{3}\Upsilon_2.$ 
The result follows from Corollary~\ref{Phione} and the expression
for $\De_{3}\Upsilon_2$ above,
since $\De_{3}$ and $\SYM_{1,1}$ commute.
\epf

Now we compare the results from the \textsf{localization} and \textsf{degeneration sides},
to obtain a proof of Faber's Conjecture for two parts.

\begin{tm}\label{Fabtwogf}
Faber's Intersection Number Conjecture is true for $n=2$.
\end{tm}

\bpf
By applying $\De_{3}\De_{2}$ to Theorem~\ref{LHS2},
and from Theorem~\ref{RHS2}, we routinely simplify to prove that
 $\De_{3}\De_{2}\Psi_2(y_1,y_2,t)$ and $\De_{3}\De_{2}\Phi_2(y_1,y_2,t)$
 (under the condition that the Faber Conjecture is true for $n=2$) are both equal to
$\oE\, \SYM_{1,1}\f{2y_1^5y_2tB_1^{-5}}{y_1-y_2}
\left(3y_1-2y_2-10y_1^3t+4y_1^2y_2t+16y_1^5t^2-8y_1^4y_2t^2\right).$
Thus we have proved that if the Faber Conjecture is true for $n=2$, then
$\De_{3}\De_{2}\Psi_2(y_1,y_2,t)=\De_{3}\De_{2}\Phi_2(y_1,y_2,t).$
Then Proposition~\ref{invDe} implies
that $\Psi_2(y_1,y_2,t)=\Phi_2(y_1,y_2,t)$, and the result follows
from Lemma~\ref{nonsing}.
\epf

\subsection{Faber's Intersection Number Conjecture for three parts}\label{Sin3p}
The proof for three parts follows the steps and strategy of the proof
for two parts, with no further technical difficulties. However, the
full details of some of the expressions are somewhat longer, and we
omit many below.

For the \textsf{localization side}, we determine that, if
Faber's Conjecture is
true for up to $3$ parts, then
\begin{eqnarray}
\Psi_3(y_1,y_2,y_3,t)&=&\tfrac{1}{2}\oE\, \SYM_{1,2}t\f{d}{dt}
t\left( 1+B_1^{-1}\right)
\f{\whY_1^8\whY_2\whY_3}{(\whY_1^2-\whY_2^2)(\whY_1^2-\whY_3^2)}\label{Psi3albe}\\
&{}&-\tfrac{1}{2}\oE\, \SYM_{1,1,1}t^{-2}\f{d}{dt}
t^4\left( 1+B_1^{-1}\right)\f{y_1^2y_2}{y_1-y_2}A_1
\f{\whY_1^6\whY_3}{\whY_1^2-\whY_3^2}\notag\\
&{}&+\tfrac{1}{2}\oE\, \SYM_{1,2}
t^{-2}\f{d}{dt}t^3\left( 1+B_1^{-1}\right)
\f{y_1^6 y_2 y_3}{(y_2-y_1)(y_3-y_1)}\whY_1^4\notag\\
&{}&+\tfrac{1}{2}\oE\, \SYM_{1,1,1}t\left( 1+B_1^{-1}\right)
\al
+\tfrac{1}{2}\oE\, \SYM_{1,2}t\left( 1+B_1^{-1}\right)
\be,\notag
\end{eqnarray}
where
\begin{eqnarray*}
\al&=&
-\f{y_3^2y_2}{y_3-y_2}A_1
\left(\f{\whY_1^4\whY_3^3}{\whY_1^2-\whY_3^2}+2\f{\whY_1^4\whY_3^5}{(\whY_1^2-\whY_3^2)^2}\right)
+2\f{y_1^2y_2}{y_1-y_2}A_1
\f{\whY_1^8\whY_3}{(\whY_1^2-\whY_3^2)^2}\\
&{}&+\f{y_1^3 y_2^4 y_3}{(y_2-y_3)(y_1-y_2)^2}A_1^2\whY_1^2\whY_2
+\f{y_1^6 y_2 y_3}{(y_2-y_1)^2(y_3-y_1)}A_1^2\whY_1^3,\\
\be&=&
(1-y_1)\whY_1^4\whY_2\whY_3
+3\f{y_1^5 y_2 y_3}{(y_2-y_1)(y_3-y_1)}A_1^2\whY_1^3.
\end{eqnarray*}
This is a large expression for $\Psi_3$, but the method of proof follows
exactly that of Theorem~\ref{LHS2}, and is routine: we apply Corollary~\ref{tSFHT4}(b)
with $m=3$, and parts (a)--(c) of Theorem~\ref{xi123}, to
obtain an expression for $\La\;\zeta_3$, and use a smoothing of Faber's
Conjecture for three parts to account for parts that are zero.  The above
expression then follows from~(\ref{Psigf}).

For the \textsf{degeneration side}, we determine that
\begin{eqnarray}
&{}&\De_{2}\Phi_3(y_1,y_2,y_3,t)=\oE\,\SYM_{1,2}y_1^5y_2y_3tB_1^{-1}
+3\oE\,\SYM_{1,2}y_1^5y_2y_3tB_1^{-5}\label{DePhi3}\\
&{}&+\tfrac{1}{2}\oE\,\SYM_{1,1,1}
\left(\f{y_1^5y_3}{y_1-y_3}\f{\pa}{\pa y_1}+\f{y_2^5y_3}{y_2-y_3}\f{\pa}{\pa y_2}\right)
t(1+B_1^{-1})
\left(\f{\whY_1^4\whY_2}{\whY_1^2-\whY_2^2}\right.
-\left. \f{y_1^2y_2}{y_1-y_2}A_1\whY_1^2\right).\notag
\end{eqnarray}
The method of proof of this result follows exactly that of Theorem~\ref{RHS2}:
we apply Lemma~\ref{symFjcutm} with $m=3$, using
the expression in Theorem~\ref{LHS2} for $\Phi_2=\Psi_2$ and
the fact that $\Ups\widbfH_3=y_1y_2y_3$, from~(\ref{singexp}), which implies that
$\Upsilon_3(y_1,y_2,y_3,t)=3\oE\,\SYM_{1,2}\, y_1^5y_2y_3tB_1^{-5},$
giving the result.

Now we compare the results from the \textsf{localization} and \textsf{degeneration sides},
to obtain a proof of Faber's Conjecture for three parts.

\begin{tm}\label{Fabthrgf}
Faber's Intersection Number Conjecture is true for $n = 3$.
\end{tm}
\bpf
Applying $\De_{2}$ to~(\ref{Psi3albe}), and from~(\ref{DePhi3}), with
routine reductions (using \textsf{Maple}), we prove that
$\De_{2}\Psi_3(y_1,y_2,y_3,t)=\De_{2}\Phi_3(y_1,y_2,y_3,t),$
under the condition that the Faber Conjecture is true for $n\le 3$.
Then Proposition~\ref{invDe} implies
that $\Psi_3(y_1,y_2,y_3,t)=\Phi_3(y_1,y_2,y_3,t)$, and the result follows
from Lemma~\ref{nonsing} and Theorem~\ref{Fabtwogf}.
\epf

\subsection{Intersection numbers with four or more parts}\label{Sin4mp}
For the case of $n$ parts,
on the \textsf{localization side}, we need to determine $\La\;\xi^{(i)}_n$, for
which we need an expression for the double Hurwitz series $\bfH^0_n$. Explicit expressions
for the latter are known for~$n=4$ \cite[Cor.~5.7]{gjv2} and $n=5$ \cite[Cor.~5.8]{gjv2}, and with the help
of \textsf{Maple}, we could in principle thus obtain expressions for $\Psi_4$ and $\Psi_5$
analogous to~(\ref{Psi3albe}).

For the \textsf{degeneration side}, we apply Lemma~\ref{symFjcutm},
and note that $\Phi_n$ is obtained recursively from $\Phi_j,\; j<n$, once
we have $\widbfH^0_n$, which is given explicitly
for all $n$ in~(\ref{singexp}).
Thus we can obtain a rational expression for $\Phi_n$ for any $n$, by
using \textsf{Maple} to help with the size of the expressions.

Together, \textit{via} Lemma~\ref{nonsing}, these would enable us to prove Faber's Conjecture
for four and five parts.
However, we are prevented from proving the result for larger
values of $n$ because we do not have explicit expressions
for the double Hurwitz series $\bfH^0_n$, and we do
not see how to obtain an explicit formula that holds 
for all $n$. In order to use our methodology to
prove the result for arbitrary $n$, we would need to
use our functional equations for the various tree series
given in Theorem~\ref{branchfg}, and then use the
Join-cut Equation for the double Hurwitz series.

\appendix
\section{Proofs of Theorems~\ref{polyfact} and~\ref{xi123}}\label{aPTF}

In this Appendix we develop the substantial amount of material that is required to
determine $\La\;\xi^{(i)}_m$ for $m=1,2,3$, and is not required elsewhere in the paper.
A key ingredient is the symmetrization of the double Hurwitz generating series. 

\subsection{Preliminaries}\label{ssPre}

\noindent \textbf{(a)} \underline{{\em Completing the symmetrization on the} \textsf{localization side}}:
In order to determine $\La\;\xi^{(i)}_m$ we have the following result, which
together with
Corollary~\ref{tSFHT4} completes the symmetrization of Theorem~\ref{branchfg} (the third form of
the Localization Tree Theorem).
It requires the symmetrized series
\begin{equation}\label{symrest2}
\left\{
\begin{array}{lll}
f_{j,m}(x_1,\ld ,x_m,u)&:=&\Xi_m f_j(z,u;\bfp)), \quad\quad
g_{j,m}(x_1,\ld ,x_m,u):=\Xi_m g_j(z,u;\bfp)),\\
\bfH^0_m(x_1,\ld ,x_m,u;\bfq)&:=&\Xi_m H^0(z,u;\bfp,\bfq),
\end{array}
\right.
\end{equation}
for $m\ge 1$, and
the substitution operator
\begin{equation}\label{eLSOm}
\Om\colon q_i\mapsto\f{i^{i-1}}{i!},\quad  i\ge 1.
\end{equation}

\begin{co}\label{remsymco}
For $m\ge 1$,
$$a)\quad \xi^{(i)}_m=\sum_{j\ge 1}\f{j^{j+i}}{j!}f_{j,m},\quad i\ge 0,$$
where, for $j\ge 1$,
$$b)\quad f_{j,m}=\sum_{\substack{k\ge 0, j_1,\ld ,j_k\geq 1, \\  i_0,\ld ,i_k\geq 1, \\  i_0+\cd +i_k=m}}
\SYM_{i_0,\ld,i_k}\f{u^{-2}}{k!}
\left(\Om\; j\f{\pa^{k+1}}{\pa q_j \pa q_{j_1}\cd \pa q_{j_k}}\bfH^0_{i_0}(x_{\rho_0},u;\bfq )\right)
\prod_{a=1}^k g_{j_a,i_a}(x_{\rho_a},u),$$
$$c)\quad g_{j,m}=\sum_{\substack{k\ge 0, j_1,\ld ,j_k\geq 1, \\  i_1,\ld ,i_k\geq 1, \\ i_1+\cd +i_k=m}}
\SYM_{i_1,\ld,i_k}
\f{(j+j_1+\cd +j_k)^{k-2}}{k!}\f{j^{j+1}}{j!}
\prod_{a=1}^k\f{j_a^{j_a}}{j_a!}f_{j_a,i_a}(x_{\rho_a},u).$$
\end{co}

\bpf
The result follows by applying $\Xi_m$ to~(\ref{Fseries}) Theorem~\ref{branchfg}(c) and~(d), for $m\ge 1$,
and using the properties
of $\Xi_m$ given in Lemmas 4.1--4.3 of~\cite{gjvai}.
\epf

We shall determine $\La\;\xi^{(i)}_m$, $m=1,2,3$, by applying Corollary~\ref{remsymco}
directly. The results that are needed for the series $\bfH^0_m$ and
its partial derivatives are given in the next two sections.

\noindent \textbf{(b)} \underline{\em The symmetrized double Hurwitz series, genus $0$}:
Now let $Q(t):=\sum_{j\ge 1}q_jt^j$, and $v$ be the unique formal power series solution of the functional 
equation
\begin{equation}\label{wxeqn}
v=xe^{uQ(v)}.
\end{equation}
Let  
\begin{equation}\label{eVVQMMI}
v_i:=v(x_i),\quad Q_i:=Q(v_i),\quad \mu(t):=(1-utQ'(t))^{-1},\quad \mu_i:=\mu(v_i).
\end{equation}
Then in~\cite{gjv2}, we proved that
\begin{equation}\label{onetwoH}
x_1\f{\pa}{\pa x_1}\bfH^0_1=uQ_1,\quad
\bfH^0_2=\log\left(\f{v_1-v_2}{x_1-x_2}\right) -uQ_1-uQ_2,\quad
\bfH^0_3=\sum_{i=1}^3 (\mu_i-1)\prod_{\substack{1\le j\le 3,\\ j\ne i}}\f{v_j}{v_i-v_j}.
\end{equation}

\noindent \textbf{(c)}  \underline{\em Partial derivatives of the symmetrized double Hurwitz series}:
Differentiating~(\ref{wxeqn}) gives
\begin{equation}\label{wxqj}
x_1\f{\pa}{\pa x_1}v_1=v_1\mu_1,\qquad\qquad
\f{\pa}{\pa q_j}v_1=uv_1^{j+1}\mu_1,
\end{equation}
and from~(\ref{onetwoH}) and~(\ref{wxqj}) we obtain
\begin{equation}\label{dH2qj}
\f{\pa}{\pa q_j}\bfH^0_2=\f{uv_1^{j+1}\mu_1-uv_2^{j+1}\mu_2}{v_1-v_2}
-uv_1^j\mu_1-uv_2^j\mu_2=\SYM_{1,1}\f{uv_1^j\mu_1v_2}{v_1-v_2}.
\end{equation}
For partial derivatives of $\bfH_1^0$ in general, we have the following result.
\begin{lm}\label{partH1}
\noindent (a) For $j\ge 1$,
$\f{\pa}{\pa q_j}\bfH_1^0 = \f{u}{j}v_1^j.$
\newline
\noindent (b) For $k\geq 1$, $j_1,\ld ,j_k\ge 1$,
$ \f{\pa^{k}}{\pa q_{j_k} \cd \pa q_{j_1}}\bfH^0_1=
u^k\left(\mu_1v_1\f{\pa}{\pa v_1}\right)^{k-2}
\left( \mu_1v_1^{j_1+\cd +j_k}\right).$
\end{lm}

\bpf
Differentiating~(\ref{onetwoH}), we obtain
\begin{equation*}
\f{\pa}{\pa q_j} x_1\f{\pa}{\pa x_1}\bfH^0_1
= uQ'(v_1)\f{\pa v_1}{\pa q_j}+uv_1^j
= uQ'(v_1)uv_1^{j+1}\mu_1+uv_1^j
=uv_1^j\mu_1,
\end{equation*}
where we have used~(\ref{wxqj}) for the second equality.
Dividing by $x_1$ and integrating, 
we obtain
$\f{\pa}{\pa q_j}\bfH_1^0=\int_0^{x_1} uv_1^j\mu_1\f{dx_1}{x_1}
= \int_0^{v_1} uv_1^{j-1}dv_1,$
where we have used~(\ref{wxqj}) to change variables for
the second equality. Part~(a) follows immediately.

We prove part~(b) by induction on $k$. For the base case, note that
the result is true for $k=1$, since this is equivalent to part~(a) above.
Now assume that the result holds for $k-1\ge 1$, and note that~(\ref{wxqj})
implies the operator identity $\mu_1v_1\pa/\pa v_1=x_1\pa/\pa x_1$, so
we have
\begin{eqnarray*}
\f{\pa^{k}}{\pa q_{j_k} \cd \pa q_{j_1}}\bfH^0_1
&=&\f{\pa}{\pa q_{j_k}}u^{k-1}\left(x_1\f{\pa}{\pa x_1}\right)^{k-2}
\left(\mu_1v_1\f{\pa}{\pa v_1}\right)^{-1}
\left( \mu_1v_1^{j_1+\cd +j_{k-1}}\right)\\
&=&u^{k-1}\left(x_1\f{\pa}{\pa x_1}\right)^{k-2}
\f{\pa}{\pa q_{j_k}}
\left( \f{v_1^{j_1+\cd +j_{k-1}}}{j_1+\cd +j_{k-1}}\right), 
\end{eqnarray*}
and part (b) follows from~(\ref{wxqj}) and the operator identity.
\epf

\noindent \textbf{(d)}  \underline{\em Polynomiality and the symmetrized double Hurwitz series}:
{}From Corollary~\ref{remsymco}(a), (b), to evaluate $\La\;\xi^{(i)}_m$ we
need to apply the substitution $\La\;\Om$ to partial derivatives of the symmetrized double 
Hurwitz series $\bfH^0_m$, which we have shown above can be expressed
in terms of $v_i$ and $\mu_i$. Thus define 
\begin{equation}\label{Vvdef}
V_i:=\La\;\Om\; v_i.
\end{equation}
In the following result, the action of the change of variables $\sC$, the
second step in our fundamental transformation, is given for expressions in $V_i$.
Note how indirectly we proceed in the proof, especially for part~(a), as referred to in the
discussion following~(\ref{Ty}).
This result is the key to polynomiality on the \textsf{localization side}.
\begin{lm}\label{Vchange}
For $i\ge 1$,
$$a)\quad \sC\;y(V_i)=\f{(1-u)y_i}{1-uy_i},\quad
b)\quad \sC\;\La\;\Om\;\mu_i =\f{1-uy_i}{1-u},\quad
c)\quad \sC\; V_i\f{\partial}{\partial V_i}=y_i^2(y_i-1)\f{1-u}{1-uy_i}\f{\partial}{\partial y_i}\sC.$$
\end{lm}

\bpf 
For part~(a), the series $w=w(x)$ defined in~(\ref{Ty}a) is the unique formal power series
solution to $w=xe^w$, which we shall call the \emph{tree equation} 
(see, e.g.,~\cite[Ex.~1.2.5]{gjce}).
Now $V_i=V(x_i)$, where, from~(\ref{wxeqn}), $V=V(x)$ is the unique formal power series
solution to
\begin{equation}\label{Vwx}
V=(1-u)^{-1}x\exp\left(-u\sum_{j\ge 1}\f{j^{j-1}}{j!}V^j\right)=
(1-u)^{-1}xe^{-uw(V)},
\end{equation}
where the second equality follows from~(\ref{Ty}a). But the tree equation
gives $w(V)=Ve^{w(V)}$, and eliminating $V$ between this equation
and~(\ref{Vwx}), we obtain $(1-u)w(V)=xe^{(1-u)w(V)}$. Comparing
this to the tree equation, we conclude that $(1-u)w(V(x))=w(x)$, and
by~(\ref{Ty}b), this implies
\begin{equation}\label{yVvx}
y(V_i)=\f{(1-u)y_i}{1-uy_i},
\end{equation}
where here $y_i=y(x_i)$ is a series in $x_i$. Part~(a) follows immediately.

For part~(b), from~(\ref{eVVQMMI}) and~(\ref{Ty}b) we have
$$\La\;\Om\;\mu_i=\f{1}{1+u\sum_{j\ge 1}\f{j^j}{j!}V_i^j}
=\f{1}{1+u(y(V_i)-1)},$$
and part~(b) follows immediately from part~(a).

For part~(c), from~(\ref{Upident}a), we have
\begin{equation}\label{VdyV}
V_i\f{\partial}{\partial V_i}=
y(V_i)^2(y(V_i)-1)\f{\partial}{\partial y(V_i)}.
\end{equation}
But, by differentiating~(\ref{yVvx}), we obtain
$$\f{\pa y(V_i)}{\pa y_i}= \f{1-u}{(1-uy_i)^2}
=\f{y(V_i)^2}{(1-u)y_i^2},$$
which implies $\sC\; y(V_i)^2\f{\pa}{\pa y(V_i)}=(1-u)y_i^2\f{\pa}{\pa y_i}\sC$, and
part~(c) follows immediately from~(\ref{VdyV}). 
\epf

\noindent \textbf{(e)}  \underline{\em Two results}:
In the proofs of Section~\ref{App2}, we make use of the identity 
\begin{equation}\label{ideninv}
\sum_{j,k\ge 1}\f{1}{j+k}\f{j^{j+1}}{j!}V_1^j\f{k^k}{k!}V_2^k
=\f{V_2}{V_2-V_1}-y(V_1)^2\f{y(V_2)-1}{y(V_2)-y(V_1)},
\end{equation}
whose proof follows from~\cite{gj0}, and the fact (an easy induction) that
\begin{equation}\label{diffyiratio}
\left(\f{y_i^3}{1-uy_i}\f{\pa}{\pa y_i}\right)^i\f{y_i}{1-uy_i}=
(2i-1)!!\left(\f{y_i}{1-uy_i}\right)^{2i+1}, \qquad i\ge 1.
\end{equation}

\subsection{The proofs}\label{App2}
\noindent\underline{\em Proofs of Theorem~\ref{polyfact}(a) ($m=1$), and Theorem~\ref{xi123}(a)}:
{}From Corollary~\ref{remsymco}(b) with $m=1$, Lemma~\ref{partH1}(a) and~(\ref{Vvdef}), we have
\begin{equation}\label{fj1}
\La\, f_{j,1}= u^{-2}\La\,\Om\; j\f{\pa}{\pa q_j}\bfH^0_1
=-u^{-1}\La\,\Om\; v_1^j=-u^{-1}V_1^j.
\end{equation}
Then from Corollary~\ref{remsymco}(a) with $m=1$, we obtain
\begin{equation*}
\La\;\xi^{(i)}_1=\sum_{j\ge 1}\f{j^{j+i}}{j!}\La\, f_{j,1}
=-u^{-1}\sum_{j\ge 1}\f{j^{j+i}}{j!}V_1^j
=-u^{-1}\left(V_1\f{\partial}{\partial V_1}\right)^i 
\left(y(V_1)-1\right),
\end{equation*}
from~(\ref{Ty}b).
But Lemma~\ref{Vchange}(a) and~(c) give
\begin{equation*}
\sC\;\La\;\xi_1^{(i)}
=-u^{-1}(1-u)^{i}\left(y_1^2\f{y_1-1}{1-uy_1}\f{\pa}{\pa y_1}\right)^{i}
\left(\f{(1-u)y_1}{1-uy_1}-1\right),
\end{equation*}
and Theorem~\ref{polyfact}(a) ($m=1$) follows immediately.
Then
$\Ups\;\La\;\xi_1^{(i)}$ $= \Ups\;u^{-1}(1-u)^{i+1}
\left(\f{y_1^3}{1-uy_1}\f{\pa}{\pa y_1}\right)^{i}$
$\f{y_1}{1-uy_1}$,
and Theorem~\ref{xi123}(a) follows immediately from~(\ref{diffyiratio}).
\epf

\noindent\underline{\em Proofs of Theorem~\ref{polyfact}(a) ($m=2$), and Theorem~\ref{xi123}(b)}:
{}From Corollary~\ref{remsymco}(c) with $m=1$, and~(\ref{fj1}), we obtain
\begin{equation}\label{gj1}
\La\; g_{j_1,1}(x_1,u)=\sum_{j_2\ge 1}\f{1}{j_1+j_2}\f{j_1^{j_1+1}}{j_1!}\f{j_2^{j_2}}{j_2!}
\La\; f_{j_2,1}(x_1,u)
=-u^{-1}\f{j_1^{j_1+1}}{j_1!}\sum_{j_2\ge 1}\f{1}{j_1+j_2}\f{j_2^{j_2}}{j_2!}V_1^{j_2},
\end{equation}
and from Corollary~\ref{remsymco}(b) with $m=2$, we obtain
$\La\, f_{j,2}=S_1+S_2,$
where
\begin{equation*}
S_1=u^{-2}\La\,\Om\; j\f{\pa}{\pa q_j}\bfH^0_2,\;\;\;\;\;\;
S_2=u^{-2}\SYM_{1,1}\sum_{j_1\ge 1}\La\; g_{j_1,1}(x_2,u)
\;\Om\; j\f{\pa^2}{\pa q_j \pa q_{j_1}}\bfH^0_1(x_1,u).
\end{equation*}
We can determine $S_1$ immediately
from~(\ref{dH2qj}), giving
$S_1=-u^{-1}j\, \SYM_{1,1}\f{V_1^j V_2}{V_1-V_2}\La\,\Om\,\mu_1.$
To determine $S_2$, use Lemma~\ref{partH1} with $k=2$ and~(\ref{gj1}),
to obtain
$ S_2=-u^{-1}j\,\SYM_{1,1}V_1^j\left(\La\,\Om\;\mu_1\right)$
$\sum_{j_1,j_2\ge 1}\f{1}{j_1+j_2}\f{j_1^{j_1+1}}{j_1!}V_1^{j_1}$
$\f{j_2^{j_2}}{j_2!}V_2^{j_2}$.
Adding $S_1$ and $S_2$, and using~(\ref{ideninv}), we get
$f_{j,2}=u^{-1}j\,\SYM_{1,1}V_1^jy(V_1)^2$ $\f{y(V_2)-1}{y(V_2)-y(V_1)}$
$\La\,\Om\;\mu_1.$
Then from Corollary~\ref{remsymco}(a) with $m=2$, we obtain
\begin{eqnarray*}
\La\;\xi^{(i)}_2&=&\sum_{j\ge 1}\f{j^{j+i}}{j!}\La f_{j,2}
=u^{-1}\SYM_{1,1}y(V_1)^2\f{y(V_2)-1}{y(V_2)-y(V_1)}\left(\La\,\Om\;\mu_1\right)
\sum_{j\ge 1}\f{j^{j+i+1}}{j!}V_1^j\\
&=&u^{-1}\SYM_{1,1}y(V_1)^2\f{y(V_2)-1}{y(V_2)-y(V_1)}\left(\La\,\Om\;\mu_1\right)
\left(V_1\f{\pa}{\pa V_1}\right)^{i+1}
\left(y(V_1)-1\right),
\end{eqnarray*}
from~(\ref{Ty}b). Now Lemma~\ref{Vchange}(a),(b),(c) give
\begin{equation*}\label{Laxi2}
\sC\;\La\;\xi_2^{(i)}
=u^{-1}(1-u)^{i+1}\SYM_{1,1}
\f{y_1^2(y_2-1)}{y_2-y_1}
\left(y_1^2\f{y_1-1}{1-uy_1}\f{\pa}{\pa y_1}\right)^{i+1}
\left(\f{(1-u)y_1}{1-uy_1}-1\right) ,
\end{equation*}
and Theorem~\ref{polyfact}(a) ($m=2$) follows immediately.
Then we have
\begin{equation*}
\Ups\;\La\;\xi_2^{(i)}=
\Ups\;u^{-1}(1-u)^{i+2}\SYM_{1,1}\f{y_1^2y_2}{y_2-y_1}
\left(\f{y_1^3}{1-uy_1}\f{\pa}{\pa y_1}\right)^{i+1}
\f{y_1}{1-uy_1},
\end{equation*}
and Theorem~\ref{xi123}(b) follows immediately from~(\ref{diffyiratio}).
\epf

\noindent\underline{\em Outline of proofs of Theorem~\ref{polyfact}(a) ($m=3$), and Theorem~\ref{xi123}(c)}:
The method of proof of these results follows those for the cases $m=1,2$ above,
with no further technical difficulties. We give only an outline.
We first determine $\La\, f_{j,3}$ from Corollary~\ref{remsymco}(b)
This requires $\La\; g_{j,2}$, which we determine from Corollary~\ref{remsymco}(c), and
the partial derivatives $\f{\pa^3}{\pa q_j \pa q_{j_1}\pa q_{j_2}}$ $\bfH^0_1$,
 $\f{\pa^2}{\pa q_j \pa q_{j_1}}\bfH^0_2$ and
 $\f{\pa}{\pa q_j}\bfH^0_3$. But rational expressions for these
partials are routinely obtained, from Lemma~\ref{partH1}(b) for $\bfH^0_1$,
and by applying the chain rule and~(\ref{wxqj})
to~(\ref{dH2qj}) and~(\ref{onetwoH}). 
We then use Corollary~\ref{remsymco}(a) with $m=3$ to obtain an explicit
expression for $\La\;\xi^{(i)}_3$, and then apply $\sC$ to prove
Theorem~\ref{polyfact}(a) ($m=3$).  Then Theorem~\ref{xi123}(c)
follows immediately.
\epf

\vfil
\newpage
\section{Glossary of Notation}
\label{glossary}\lremind{glossary}If a symbol in column~1 has an argument, it appears as the first item of the entry in column~3.
In column~2,``C'' denotes ``Corollary'',  ``Cn'' denotes ``Conjecture'',  ``D'' denotes ``Definition'', 
``P'' denotes ``Proposition''  ``T'' denotes ``Theorem'',  ``b'' denotes ``immediately before'',
``f'' denotes ``immediately following'' and ``I'' denotes ``Introduction to''.

\vspace{0.25in}
\begin{minipage}{3.20in}  
\tracingtabularx
\begin{tabularx}{\linewidth}%
{  >{\setlength{\hsize}{0.50\hsize}}X%
   >{\setlength{\hsize}{0.35\hsize}}X%
   >{\setlength{\hsize}{2.15\hsize}}X|
}
${}^\ddagger$   &f(\ref{Aonly})  &  Removal of  contribution of root from  corresponding product \\
$\bullet$            &  D\ref{Dloctr} & Root-vertex of tree $\rt$ \\
$\star$                & D\ref{starprod} & Combinatorial product   \\
$[x^k]f$                & -                     & Coefficient of $x^k$ in a formal power series $f$ \\
$[.,.]_{\eta}$ & D\ref{expgf} & $\sA ,\wt:$ Exponential generating series \\
$|\al|$                  &  I\S\ref{DandL}     & Sum of parts of  $\al$ \\
$\laFab  \cdots \raFab$  &(\ref{eWS})   &  $\tau_{a_1}\cdots\tau_{a_n}\la_k:$ Faber symbol for 
                                                         $\pi_* \left( \psi_1^{a_1} \cdots \psi_n^{a_n} \la_k \right)$ \\[5pt]
$\langle \cdots \rangle_g^{\Faber}$ & (\ref{Fabrational}) & $\tau_{a_1}\cdots\tau_{a_n}\la_k:$ Faber number\\
$\al$                    & I\S\ref{DandL}     &Ptn. of $d$:  branching over $\infty$ \\
$\al$                    &(\ref{spal})       & $\coprod_{v\in\sV_\infty}\gamma^v$ \\    
$\al(\rt)$   & bD\ref{starprod} & Relabelled localization tree \\
$\be$                   & \S\ref{sssRTPT}    & Ptn. of $d$:  branching  over $0$ \\
$\be ^v$              & D\ref{Dloctr} & List of weights on $0\infty$-edges incident with  $\infty$-vertex $v$ \\
$\ga ^v$              &D\ref{Dloctr}  &  List of weights on all $\infty t$-edges incident with $\infty$-vertex $v$ \\
$\De_{k}$     &  (\ref{eDOp})  & Partial differential operator \\   
$\de^v$                 &  D\ref{Dloctr} & List of weights on $0\infty$-edges incident with $0$-vertex $v$ \\
$\epsilon(e)$      &   D\ref{Dloctr}    &       Weight on a $0\infty$-edge of $\rt$ \\
$\eta_0(\rt)$            & D\ref{Dloctr} & No.\ of  non-root $0$-vertices in $\rt$ \\
$\eta_k(\rt)$            & bD\ref{expgf} & No.\ of $t$-vertices in $\rt$ incident with edge of weight $k$ \\
$\Upsilon_m$  &  (\ref{ser2g})  &    $(y_1,\ldots, y_m)$ \\  
$\ka_i$           & \S\ref{ssTR}  & ``$\ka$-class'' \\
$\La$              & (\ref{eLSO}) & Substitution operator  $u\mapsto -u,$ $x_i\mapsto(1-u)^{-1} x_i$ \\
$\la_i$           &  \S\ref{ssTR} & $\la$-class \\
$\mu,\mu_i$             &(\ref{eVVQMMI}) & Used for double Hurwitz series  \\
$\Xi_m$      & (\ref{eXi})  & Symmetrization operator \\
$\xi^{(i)}$   & (\ref{Fseries}) & $(z,u;\bfp ):$ Localization tree series \\
$\xi^{(i)}_m$ &  (\ref{symrest1}) & $(x_1,\ldots,x_m):$ Symmetrized localization tree series \\ 
$\Phi_m$  & (\ref{Phigf})  & $(y_1,\ldots,y_m,t):$ Generating series for     $\Ups\,F_m^g(x_1,\ld ,x_m)$  \\ 
\end{tabularx}
\end{minipage}
\begin{minipage}{3.2in}  
\tracingtabularx
\begin{tabularx}{\linewidth}%
{  >{\setlength{\hsize}{0.5\hsize}}X%
   >{\setlength{\hsize}{0.35\hsize}}X%
   >{\setlength{\hsize}{2.15\hsize}}X
}
%
$\rho_{j}$       &  (\ref{eRHOJ})           &  Used in $\Xi_m$   \\
$\Psi_m$  & (\ref{Psigf})  & $(y_1,\ldots,y_m,t):$ Generating series for
                                      $\Ups\,\La\,\zeta_m^g/\G_{g,1}$ \\
$\psi_i$              & \S \ref{ssTR1}  &  First Chern class \\
$\Omega$  &  (\ref{eLSOm})  &   Substitution operator $q_i\mapsto\f{i^{i-1}}{i!}$  \\
$\ombr,\ombr_0$ & \S\ref{branchdec} & Variants of branch decompositions for trees   \\[-13pt]     
$\ombr_{\infty}$ & & \\
$A_i$    &  (\ref{AiBi})  &  -  \\
$A(\rt )$ &(\ref{Aonly}) & Weight function for $\rt$\\
$B(\rt )$ &(\ref{ABCD}) & Weight function for $\rt$\\   
$B_i$    &  (\ref{AiBi})  & - \\          
$ br$                   &(\ref{branchmorphism})  & Fantechi-Pandharipande branch morphism \\
$C(\rt )$ &(\ref{ABCD}) & Weight function for $\rt$\\
$\sC$               &  f(\ref{Ty}) & Change of variables  $x_i\leadsto y_i$ \\
$D(\rt )$ &(\ref{ABCD}) & Weight function for $\rt$\\
$d$                      &  \S\ref{sssROPT}   &  Degree of a cover \\
$\oE$       & b(\ref{AiBi})  &  Even subseries operator \\   
$\sE_{0\infty}(\rt)$  & D\ref{Dloctr} & Set of edges  $0\infty$-edges of $\rt$ \\ 
$\sE_{\infty t}(\rt)$  & D\ref{Dloctr} & Set of edges  $\infty t$-edges of  $\rt$ \\                            
$\F^{g,\al}$         &(\ref{fhclass})  & Faber-Hurwitz class \\
$F^g_{\al}$         &(\ref{ianeq})   &  Faber-Hurwitz number,  $\in\Q;$ note $\F^{g,\al}= F^g_{\al} \G_{g,1}$\\
$F^g$      & (\ref{gsFnum}) & $(z;\bfp)$: Generating series for Faber numbers $F^g_{\al}$ \\
$F^g_m$ &  (\ref{symrest1}) & $(x_1,\ldots,x_m)$ Symmetrized $F^g(z;\bfp)$ \\
$f_j$        & T\ref{branchfg}     & Generating series for $\sT_{0,j}$ \\
$f_{j,m}$  &  (\ref{symrest2})  & $(x_1,\ld ,x_m,u)$:  Symmetrized series for $\sT_{0,j}$ \\
$\G_{g,1}$           &  T\ref{tGgd}   &  A ``natural generator'' of $R_{2g-1}(\cm_{g,1}),$
                                  C\ref{cgformula} \\
$\G_{g,d}$          &T\ref{tGgd} & - \\
$\G_{g,d, \sim}$ & T\ref{tGgd} & - \\
$g$    &   \S\ref{ssTR1}  &  Genus (of a curve) \\
$g_j$        & T\ref{branchfg}     & Generating series for $\sT_{\infty,j}$ \\
$g_{j,m}$  &  (\ref{symrest2})  & $(x_1,\ld ,x_m,u)$ symmetrized series for $\sT_{\infty,j}$ \\
$H^0_{\al}$         &(\ref{singHur})  & Genus $0$ (single) Hurwitz no.\ \\[5pt]
$H^0_{\al,\be}$   & (\ref{eHalbe})  & Genus $0$ double Hurwitz no.\ \\[5pt]
\end{tabularx}
\end{minipage}

\begin{minipage}{3.20in}  
\tracingtabularx
\begin{tabularx}{\linewidth}%
{  >{\setlength{\hsize}{0.85\hsize}}X%
   >{\setlength{\hsize}{0.35\hsize}}X%
   >{\setlength{\hsize}{1.80\hsize}}X|
}
$\H^{g,\al}_j $      & \S\ref{hcdef}      & Hurwitz class  \\[5pt]
$\widH^0$ & (\ref{sHser})    & $(z;\bfp)$: Gen. series for genus $0$  single Hurwitz nos. $H^0_{\al}$ \\
$H^0$ &  (\ref{dHser})    & $(z,u;\bfp;\bfq)$: Gen. series for genus $0$ double Hurwitz 
                                                                       nos. $H^0_{\al,\be}$ \\
$ \widbfH^0_m$  &(\ref{symrest1})  &  $(x_1,\ld,x_m)$: Symmetrized single Hurwitz series \\  
$\bfH^0_m$ &(\ref{symrest2})  & $(x_1,\ld ,x_m,u;\bfq)$: Symmetrized double Hurwitz series, genus $0$ \\
$L$ & P\ref{bigbusiness} & Class in $\Sym^{r^g_{\al ,\be}}(\proj^1)$ \\
$\l(\al)$               &  I\S\ref{DandL}         & Number of parts of $\al$ \\  
$\cm_g$             & \S \ref{ssTR1} & Moduli space of smooth curves \\
$\cmbar_{g,n}$ & \S \ref{ssTR1} & Moduli space  of stable $b$-pointed genus $g$ curves\\
$\cmbar_{g,\al}(\proj^1)$  & \S\ref{sssROPT}  & Moduli space of genus $g$ relative stable maps to
                              $\proj^1$ relative to one pt. $\infty$ \\ 
$\cmbar_{g, \al, \be}(\proj^1)$ & \S\ref{sssRTPT} &  Moduli space of relative stable maps
to $\proj^1$ relative to two points $0$ and $\infty$ \\
$\cmbar_{g, \al, \be}(\proj^1)_{\sim}$ & \S\ref{sssRTPT} & Moduli space of ``rubber maps''\\
$\cm_{g,\al}(\proj^1)^{rt}$ & \S\ref{sssSRMRT} & Moduli space of relative stable  maps with
                               rational tails  \\
$\cm_{g,\al,\be}(\proj^1)^{rt}$ & \S\ref{sssSRMRT} & Similar to the  above \\ 
$\cm_{g,\al}(\proj^1)_{\bullet}$ & \S\ref{secpossdis} & Possibly disconnected moduli space \\
$\P^g_m(\dots)$  &(\ref{Fabpoly2}) & $(\al_1,\ldots,\al_m)$: Faber polynomial \\ 
$\sP$                  &   I\S\ref{DandL}          &  Set of all non-empty partitions \\
$\bfp$                 &(\ref{dHser})    & $(p_1,p_2,\ld )$ \\
$p_i$                  & f(\ref{sHser}) & Indeterminate marking a part $i$ in $\al$ \\
$Q(t)$  & (\ref{eVVQMMI})  & Generating series for indeterminates $q_i$ \\
$\bfq$                  &(\ref{dHser})    & $(q_1,q_2,\ld )$ \\
$q_i$                  &  f(\ref{dHser}) & Indeterminate marking a part $i$ in $\be$ \\
\end{tabularx}
\end{minipage}
\begin{minipage}{3.20in} 
\tracingtabularx
\begin{tabularx}{\linewidth}%
{  >{\setlength{\hsize}{0.70\hsize}}X%
   >{\setlength{\hsize}{0.55\hsize}}X%
   >{\setlength{\hsize}{1.75\hsize}}X
}
$R^*(\cm^{rt}_{g,n})$ & I\S\ref{sor} & Tautological ring of $\cm^{rt}_{g,n}$ \\
$r^g_{\al}$         & (\ref{rgadef})  &  ``Expected'' no.\ of  branch pts.\ away from $\infty$ \\
$r^g_{\al,\be}$   & (\ref{rgabdef}) &  ``Expected'' no.\ of  branch pts.\ 
                                                                 away from $0$ and $\infty$ \\ 
$r^{\Faber}_{g,\al} $     & (\ref{rFgadef})  &  Number of fixed branch pts. on $\proj^1$ \\ 
$r_{\infty}$          &(\ref{rrr})    &  -  \\
$R_j(\cm)$         &S\ref{ssTR} &  $R^{\dim \cm - j}(\cm)$ for open subset  $\cm$ of $\cmbar_{g,n}$ \\
$\sS_i$ & fT\ref{polyfact} & Operator for terms of total degree $i$ \\
$\SYM_{i_1,.,i_k}$ &f(\ref{eRHOJ}) &    Summation operator for symmetrization \\
$\sT_{g,m},\sT_{0 ,j}$        & D\ref{Dloctr} & Sets of localization trees \\
$\sT_{\infty ,j}$  & & \\
$\Ups$, $\Ups'$      & \S\ref{Spft}  & Operators for maximum degree  \\
$\Ups\;\Xi_m$ & \S\ref{sssSOD}  &  Fundamental transformation \\
$\zeta^g$        &T\ref{branchfg}(b)     & $(z,u;\bfp)$: Tree sum related to generating series for $\F_{g,\al}$  \\
$\zeta^g_m$ & (\ref{symrest1}) & $(x_1,\ld ,x_m,u)$: Symmetrized $\zeta^g(z,u;\bfp)$ \\
$\rt$                      & D\ref{Dloctr} &  Localization tree \\
$\sV_0(\rt)$       &  D\ref{Dloctr}& Set of $0$-vertices in $\rt$ \\
$\sV_\infty(\rt)$ &  D\ref{Dloctr}& Set of  $\infty$-vertices in $\rt$ \\
$\sV_t(\rt)$       &  D\ref{Dloctr}& Set of $t$-vertices in $\rt$ \\
$\cU_k(\cS)$ & fL\ref{prodlem} & Unordered list \\
$v,v_i$    &  (\ref{wxeqn}),(\ref{eVVQMMI}) & Implicitly defined indets. for double Hurwitz series\\
$\rv_i$ & \S\ref{branchdec}& $0$-vertex labelled $i$ \\
$V_i$  & (\ref{Vvdef})  & Transform of $v_i$\\
$w, w_i$  & (\ref{Ty}a) &  Rooted tree series \\
$\rw_{\ga}$ & \S\ref{branchdec}& Star of $t$-vertices rooted at an $\infty$-vertex \\
$x_i$  & (\ref{eXi}) & Indeterminates for symmetrized series \\
$\whY_m$     &             (\ref{whYA})  &     -  \\
$Y_i(u)$                   & (\ref{Yiu}) &    -  \\
$y, y_i$  & (\ref{Ty}b) & Implicitly defined indets. for symmetrized series  \\
%
%
\end{tabularx}
\vspace{0.0in}
\end{minipage}

%

 \end{document}